\documentclass[a4paper,twoside,10pt]{amsart}

\usepackage[utf8x]{inputenc}

\usepackage{amssymb}
\usepackage{amsmath}
\usepackage{epsfig}

\usepackage[T1]{fontenc}

\newtheorem{The}{Theorem}[section]
\newtheorem{Lem}[The]{Lemma}
\newtheorem{Cor}[The]{Corollary}
\newtheorem{Pro}[The]{Proposition}

\newtheorem{Rem}[The]{Remark}
\newtheorem{Defi}[The]{Definition}

\newtheorem{Exam}[The]{Example}

\def\A{\mathbb A}

\def\C{\mathbb C}
\def\D{\Delta}
\def\d{\delta}

\def\g{\gamma}

\def\l{\lambda}

\def\N{\mathbb N}

\def\Q{\mathbb Q}
\def\R{\mathbb R}
\def\r{\rho}

\def\s{\sigma}
\def\t{\tau}

\def\Z{\mathbb Z}
\def\z{\zeta}

\begin{document}

\date{\today}

\title{Jet schemes of quasi-ordinary surface singularities}

\author{Helena Cobo}


 \email{helenacobo@gmail.com}

\author{Hussein Mourtada}

\address{Institut de Mathématiques de Jussieu-Paris Rive Gauche, Université Paris 7, B\^atiment Sophie Germain, 75013 Paris, France.}

\email{hussein.mourtada@imj-prg.fr}

\begin{abstract}
In this paper we give a complete description of the irreducible components of the jet schemes (with origin in the singular locus) of a  two-dimensional quasi-ordinary hypersurface singularity.
We associate with these components and with their codimensions and embedding dimensions, a weighted graph. We prove that the data of this weighted graph is equivalent to the data of the
topological type of the singularity. We also determine a component of the jet schemes (or equivalently, a divisor on $\A^3$), that computes the log canonical threshold
 of the singularity embedded in $\A^3$. This provides us with pairs $X\subset\A^3$ whose log canonical thresholds are not contributed by monomial divisorial valuations.
 Note that for a pair $C\subset\A^2$, where $C$ is a plane curve, the log canonical threshold is always contributed by a monomial divisorial valuation (in suitable coordinates of $\A^2$).
\end{abstract}



\maketitle

\section{Introduction}

A quasi-ordinary singularity $(X,0)$ of dimension $d,$ comes with a finite projection $p:X\longrightarrow \mathbb{A}^d,$ whose discriminant is a normal crossing
divisor. These singularities appear in the Jungian approach to resolution of singularities (see \cite{PP-HJ}). We are interested in irreducible quasi-ordinary
hypersurfaces. Thanks to Abhyankar-Jung theorem, we know that a hypersurface of this type can be parametrized by a Puiseux series (i.e an element in
$\C[[x_1^{\frac{1}{n}},\ldots,x_d^{\frac{1}{n}}]]).$ Moreover,  some special exponents (called the characteristic exponents) which belong to the support of
this series, are complete invariants of the topological type of the singularity (see \cite{Gau}). In particular, they determine invariants which come from
resolution of singularities, like the log canonical threshold or the Motivic zeta functions (\cite{BGG}, \cite{ACLM}, \cite{CoGPqo}). They also give
insights about the construction of a resolution of singularities (\cite{BMc0}, \cite{BMc}, \cite{GP4}, \cite{Vill}).

\vspace{3mm}

Our aim is to construct some comparable complete invariants for all type of singularities. We search for such invariants in the jet schemes. For $m\in \mathbb{N},$ the
$m$-jet scheme, denoted by $X_m$, is a scheme that parametrizes morphisms $\mbox{Spec }\C[t]/(t^{m+1})\longrightarrow X.$ Intuitively we can think of it as
parametrizing arcs in an ambient space, which have a large contact, depending on $m$, with $X.$ We know already that some invariants which come from resolution of singularities
are encoded in jet schemes (\cite{Mus}, \cite{EinMustata}).

\vspace{3mm}

We want to extract from the jet schemes information about the singularity, which can be expressed in terms of invariants of  resolutions of singularities.
With this, our next goal is to construct a resolution of singularities by using invariants of jet schemes. For specific types of singularities, the knowledge of the irreducible components of the jet schemes $X_m$
of a singular variety X, together with some invariants of them, like dimension or embedding dimension, permits to determine deep invariants of the singularity of $X$: the topological type in the case of curves (see \cite{Hcur}),
and the analytical type in the case of normal toric surfaces (see \cite{HCRAS} and \cite{Htor}). Moreover, in the case of irreducible plane curves, the minimal embedded resolution can be constructed from the jet schemes (\cite{LMR}), and the same for rational double point singularities (\cite{MP}).

\vspace{3mm}

Notice that, understanding the structure of jet schemes for particular singularities, remains a difficult problem. These structures have been studied in \cite{Yuen}
and \cite{DoC} for determinantal varieties, in \cite{Hcur}
for plane curve singularities, \cite{HCRAS} and \cite{Htor} for normal toric surfaces, in \cite{Hrat} for rational double point surface singularities, and in
\cite{SS} for commuting matrix pairs schemes. In the toric case, no result is known for dimension bigger than two.

More general jet schemes can be defined, and there are relations with invariants of singularities (see \cite{Mus14}).

\vspace{3mm}

In this paper, we study jet schemes of a two-dimensional quasi-ordinary hypersurface singularity $X$. We will give a combinatorial description of the irreducible components of the set of $m$-jets with center
in the singular locus of $X$, in terms of the following invariants of the singularity: the lattices $N_0,N_1,\ldots,N_g,$ the minimal system of generators $\gamma_1,\ldots,\gamma_g$ of the semigroup $\Gamma$ of $X,$
and the numerical data attached to them, $n_1,\ldots,n_g$ and $e_1,\ldots,e_g$ (see Section \ref{secQO} for definitions). Given $h\in\C[x_1,\ldots,x_n]$, an algebraic variety $X$, and $p,m\in\Z_{>0}$ with $p\leq m$,
let $\mbox{Cont}_X^p(h)_m$ be the locally closed set  defined as $\mbox{Cont}_X^p(h)_m=\{\gamma\in X_m\ |\ \mbox{ord}_t(h\circ\gamma)=p\}$. Then we associate with any lattice point $\nu\in N_0$ the constructible set
\[D_m^\nu=(\mbox{Cont}_X^{\nu_1}(x_1)_m\cap\mbox{Cont}_X^{\nu_2}(x_2)_m)_{red},\]
and its closure $C_m^\nu=\overline{D_m^\nu}$. The sets $C_m^\nu$ are the candidates to be the irreducible components, but there are many inclusions among them. We study these inclusions by defining, on the lattices
associated with $X$, subtle relations which depend on the singular loci of the quasi-ordinary surfaces defined by the approximated roots.
For some components the relation is very easy, is given by the product ordering $\leq_p$, and we can prove that
\[\mbox{ if }\nu\leq_p\nu'\mbox{ then }C_m^{\nu'}\subset C_m^\nu.\]
For other components the relation is more complicated, we have that
\[\mbox{ if }\nu'-\nu\in\s_{Reg,j'(m,\nu)}\mbox{ then }C_m^{\nu'}\subset C_m^\nu\]
where $j'(m,\nu)$ is the integer $j\in\{0,\ldots,g\}$ defined by
\[n_je_j\langle\nu,\gamma_j\rangle+e_j\leq m<n_{j+1}e_{j+1}\langle\nu,\gamma_{j+1}\rangle+e_{j+1},\]
and $\s_{Reg,j'(m,\nu)}$ is defined to keep track of the singular locus of the $j'(m,\nu)$-th approximated root as follows.
Since the normalization of any quasi-ordinary singularity is a toric variety, let $\nu_j: Z_{\s,N_j}\longrightarrow V(f_j)$ be the normalization of $V(f_j)$ ($f_j$ being the $j$-th approximated root).  Then $Z_{\s,N_j}\setminus\nu_j^{-1}(Sing(V(f_j)))$ is an open set of the toric variety $Z_{\s,N_j}$, it is a union of orbits, and we denote by $\s_{Reg,j}$ the fan formed by the faces of $\s$ corresponding to the orbits.
Then, with the minimal elements with respect to these relations we define a set $F_m\subset\Z^2$, and for any $\nu\in F_m,$ we have a component $C_m^\nu\subset X_m$. We prove that these are the irreducible components.

\begin{The}
For any $m\in\Z_{>0}$ we have that the space of $m$-jets of $X$ with center in the singular locus has the following decomposition into irreducible components
\[\pi_m^{-1}(X_{Sing})=\cup_{\nu\in F_m}C_m^\nu.\]
\end{The}

Hence, for
every $m \in \Z_{>0}$, we determine the irreducible components of $X_m$ with center in the singular locus of $X$. Moreover, we give a formula for the codimension of $C_m^\nu$ in Proposition \ref{Prop1}.

To prove the theorem we need to understand the geometry of the sets  $D_m^\nu$. This is done by comparing their geometries with
the geometries of the corresponding subsets for the approximate roots. The following proposition describes this comparison.

\begin{Pro}
For $m\in\Z_{>0}$ and $\nu\in\s\cap N_j$ such that $n_1\cdots n_g\langle\nu,\gamma_1\rangle\leq m$, we have that
\[D_m^\nu=(\pi_{m,[\frac{m}{e_j}]}^a)^{-1}(D_{j,[\frac{m}{e_j}]}^\nu)\]
where the integer $j\in\{1,\ldots,g\}$ satisfies the relation
\[n_j\cdots n_g\langle\nu,\gamma_j\rangle\leq m<n_{j+1}\cdots n_g\langle\nu,\gamma_{j+1}\rangle,\]
by $D_{j,m}^\nu$ we denote the corresponding set $D_m^\nu$ for the quasi-ordinary surface defined by the $j$-th approximated root, and for $q>p$, $\pi_{q,p}^a:\A^3_q\longrightarrow\A^3_p$ is the projection of the jet schemes of the affine ambient space.
\end{Pro}

This explains in part how the singular locus of the approximated roots play a role in our problem.

\vspace{3mm}

The irreducible components of the jet schemes fit in natural projective systems, to which we associate a weighted graph. The vertices of the graph correspond to irreducible components, and to every vertex we attach the corresponding codimension and embedded dimension. We will prove the following result.
\begin{The}
The weighted graph determines and it is determined by the topological type of the singularity.
\end{The}
This theorem achieves one of our goals for this type of singularities: constructing a complete invariant of the singularity from its jet schemes. Note that other invariants involving
arcs and jets, like motivic zeta functions, do not determine the topological type in the case of quasi ordinary singularities, see \cite{CoGPqo} and \cite{Nuelo}.

\vspace{3mm}

In another direction, using Mustata's formula (\cite{Mus}), we will determine an irreducible component of an $m$-jet scheme, or equivalently a divisor on the ambient space $\A^3$, which contributes the log canonical threshold
of the pair $X\subset\A^3$ (note that the log canonical threshold for such a pair has been computed in  \cite{BGG}, by looking at the poles of the motivic zeta function). This provides us
with pairs $X\subseteq\A^3$ whose log canonical threshold is not contributed by a monomial divisorial valuation. For instance, for the quasi-ordinary surface defined by $f=(z^2-x_1x_2)^2-x_1^3x_2z$,
the log canonical thresholds satisfy this property. Note that for a pair $C\subseteq\A^2$, where $C$ is a plane curve, the log canonical threshold is always contributed by a monomial valuation. See \cite{AN} and \cite{ACLM2} for the computation of the log canonical threshold for plane curves.

\vspace{3mm}

Along with the same ideas of \cite{Gen}, we are working to construct an embedded resolution of singularities of $X$ from the data of the graph
 constructed in this paper. We think that such a resolution puts light on the resolution of singularities obtained by González Pérez in \cite{GP4}, and give in
 some sense an answer to the question of Lipman (\cite{Lipman-Eq}) on the construction of a canonical resolution of singularity of a quasi-ordinary hypersurface from the
 characteristic exponents.
Moreover, understanding the surface case is an important step  in the understanding of the general case.

\vspace{3mm}

The structure of the paper is as follows. In Section \ref{Sec2} we introduce jet schemes. A brief exposition on quasi-ordinary singularities is given in Section
\ref{secQO}, together with some useful definitions at the end of the section. Section \ref{Sec4} is the heart of the paper; it is devoted to the study of the irreducible components of the
jet schemes of quasi-ordinary surface singularities. In Section \ref{Proofs} we state and proof some results which are useful but technical, and moreover we leave the proofs of some previous results, to make Section \ref{Sec4} more readable.

\vspace{3mm}

{\bf Acknowledgments.}
We thank Pedro González Pérez and the referees, for comments and suggestions which improved substantially the presentation of this paper.
The beginning of this work was done during a stay of HC supported by the
European Research Council under the European Community’s Seventh Framework
Programme (FP7/2007-2013) / ERC Grant Agreement nr. 246903 NMNAG.  She thanks Institut Mathématique de
Jussieu for hospitality.
HM is partially supported by ANR-12-JS01-0002-01 SUSI.

\section{Jet schemes}
\label{Sec2}

In this section we define jet schemes of an affine scheme $X,$ see \cite{Ishii-07} for details.
Let $X=\mbox{Spec }\C[x_1,\ldots,x_n]/I$ be an affine scheme of finite type. For $m\in\Z_{>0}$ the functor
$F_m:\ \C\mbox{-Schemes}\longrightarrow \mbox{Sets}$ which, with an affine scheme defined by a $\C$-algebra $A$, associates
\[F_m(\mbox{Spec}(A))=\mbox{Hom}_\C(\mbox{Spec}(A[t]/(t^{m+1})),X),\]
is representable by a $\C$-scheme, denoted by $X_m$. This is the scheme of $m$-jets. Its closed points are morphisms of the form
\[\gamma:\mbox{Spec}(\C[t]/(t^{m+1}))\longrightarrow X.\]
Such a morphism $\gamma$ is equivalent to a $\C$-algebra homomorphism
\[\gamma^\ast:\C[x_1,\ldots,x_n]/I\longrightarrow \C[t]/(t^{m+1}).\]
If we fix a set of generators $f_1,\ldots,f_r$ for the ideal $I$, the map $\gamma^\ast$ is determined by the image of the $x_i^{'s}$
\[x_i\mapsto x_i^{(0)}+x_i^{(1)}t+\cdots+x_i^{(m)}t^m,\ 1\leq i\leq n,\]
where the relations
\begin{equation}
f_i(x_1^{(0)}+\cdots+x_1^{(m)}t^m,\ldots,x_n^{(0)}+\cdots+x_n^{(m)}t^m)\equiv 0\mbox{ mod }t^{m+1}
\label{relation}
\end{equation}
must hold for each $f_i$, with $1\leq i\leq r$.
If we write equations in (\ref{relation}) as
\begin{equation}
\begin{array}{cc}
f_i(x_1^{(0)}+x_1^{(1)}t+\cdots+x_1^{(m)}t^m,\ldots,x_n^{(0)}+x_n^{(1)}t+\cdots+x_n^{(m)}t^m)=\\
\\
=\sum_{j=0}^mF_i^{(j)}(x_1^{(0)},\ldots,x_1^{(j)},\ldots,
x_n^{(0)},\ldots,x_n^{(j)})\ t^j\mbox{ mod }t^{m+1},\\
\end{array}
\label{expf}
\end{equation}
we have that giving a closed point of $X_m$ is equivalent to giving a point in $V(F_l^{(j)})_{0\leq j\leq m,1\leq l\leq r}\subset\A_m^n$, where
$\A_m^n=\mbox{Spec}(\C[x_i^{(0)},\ldots,x_i^{(m)}]_{i=1,\ldots,n})$. Hence we can make the
following identification
\begin{equation}
X_m=\mbox{Spec }\left(\frac{\C[x_i^{(0)},\ldots,x_i^{(m)}]_{i=1,\ldots,n}}{(F_l^{(j)})_{0\leq j\leq m,\ 1\leq l\leq r}}\right).
\label{Desc1}
\end{equation}

\begin{Exam}
Let $X$ be the quasi-ordinary surface defined by the polynomial $f=z^3-x_1^3x_2^2$.  The equations defining the $3$-jets are
\[\begin{array}{ll}
F^{(0)} & = {z^{(0)}}^3-{x_1^{(0)}}^3{x_2^{(0)}}^2\\
\\
F^{(1)} & =3{z^{(0)}}^2z^{(1)}-3{x_1^{(0)}}^2x_1^{(1)}{x_2^{(0)}}^2-2{x_1^{(0)}}^3{x_2^{(0)}}x_2^{(1)}\\
\\
F^{(2)} & =3{z^{(0)}}^2z^{(2)}+3z^{(0)}{z^{(1)}}^2-6{x_1^{(0)}}^2{x_1^{(1)}}x_2^{(0)}x_2^{(1)}-2{x_1^{(0)}}^3x_2^{(0)}x_2^{(2)}-3{x_1^{(0)}}^2x_1^{(2)}{x_2^{(0)}}^2\\
 & \ \ -{x_1^{(0)}}^3{x_2^{(1)}}^2-3x_1^{(0)}{x_1^{(1)}}^2{x_2^{(0)}}^2\\
\\
F^{(3)} & = {z^{(1)}}^3+6z^{(0)}z^{(1)}z^{(2)}+3{z^{(0)}}^2z^{(3)}-2{x_1^{(0)}}^3x_2^{(0)}x_2^{(3)}-2{x_1^{(0)}}^3x_2^{(1)}x_2^{(2)}-6{x_1^{(0)}}^2x_1^{(1)}x_2^{(0)}x_2^{(2)}\\
 & \ \ -3{x_1^{(0)}}^2x_1^{(1)}{x_2^{(1)}}^2-6{x_1^{(0)}}^2x_1^{(2)}x_2^{(0)}x_2^{(1)}-6x_1^{(0)}{x_1^{(1)}}^2x_2^{(0)}x_2^{(1)}-3{x_1^{(0)}}^2x_1^{(3)}{x_2^{(0)}}^2\\
 & \ \ -6x_1^{(0)}x_1^{(1)}x_1^{(2)}{x_2^{(0)}}^2-{x_1^{(1)}}^3{x_2^{(0)}}^2\\

\end{array}\]
\label{Ex0}
\end{Exam}

\begin{Rem}
Every polynomial $F^{(l)}$ is non-zero and quasi-homogeneous of degree $l$ in the variables $x_i^{(0)},\ldots,x_i^{(l)}$, $z^{(0)},\ldots,z^{(l)}$ for $i=1,2$. In $F^{(0)},\ldots,F^{(l)}$ the variables $x_1^{(l)},\ x_2^{(l)}$ and $z^{(l)}$ appear only in $F^{(l)}$.
\label{RemTonto}
\end{Rem}

\section{Quasi-ordinary surface singularities}
\label{secQO}

In this section we collect some well known facts about quasi-ordinary hypersurface singularities of dimension two, and we prove some lemmas which will be used in the next section. We state everything for the case of dimension two, though the definitions and results hold in any dimension.

An equidimensional germ $(X, 0)$, of dimension $2$, is {\em quasi-ordinary} (q.o. for short) if there exists a finite projection $p: (X, 0) \rightarrow (\C^2,0)$ which is a local isomorphism outside a normal crossing divisor. If $(X, 0)$ is a hypersurface there is an embedding $(X, 0) \subset (\C^{3},0) $, where $X$ is
defined by an equation $f= 0$,  and $f \in \C \{ x_1,x_2  \} [z]$ is a {\it quasi-ordinary ~polynomial}; that is, a Weierstrass polynomial with discriminant
$\D_z f$ of the form $\D_z f = x_1^{\d_1}\cdot x_2^{\d_2} \epsilon$ for a unit $\epsilon$ in the ring $ \C \{ x_1,x_2 \}$ of convergent power series and
$(\d_1,\d_2) \in \Z^2_{\geqslant 0}$. In these coordinates the projection $p$ is the restriction of the projection
\[\C^{3}\rightarrow\C^2,\ \ \ \ (x_1,x_2,z)\mapsto(x_1,x_2).\]
From now on we assume that $(X,0)$ is analytically irreducible, that is $f \in \C \{ x_1,  x_2  \} [z]$ is irreducible (see \cite{Ass} and \cite{Ev} for  criteria of irreducibility of q.o. polynomial). The Jung-Abhyankar theorem guarantees
that the roots of a q.o.~ polynomial $f$, called {\it q.o.~branches}, are fractional power series in $\C \{ x_1^{1/n},x_2^{1/n}\}$, for $n =\deg f$ (see
\cite{Abhyankar}). The difference $\zeta^{(i)}-\zeta^{(j)}$ of two different roots of $f$ divides the discriminant of $f$ in the ring $\C\{x_1^{1/n},
x_2^{1/n}\}$. Therefore $\zeta^{(i)}-\zeta^{(j)}=x_1^{\l_{ij}^{(1)}}x_2^{\l_{ij}^{(2)}}u_{ij}$ where $u_{ij}$ is a unit in $\C\{x_1^{1/n},x_2^{1/n}\}$. The exponents $\l_{ij}$
are characterized in the following Lemma:

\begin{Lem} \label{expo}     {\rm (see  \cite{Gau}, Prop. 1.3)}
Let  $f \in \C \{ x_1, x_2 \} [z] $ be an irreducible
q.o.~polynomial. Let $\z$ be a root of $f$ with expansion:
\begin{equation} \label{expan}
\z = \sum \beta_{\l} \bold{x}^\l.
\end{equation}
There exists $0 \ne \l_1, \dots, \l_g \in \Q^2_{\geqslant 0}$
such that $\l_1\leq\l_2\leq\cdots\leq\l_g$, and if $M_0 :=\Z^2 $ and $M_j := M_{j-1} + \Z \l_j$ for
$j=1, \dots, g$, then:
\begin{enumerate}
\item [(i)] $\beta_{\l_i} \ne 0$ and if $\beta_{\l} \ne 0$  then $\l \in M_j$ where
    $j$ is the unique integer such that  $\l_j \leqslant \l$ and
    $\l_{j+1} \nleq \l$ (where $\leqslant$ means coordinate-wise and we convey that $\l_{g +1} = \infty$).
\item [(ii)] For $j=1, \dots, g$, we have   $\l_j \notin M_{j-1}$,   hence the index
$n_j = [M_{j-1} : M_j]$ is $> 1$.
\end{enumerate}
Moreover if $\zeta\in\C\{x_1^{1/n},x_2^{1/n}\}$ is a fractional power series satisfying the conditions above, then $\zeta$ is a quasi-ordinary branch.
\end{Lem}

\begin{Defi}
The exponents $\l_1 , \dots, \l_g $ in Lemma \ref{expo} are called {\em characteristic exponents} of the q.o.~branch $\z$. We denote by  $M$  the lattice
$M_g$ and we call it the lattice associated to the q.o.~branch $\z$. We denote by $N$ (resp. $N_i$) the dual lattice of $M$ (resp. $M_i$ for $i=1,\ldots,g$). For
convenience we denote $\l_0:=(0,0)$ and $n_0:=1$.
\label{defExp}
\end{Defi}

In \cite{Gau} Gau proved that the characteristic exponents determine and are determined by the embedded topological type of $(X,0)$.
As a consequence of Lemma \ref{expo} we have the following result:

\begin{Lem}
If $\zeta$ is a quasi-ordinary branch of the form (\ref{expan}) then the series $\zeta_{j-1}:=\sum_{\lambda\not\geq\lambda_j}\beta_\lambda\bold{x}^\lambda$ is a
quasi-ordinary branch with characteristic exponents $\l_1,\ldots,\l_{j-1}$, for $j=1,\ldots,g$.

\label{LemConseq}
\end{Lem}

\begin{Defi}
For $0\leq j\leq g-1$ we have the germ of quasi-ordinary hypersurface $(X^{(j)},0)$, where $X^{(j)}$ is parametrized by the branch $\zeta_j$. For convenience we also denote $\zeta$ by $\zeta_g$ and $X$ by $X^{(g)}$.
\label{SemiX}
\end{Defi}

Without loss of generality we relabel the variables $x_1, x_2$ in such a way that if  $\l_j = ( \l_j^{(1)}, \l_j^{(2)}) \in
\Q^2$ for $j=1, \dots, g$, then we have:
\begin{equation}
(\l_1^{(1)}, \dots, \l_g^{(1)}) \geqslant_{\mbox{\rm lex}}  (\l_1^{(2)}, \dots, \l_g^{(2)}),
\label{lex}
\end{equation}
where $\geqslant_{\mbox{\rm lex}}$ is lexicographic order. The
q.o.~branch $\z$
 is said to be normalized if
$\l_1$ is not of the form $(\l_1^{(1)}, 0)$ with $\l_1^{(1)} <
1$. Lipman proved that the germ $(X,0)$ can be parametrized by a
normalized q.o.~branch (see  \cite{Gau}, Appendix).
 We assume from now on that the q.o.~branch $\z$ is normalized.

The semigroup $\Z^2_{\geqslant 0}$ has a minimal set of generators
$v_1,v_2$, which is a basis of the lattice $M_0$. The dual
basis, $\{w_1,w_2\}$, is a basis of the dual lattice $N_0$, and spans a regular cone $\s$ in
$N_{0,\R}=N_0\otimes_\Z\R$. It follows that $\Z^2_{\geqslant 0} = \s^\vee \cap
M_0$, where $\s^\vee = \R^2_{\geqslant 0}$ is the dual cone of
$\s$. The $\C$-algebra $\C \{  x_1,x_2 \} $ is isomorphic
to $\C \{ \s^\vee \cap M_0 \} $.  The
local algebra ${\mathcal O}_X = \C\{x_1, x_2\}[z]/(f)$ of the singularity $(X,0)$ is isomorphic to
$\C\{\s^\vee \cap M_0\}[\z]$. By Lemma \ref{expo} the series  $\z
$ can be viewed as an element $\sum \beta_\l \bold x^\l$ of the algebra
$ \C\{\s^\vee \cap M\}$.

\begin{Lem} (see \cite{GP4})
The homomorphism $\mathcal O_X\longrightarrow\C\{\s^\vee\cap M\}$ is the inclusion of $\mathcal O_X$ in its integral closure in its field of fractions.
\label{NormTor}
\end{Lem}

This Lemma shows that the normalization of a quasi-ordinary hypersurface $(X,0)$ is the germ of the toric variety $X(\s,N)=Z^{\s^\vee\cap M}$ at the distinguished point.

The elements of $M$ defined by:
\begin{equation}
{\g}_1 =  \l_1\mbox{ and } {\g}_{j+1}- n_j {\g}_{j} = \l_{j+1} -  \l_{j}\mbox{ for } j= 1, \dots, g-1,
\label{defSem}
\end{equation}
span the semigroup $\Gamma := \Z^2_{\geqslant 0} +  \g_1 \Z_{\geqslant
0} + \cdots + \g_g \Z_{\geqslant 0} \subset \s^\vee \cap M$. For convenience we denote $\gamma_0:=0$. The
semigroup $\Gamma$ defines an analytic invariant of the germ
$(X,0)$ (see \cite{Jussieu},\cite{PPP04},\cite{KM}).

\begin{Defi}
The monomial variety associated to $(X,0)$ is the toric variety
\[X^\Gamma:=\mbox{Spec }\C[\Gamma].\]
\label{defMV}
\end{Defi}
 Moreover we associate with the characteristic exponents the following sequence of semigroups:
\[\Gamma_j=\s^\vee\cap M+\gamma_1\Z_{\geq 0}+\cdots+\gamma_j\Z_{\geq 0},\mbox{ for }j=0,\ldots,g.\]

And we have the corresponding monomial varieties associated to $\Gamma_j$. We denote by $e_{i-1}:=n_i\cdots n_g$ for $1<i\leq g$ and set $e_g:=1$.
Notice that, by (\ref{lex}) and the definition of $\gamma_1,\ldots,\gamma_g$, we deduce that
\begin{equation}
(\gamma_1^{(1)}, \dots, \gamma_g^{(1)}) \geqslant_{\mbox{\rm lex}}(\gamma_1^{(2)}, \dots, \gamma_g^{(2)}).
\label{lexGamma}
\end{equation}

The following Lemma gathers some important facts about the generators $\gamma_j$ and the semigroups $\Gamma_j$.

\begin{Lem}(see Lemma 3.3 in \cite{Jussieu})
\begin{enumerate}
\item[(i)] We have that $\gamma_j>n_{j-1}\gamma_{j-1}$ for $j=2,\ldots,g$, where $<$ means $\neq$ and $\leq$ coordinate-wise.
\\
\item[(ii)] If a vector $u_j\in\s^\vee\cap M_j$, then we have $u_j+n_j\gamma_j\in\Gamma_j$.
\\
\item[(iii)] The vector $n_j\gamma_j$ belongs to the semigroup $\Gamma_{j-1}$ for $j=1,\ldots,g$. Moreover, we have a unique relation
\begin{equation}
n_j\gamma_j=\alpha^{(j)}+r_1^{(j)}\gamma_1+\cdots+r_{j-1}^{(j)}\gamma_{j-1}
\label{relgam}
\end{equation}
such that $0\leq r_i^{(j)}\leq n_i-1$ and $\alpha^{(j)}\in M_0$ for $j=1,\ldots,g$.
\end{enumerate}
\label{LemaPedro}
\end{Lem}

\begin{Defi}
Given two irreducible quasi-ordinary polynomials $f$ and $g$ in $\C\{x_1,x_2\}[z]$ such that $fg$ is a quasi-ordinary polynomial, we say that $f$ and $g$
have order of coincidence $\alpha\in \Q^2$ if $\alpha$ is the largest exponent on the set
\[\{\l_{ij}\ |\ f(\zeta^{(i)})=g(\zeta^{(j)})=0\},\]
where $\zeta^{(i)}$ and $\zeta^{(j)}$ are roots of $fg$.
\label{OrdCon}
\end{Defi}

\begin{Defi}
 We associate to $f$ a set of semi-roots
\[z=f_0,f_1\ldots,f_g=f\in \C\{x_1,x_2\}[z].\]
Every $f_j$ is an irreducible quasi-ordinary polynomial of degree $n_0\cdots n_j$ with order of coincidence with $f$ equal to $\l_{j+1}$ for $j=0,\ldots,g$.
\label{defAroots}
\end{Defi}

They are parametrized by truncations of a root $\zeta(x_1^{1/n},x_2^{1/n})$ of $f$ in the following sense:

\begin{Pro}
(see \cite{Jussieu}) Let $q\in \C\{x_1,x_2\}[z]$ be a monic polynomial of degree $n_0\cdots n_j$. Then $q$ is a $j$-th semi-root of $f$ if and only if
$q(\zeta)=\bold{x}^{\gamma_j}\epsilon_j$ for a unit $\epsilon_j$ in $\C\{x_1,x_2\}[z]$.

\end{Pro}

\begin{Cor}
The quasi-ordinary polynomials $f_j\in\C\{x_1,x_2\}[z]$ defining $X^{(j)}$ (see Definition \ref{SemiX}) for $j=0,\ldots,g$ form a system of semiroots of $f$.
\label{CorSR}
\end{Cor}

In what follows we state some results about quasi-ordinary polynomials and approximated roots. Moreover we give some definitions and notations that will be
used in the next section.

Approximated roots play an important role in the understanding of quasi-ordinary singularities. We have the following expansions of the semiroots in terms of the previous ones:

\begin{Lem} (See Lemma 35 in \cite{GP4})
The expansion of the approximated roots is of the following form:
\begin{equation}
c_j^*f_j=f_{j-1}^{n_j}-c_jx_1^{\alpha_1^{(j)}}x_2^{\alpha_2^{(j)}}f_0^{r_1^{(j)}}\cdots f_{j-2}^{r_{j-1}^{(j)}}+\sum c_{\underline\alpha,\underline{r}}x_1^{\alpha_1}x_2^{\alpha_2}f_0^{r_1}\cdots f_{j-1}^{r_j},
\label{expSR}
\end{equation}
where $c_j^*,c_j\in\C^*$, $0\leq r_i^{(j)},r_i<n_i$ for $i=1,\ldots,j$, and \[n_j\gamma_j=(\alpha_1^{(j)},\alpha_2^{(j)})+r_1^{(j)}\gamma_1+\cdots+r_{j-1}^{(j)}\gamma_{j-1}<(\alpha_1,\alpha_2)+r_1\gamma_1+\cdots+r_j\gamma_j.\]
\label{Lema35}
\end{Lem}
Let $(X,0)\subset(\C^3,0)$ be a germ of quasi-ordinary surface with characteristic exponents $\l_1,\ldots,\l_g$. We denote by $\l_1=(\frac{a_1}{n_1},\frac{b_1}{n_1})$ the first
characteristic exponent.
Notice that, by (\ref{lex}), we have that $a_1\geq b_1\geq 0$.

\begin{Lem}
We have that
\begin{equation}
f_1=z^{n_1}-c_1x_1^{a_1}x_2^{b_1}+\sum c_{ijk}x_1^ix_2^jz^k,
\label{eqSR1}
\end{equation}
with $(i,j)+k\gamma_1>n_1\gamma_1$ and $k<n_1$ whenever $c_{ijk}\neq 0$.
And for $1\leq l\leq g-1$ we have
\begin{equation}
f=f_l^{e_l}+\sum_{(i,j)+k\gamma_1>n_le_l\gamma_l}c_{ijk}^{(l)}x_1^ix_2^jz^k.
\label{eqSRpqprima}
\end{equation}
Moreover, the following expansions will be useful. For $0\leq j<g-1$
\begin{equation}
f =f_j^{e_j}+\sum c_{\underline{\alpha},\underline{r}}x_1^{\alpha_1}x_2^{\alpha_2}z^{r_1}\cdots f_j^{r_{j+1}}
\label{grgr1}
\end{equation}
where $(\alpha_1,\alpha_2)+r_1\gamma_1+\cdots+r_{j+1}\gamma_{j+1}>n_{j+1}e_{j+1}\gamma_{j+1}$ whenever $c_{\underline{\alpha},\underline{r}}\neq 0$.
\begin{equation}
f =f_j^{e_j}+\sum c_{\underline{\alpha},\underline{r}}x_1^{\alpha_1}x_2^{\alpha_2}z^{r_1}\cdots f_{j-1}^{r_j}
\label{grgr2}
\end{equation}
where $(\alpha_1,\alpha_2)+r_1\gamma_1+\cdots+r_j\gamma_j>n_je_j\gamma_j$ whenever $c_{\underline{\alpha},\underline{r}}\neq 0$.
\label{LemExpSR}
\end{Lem}
{\em Proof.} By applying recursively Lemma \ref{Lema35} and using Lemma \ref{LemaPedro} (i).\hfill $\Box$

\subsection{Some toric geometry}
\label{toric}
See \cite{F} for a reference on toric geometry.

Given a lattice $N$ we denote by $N_\R$ the vector space spanned by $N$ over the field $\R$. We denote by $M$ the dual lattice, $M=\mbox{Hom}(N,\Z)$, and by $\langle,\rangle : N\times M\longrightarrow\Z$ the duality pairing between the lattices $N$ and $M$.
 A rational convex polyhedral cone (or simply a cone) is the set of non-negative linear combinations of vectors $v_1,\ldots,v_r\in N$. A cone is strictly convex if it contains no lines. The dual cone of $\s$, denoted by $\s^\vee$, is the set
 \[\s^\vee=\{\omega\in M_\R\ |\ \langle\omega,u\rangle\geq 0\ \forall u\in\s\},\]
 and the orthogonal of $\s$, denoted by $\s^\bot$, is
 \[\s^\bot=\{\omega\in M_\R\ |\ \langle\omega,u\rangle=0\ \forall u\in\s\}.\]
 We denote by $\stackrel{\circ}{\s}$ the relative interior of the cone $\s$. A fan $\Sigma$ is a family of strictly convex cones in $N_\R$ such that for any $\s\in\Sigma$ any face of $\s$ belongs to $\Sigma$, and for any $\s,\t\in\Sigma$, the intersection $\s\cap\t$ is a face of both.
  The relation $\t\leq\s$ denotes that $\t$ is a face of $\s$. The support of the fan $\Sigma$ is the set $|\Sigma|:=\cup_{\t\in\Sigma}\t\subset N_\R$.

  Let $\t$ be a strictly convex cone, rational for the lattice $N$. By Gordan's Lemma the semigroup $\s^\vee\cap M$ is finitely generated. We denote by $\C[\s^\vee\cap M]$ the semigroup algebra of $\s^\vee\cap M$ with complex coefficients. The toric variety $Z(\t,N):=\mbox{Spec }\C[\t^\vee\cap M]$ is normal. The torus $T_N:=Z(N)$ is an open dense subset of $Z(\t,N)$ which acts on $Z(\t,N)$ and the action extends the action of the torus on itself by multiplication. There is a one to one correspondence between the faces $\theta$ of $\t$ and the orbits orb$(\theta)$ of the torus action on $Z(\t,N)$, which reverses the inclusions of their closures. The closure of orb$(\theta)$ is the toric variety $Z((\s^\vee\cap\theta^\bot)^\vee,N)$ for $\t\leq\s$.

\subsection{Definitions and Notations}
\label{defsnots}
We introduce now some definitions and notations which will be used throughout the paper.

The singular locus of a quasi-ordinary singularity is determined, after Lipman, by its characteristic exponents (see \cite{Lipman2} and \cite{PPP04}).

\begin{Defi} We define
\[\begin{array}{ll}
Z_i=X\cap\{x_i=0\}, & \mbox{for }i=1,2\\
\\
Z_{12}=X\cap\{x_1=x_2=0\}. & \\
\end{array}\]
Moreover, the smallest number $c\in\{1,2\}$ with the property that
\[\l_i^{(j)}=0,\mbox{ for all }1\leq i\leq g\mbox{ and }c+1\leq j\leq 2\]
is called the equisingular dimension of the quasi-ordinary projection $p$.
\label{equi}
\end{Defi}
In \cite{Lipman2}
Lipman proved that the spaces $Z_1,Z_2$ and $Z_{12}$ are irreducible.
By condition (\ref{lex}) we have that $c$ gives the number of variables appearing in the monomials $\bold x^{\l_1},\ldots,\bold x^{\l_g}$.

\begin{Defi}
Let $X$ be a quasi-ordinary surface singularity with $g\geq 1$ characteristic exponents. We define the integers $g_1\geq 0$ and $g_2\in\{g_1,g_1+1\}$ as follows
\[\mbox{ if }c=1\mbox{ we set }g_1=g_2=g+1,\]
otherwise
\[\begin{array}{l}
\gamma_{g_1}^{(2)}=0\mbox{ and }\gamma_{g_1+1}^{(2)}\neq 0,\\
\\
g_2=\left\{\begin{array}{cl}
g_1+1 & \mbox{ if }\gamma_{g_1+1}^{(2)}=\frac{1}{n_{g_1+1}}\\
\\
g_1 & \mbox{ otherwise}\\
\end{array}\right.\\
\end{array}\]
\label{defg1}
\end{Defi}

Note that these integers can be defined with the same property for the characteristic exponents. Now we use them to describe the singular locus of $X$.
Lipman's theorem describes the singular locus $X_{Sing}$ of a quasi-ordinary hypersurface $X$. We state it here in the particular case of surfaces.

\begin{The} (See Theorem 7.3 in \cite{Lipman2})
Let $X$ be a quasi-ordinary surface singularity with characteristic exponents $\l_1,\ldots,\l_g$. Then we have:
\begin{enumerate}
\item[(i)] $X_{Sing}=Z_{12}$ if and only if $g=1$ and $\l_1=(\frac{1}{n},\frac{1}{n})$.
\item[(ii)] If $c=1$ then $X_{Sing}=Z_1$.
\item[(iii)] Otherwise $c=2$, and since $\l_1^{(1)}\neq 0$, $Z_1\subset X$ is a component of $X_{Sing}$. Moreover, if we do not have
simultaneously $\l_k^{(2)}=0$ for all $1\leq k\leq g-1$ and $\l_g^{(2)}=\frac{1}{n_g}$ then the singular locus is reducible of the form $X_{Sing}=Z_1\cup Z_2$.
\end{enumerate}
\label{CorSing}
\end{The}

\begin{Rem}
Notice that
\[\begin{array}{l}
Z_1=\{x_1=z=0\}\\
\\
Z_2=\{x_2=f_{g_1}=0\}\\
\\
Z_{12}=\{(0,0,0)\}\\
\end{array}\]
and hence the singular locus of a quasi-ordinary surface singularity is either a point, or a line, or two lines, or a line and a singular curve.
\end{Rem}

Then, geometrically, the meaning of the integer $g_2$ is to measure the irreducibility of the singular locus of the approximated roots. Indeed,
\[\begin{array}{l}
X_{Sing}^{(j)} \mbox{ is irreducible, for }1\leq j\leq g_2\\
\\
X_{Sing}^{(j)} \mbox{ is reducible, for }g_2<j\leq g\\
\end{array}\]

Now we define a sequence of semi-open cones keeping track of the singular locus of the quasi-ordinary hypersurfaces $X^{(j)}$ for $j=1,\ldots,g$ (see Definition
\ref{SemiX}). Let
\[\nu_j:X(\s,N_j)\longrightarrow X^{(j)}\]
be the normalization of $X^{(j)}$ (see Lemma \ref{NormTor}). Consider $\nu_j^{-1}(X_{Sing}^{(j)})\subseteq X(\s,N_j)$, it is a disjoint union of orbits
\[\nu_j^{-1}(X_{Sing}^{(j)})=\bigsqcup_\t\mbox{orb}(\t),\ \mbox{ for some }\t\mbox{ faces of }\s.\]
 We also have that the complement of $\nu_j^{-1}(X_{Sing}^{(j)})$ in the toric variety $X(\s,N_j)$ is a union of orbits.
\begin{Defi}
We define, for $j=1,\ldots,g$
\[\begin{array}{ll}
\s_{Sing,j} & \mbox{ the fan associated to }\nu_j^{-1}(X_{Sing}^{(j)})\\
\\
\s_{Reg,j} & \mbox{ the fan associated to }X(\s,N_j)\setminus\nu_j^{-1}(X_{Sing}^{(j)})\\
\end{array}\]
For $j=g$ we just denote them by $\s_{Sing}$ and $\s_{Reg}$ respectively.
\label{DefTSigmas}
\end{Defi}

\begin{Rem}
Recall that $\s=\R_{\geq 0}^2$, let $\r_1$ and $\r_2$ be its one-dimensional faces. For $1\leq j\leq g$ we have
\[\s_{Sing,j}=\left\{\begin{array}{cll}
\s\setminus(\r_1\cup\r_2) & \mbox{ if } & X_{Sing}^{(j)}=Z_{12}\\
\\
\s\setminus\r_2 & \mbox{ if } & X_{Sing}^{(j)}=Z_1\\
\\
\s\setminus\{(0,0)\} & \mbox{ if } & X_{Sing}^{(j)}=Z_1\cup Z_2\\
\end{array}\right.\]
and by definition $\s_{Reg,j}=\s\setminus\s_{Sing,j}$.
\label{defSigmas}
\end{Rem}

The fan $\s_{Sing}$ will turn out to be necessary in our description of $\pi_m^{-1}(X_{Sing})$ (see Lemma \ref{LemGamma}), while the fans $\s_{Reg,j}$ will be important in the description of the irreducible components (see Proposition \ref{PropC1}).

The sequence
$\{\s_{Reg,1},\ldots,\s_{Reg,g}\}$ is not very complicated, in the sense that most of the elements are the same.
Since by definition $\gamma_{g_1+1}^{(2)}=\lambda_{g_1+1}^{(2)}$ then, by Theorem \ref{CorSing}, we deduce
\begin{equation}
\begin{array}{ll}
\mbox{for }1\leq j\leq g_2 & \sigma_{Reg,j}=\left\{\begin{array}{cl}
                                        \r_1\cup\r_2 & \mbox{ if }\gamma_1=(\frac{1}{n_1},\frac{1}{n_1})\mbox{ and }j=1\\
                                        \r_2 & \mbox{ otherwise}\\
                                        \end{array}\right.\\
\\
\mbox{for }g_2+1\leq j\leq g & \sigma_{Reg,j}=\{(0,0)\}\\
\end{array}
\label{SigRegj}
\end{equation}

\vspace{3mm}

Moreover notice that, by definition, we have $\s_{Sing,j}\subseteq\s_{Sing,j+1}$.

\vspace{3mm}

\begin{Defi}
Given $\nu\in\s\cap N_0$, we define the following sequence of real numbers
\[l_i(\nu):=n_ie_i\langle\nu,\gamma_i\rangle,\mbox{ for }1\leq i\leq g.\]
Set $l_0(\nu)=0$ and $l_{g+1}(\nu)=\infty$ for any $\nu\in\Z^2$.
Moreover, we define
\[i(\nu)=\left\{\begin{array}{cl}
g+1 & \mbox{ if }\nu\in N_g\\
\\
\mbox{min }\{i\in\{1,\ldots,g\}\ |\ \nu\notin N_i\} & \mbox{ otherwise}\\
\end{array}\right.\]
and
\[\begin{array}{l}
c(\nu)=\mbox{max }\{0\leq i\leq g\}\ |\ \langle\nu,\gamma_i\rangle=0\}\\
\\
m(\nu)=\mbox{min }\{1\leq i\leq g\ |\ n_i\langle\nu,\gamma_i\rangle<\langle\nu,\gamma_{i+1}\rangle\}\\
\end{array}\]
\label{defs}
\end{Defi}

Notice that, by definition, $l_1(\nu)=l_2(\nu)=\cdots=l_{m(\nu)}(\nu)<l_{m(\nu)+1}(\nu)$. Moreover, for $\nu\in\s_{Sing}\cap N_0$, we have that
\[\begin{array}{ccc}
c(\nu)=\left\{\begin{array}{cc}
g_1 & \mbox{ if }g_1>0\mbox{ and }\nu\in\r_2\\
\\
0 & \mbox{ otherwise}\\
\end{array}\right.
&
\mbox{ and }
&
m(\nu)=\left\{\begin{array}{cc}
\in\{1,\ldots,g\} & \mbox{ if }\nu\in\r_1\\
\\
g_1 & \mbox{ if } \nu\in\r_2\mbox{ and }g_1>0\\
\\
1 & \mbox{ otherwise}\\
\end{array}\right.
\end{array}\]
It is straightforward to check that $c(\nu)\leq m(\nu)\leq i(\nu)$.

\begin{Lem}
 For any $\nu\in\s_{Sing}\cap N_0$, we have that
 \[l_i(\nu)\in\Z \mbox{ for }1\leq i\leq\mbox{min }\{i(\nu),g\}.\]
 Moreover the integers $l_i(\nu)$ are ordered as
 \begin{equation}
 0= l_{c(\nu)}(\nu)< l_{c(\nu)+1}(\nu)\leq\cdots \leq l_g(\nu)<l_{g+1}(\nu)=\infty,
 \label{ineq}
 \end{equation}
and we have the equality $l_i(\nu)=l_{i+1}(\nu)$ if and only if $\gamma_{i+1}^{(j)}=n_i\gamma_i^{(j)}$ for $j$ either $1$ or $2$, and $\nu\in\r_j$.

For $i>c(\nu)$ we have that the following statements are equivalent:
\begin{enumerate}
\item[(i)] $\nu\in N_i$

\item[(ii)] for $1\leq j\leq i-1,\ \nu\in N_j$ and the number $\frac{l_{j+1}(\nu)-l_j(\nu)}{e_j}$ is a positive integer.
\end{enumerate}
\label{l-ord}
\end{Lem}

{\em Proof.}
For $1\leq i<i(\nu)$ we have that $\nu\in N_i$, therefore $\langle\nu,\gamma_i\rangle\in\Z$ and $l_i(\nu)$ is an integer. If $i(\nu)<g+1$, then
$\nu\in N_{i(\nu)-1}$ and by (\ref{relgam}) we deduce that $l_{i(\nu)}(\nu)$ is an integer.

By definition of $\gamma_i$ it follows that $e_i\gamma_{i+1}=e_{i-1}\gamma_i+e_i(\l_{i+1}-\l_i)$. Then
\[l_{i+1}(\nu)=l_i(\nu)+e_i\langle\nu,\l_{i+1}-\l_i\rangle,\]
and the second claim of the Lemma follows because the exponents $\l_i$ are ordered lexicographically as (\ref{lex}).

The equivalence follows by the fact that
\[l_j(\nu)-l_{j-1}(\nu)=e_{j-1}(\langle\nu,\gamma_j\rangle-n_{j-1}\langle\nu,\gamma_{j-1}\rangle).\]
\hfill $\Box$

\section{Jet schemes of quasi-ordinary surface singularities}
\label{Sec4}
In this section we describe the irreducible components of $\pi_m^{-1}(X_{Sing})\subset X_m$. We begin with an overview of the section.

We will associate, to any $\nu=(\nu_1,\nu_2)\in\s_{Sing}\cap N_0$ with $0\leq\nu_i\leq m$, a family of $m$-jets that we call
$C_m^\nu$. Roughly speaking, it is the Zariski closure of the set of $m$-jets whose order of contact with the hyperplane coordinate $x_i$ is bigger or equal to $\nu_i$, for $i=1,2$.
 We divide the sets $C_m^\nu$ into two types. Sets defined by the annihilation of hyperplane coordinates in $\A_m^3=
 \mbox{Spec }\C[x_1^{(i)},x_2^{(i)},z^{(i)}]_{i=0,\ldots,m}$, which are associated with $\nu$ in a certain set $H_m\subset\sigma_{Sing}\cap N_0$ (where $H$ stands for hyperplane). The $C_m^\nu$ of the second type have a more complicated geometry and are associated with $\nu$ in a certain set $L_m$ (where $L$ stands for
 lattice, because for such components, $\nu$ belongs to one of the lattices $N_i$ for $1\leq i\leq g$).

 The geometry of $C_m^\nu$ for $\nu\in H_m$ is evident and the possible inclusions among different $C_m^\nu$ are determined only by looking at the
 product ordering of the associated vectors $\nu$. For $\nu\in L_m$, to understand the geometry of $C_m^\nu$, we consider the order of contact of the generic point
 of $C_m^\nu$ with the different approximated roots. This allows us to detect certain dense subset of $C_m^\nu$, which is isomorphic to the cartesian product of an open set of the spectrum of the graded algebra (associated with the last approximated root appearing in the equations defining $C_m^{\nu}$) and an affine
 space. This permits to prove that each $C_m^\nu$ is irreducible and to compute its codimension (see Proposition \ref{Prop1}).

 The inclusions among the $C_m^\nu$ for $\nu\in H_m\cup L_m$ are more delicate. We introduce in Definition \ref{orden2} a new relation to detect such
 inclusions. In this definition, the singular locus of the last approximated root affecting the geometry of $C_m^\nu$, plays a crucial role.
Finally, with the collection of sets $C_m^\nu$ left, we prove in Theorem \ref{TheCaso1} that they are the irreducible components of the $m$-jets
through the singular locus.

\begin{Defi}
Let $h\in \C[x_1,\ldots,x_n]$ and let $X$ be an algebraic variety. For $p,m\in\Z_{>0}$ with $p\leq m$ we set
\[\mbox{Cont}_X^p(h)_m:=\{\gamma\in X_m\ |\ \mbox{ord}_t(h\circ\gamma)=p\}.\]
And, for $m\in\Z_{>0}$ and any $\nu=(\nu_1,\nu_2)\in\s_{Sing}\cap N_0$ with $\nu_i\leq m$ we define the constructible set
\[D_m^\nu=(\mbox{Cont}_X^{\nu_1}(x_1)_m\cap\mbox{Cont}_X^{\nu_2}(x_2)_m)_{red},\]
where we consider the reduced structure, and
\[C_m^\nu=\overline{D_m^\nu},\]
its Zariski closure. We denote by $D(f)$ the open set
\[D(f)=\mbox{Spec }R_f\]
where $R$ is the ring $R=\C[x_1^{(j)},x_2^{(j)},z^{(j)}]_{j\geq 0}$.
\label{DefDC}
\end{Defi}

Given a jet $\gamma\in X_m$, if $x_i\circ\gamma\neq0$ for $i=1,2$, the vector $\nu=(\mbox{ord}_t(x_1\circ\gamma),\mbox{ord}_t(x_2\circ\gamma))$ belongs to $\s\cap N_0$ and $0\leq\nu_i\leq m$. Moreover it is trivial that $\gamma\in C_m^\nu$. Now we look at $m$-jets with origin at the singular locus.

\begin{Lem}
Given $\gamma\in\pi_m^{-1}(X_{Sing})$, there exists $\nu\in\s_{Sing}\cap N_0$ with $0\leq\nu_i\leq m$ for $i=1,2$, such that $\gamma\in C_m^\nu$.
Moreover
\[\pi_m^{-1}(X_{Sing})=\bigcup_{\nu\in[0,m]^2\cap\s_{Sing}\cap N_0}C_m^\nu.\]
\label{LemGamma}
\end{Lem}

{\em Proof.}
 Given $\gamma\in\pi_m^{-1}(X_{Sing})$, suppose first that $x_i\circ\gamma\neq 0$ for $i=1,2$. Then $\nu:=(\mbox{ord}_t(x_1\circ\gamma),\mbox{ord}_t(x_2\circ\gamma))\in[0,m]^2$ and obviously $\gamma\in D_m^\nu\subseteq C_m^\nu$. We have to prove that $\nu\in\s_{Sing}\cap N_0$, and this follows by Remark \ref{defSigmas}, since:

\begin{enumerate}
\item[(i)] If $X_{Sing}=\{(0,0,0)\}$, then $\gamma(0)=(0,0,0)$ and ord$_t(x_i\circ\gamma)>0$ for $i=1,2$.

\item[(ii)] If $X_{Sing}=Z_1$, then $\gamma(0)=(0,x_2(0),0)$, and $\mbox{ord}_t(x_1\circ\gamma)>0,\ \mbox{ord}_t(x_2\circ\gamma)\geq 0$.

\item[(iii)] If $X_{Sing}=Z_1\cup Z_2$, then
\[\mbox{if }\gamma(0)\in Z_1, \mbox{ we have ord}_t(x_1\circ\gamma)>0\mbox{ and ord}_t(x_2\circ\gamma)\geq 0\]
\[\mbox{if }\gamma(0)\in Z_2, \mbox{ we have ord}_t(x_1\circ\gamma)\geq 0\mbox{ and ord}_t(x_2\circ\gamma)>0\]
\end{enumerate}

To deal with the other cases, notice that
\[C_m^\nu=\{\gamma\in X_m\ |\ \mbox{ord}_t(x_i\circ\gamma)\geq\nu_i,\ i=1,2\}.\]

If $x_i\circ\gamma=0$ for $i=1,2$, then $\gamma\in C_m^\nu$ for any $\nu\in\s_{Sing}\cap N_0$ with $0\leq\nu_i\leq m$ for $i=1,2$.

\vspace{2mm}

If $x_1\circ\gamma=0$ and $x_2\circ\gamma\neq 0$, then we denote $\alpha:=\mbox{ord}_t(x_2\circ\gamma)$. We have $0\leq\alpha\leq m$, and $\gamma\in C_m^\nu$
for any $\nu\in\s_{Sing}\cap N_0$, with $0\leq\nu_i\leq m$ for $i=1,2$, and $\nu_2\leq\alpha$.

\vspace{2mm}

The left case $x_1\circ\gamma\neq 0$ and $x_2\circ\gamma=0$ is analogous to the last one.

\vspace{2mm}

We prove the other inclusion. If $\gamma\in X_m\setminus\pi_m^{-1}(X_{Sing})$, then $\gamma(0)\notin X_{Sing}$. Again distinguishing cases depending on the singular locus, we can prove that $\nu=(\mbox{ord}_t(x_1\circ\gamma),\mbox{ord}_t(x_2\circ\gamma))\notin\s_{Sing}$.
\hfill$\Box$

\subsection{Description of the sets $C_m^\nu$}

The sets $C_m^\nu$ are the candidates to be irreducible components of $\pi_m^{-1}(X_{Sing})$. We proceed to study these sets. Notice that, by definition, it follows that
\[C_m^\nu\subset V(x_i^{(0)},\ldots,x_i^{(\nu_i-1)},i=1,2).\]

\begin{Defi}
For  $\nu\in\s\cap N_0$ and $m\in\Z_{>0}$ we define the ideals
\[I^\nu=\left(x_i^{(0)},\ldots,x_i^{(\nu_i-1)}\right)_{i=1,2},\]
\[J_m^\nu=Rad\left((F^{(i)}\mbox{ mod }I^\nu)_{0\leq i\leq m}\right).\]
Moreover we define the integers $j(m,\nu)\in\{0,\ldots,i(\nu)-1\}$ and $j'(m,\nu)\in\{0,\ldots,j(m,\nu)\}$ by the
inequalities
\[\begin{array}{c}
l_j(\nu)\leq m< l_{j+1}(\nu),\\
\\
l_j(\nu)+e_j\leq m<l_{j+1}(\nu)+e_{j+1},\\
\end{array}\]
respectively.
\label{IJ}
\end{Defi}

Then  we have that
\begin{equation}
D_m^\nu=V(I^\nu,J_m^\nu)\cap D(x_1^{(\nu_1)})\cap D(x_2^{(\nu_2)}),
\label{eqast}
\end{equation}
where the fact of taking the radical in the definition of $J_m^\nu$ corresponds to taking the reduced structure in the definition of $D_m^\nu$.
We have to study the polynomials $F^{(j)}\mbox{ mod }I^\nu$ for $0\leq j\leq m$. Thanks to the identification in (\ref{Desc1}), these polynomials have to be seen as the defining equations of $D_m^\nu$.

\begin{Exam}
Let $X$ be the quasi-ordinary surface defined by the polynomial $f=z^3-x_1^3x_2^2$. The singular locus of $X$ is reducible
\[X_{Sing}=\{x_1=z=0\}\cup\{x_2=z=0\}.\]
We described in Example \ref{Ex0} the equations of the $3$-jets, and if we look at $3$-jets with origin in the singular locus, then we have to add the condition either $x_1^{(0)}=z^{(0)}=0$ or $x_2^{(0)}=z^{(0)}=0$. This is equivalent to consider the equations modulo the ideal $I^\nu$.
Then
\[\begin{array}{lcl}
(\pi_1^{-1}(X_{Sing}))_{red} & = & \pi_1^{-1}(\{x_1=z=0\})\cup\pi_1^{-1}(\{x_2=z=0\})\\
& = & V(x_1^{(0)},z^{(0)})\cup V(x_2^{(0)},z^{(0)})=C_1^{(1,0)}\cup C_1^{(0,1)}\subset\A_1^3.\\
\\
(\pi_2^{-1}(X_{Sing}))_{red} & = & \pi_2^{-1}(\{x_1=z=0\})\cup\pi_2^{-1}(\{x_2=z=0\})\\
& = & V(x_1^{(0)},z^{(0)})\cup (V(x_2^{(0)},z^{(0)},{x_1^{(0)}}^3{x_2^{(1)}}^2))_{red}=\\

 & = & V(x_1^{(0)},z^{(0)})\cup V(x_1^{(0)},x_2^{(0)},z^{(0)})\cup V(x_2^{(0)},x_2^{(1)},z^{(0)})=\\

 & = & V(x_1^{(0)},z^{(0)})\cup V(x_2^{(0)},x_2^{(1)},z^{(0)})=C_2^{(1,0)}\cup C_2^{(0,2)}\subset\A_2^3,\\
 \end{array}\]
since $V(x_1^{(0)},x_2^{(0)},z^{(0)})\subset V(x_1^{(0)},z^{(0)})$. And
\[(\pi_3^{-1}(X_{Sing}))_{red}=V(x_1^{(0)},z^{(0)},{z^{(1)}}^3-{x_1^{(1)}}^3{x_2^{(0)}}^2)\cup V(x_2^{(0)},x_2^{(1)},z^{(0)},z^{(1)})=C_3^{(1,0)}\cup C_3^{(0,2)}
\subset\A_3^3.\]
\label{Ex1}
\end{Exam}
In this example we see how the components are defined by hyperplane coordinates for $m<3$, and at level $m=3$ the equation $f$ starts playing a role. When there are more than one approximated root, the approximated roots affect the geometry of $C_m^\nu$ one after the other as $m$ grows. This will be explained in Proposition \ref{Cgeom}. We illustrate this with another example.

\begin{Exam}
Consider the quasi-ordinary surface $f=(z^2-x_1^3)^3-x_1^{10}x_2^4$. The generators of the semigroup are $\gamma_1=(\frac{3}{2},0)$ and $\gamma_2=(\frac{10}{3},
\frac{4}{3})$, and the singular locus is $X_{Sing}=\{x_1=z=0\}\cup\{x_2=z^2-x_1^3=0\}$.
If we lift the component of the singular locus $Z_2=\{x_2=f_1=0\}$ at level $3$, we have that
$(\pi_3^{-1}(Z_2))_{red}=V(x_2^{(0)},F_1^{(0)},F_2^{(1)},F_2^{(2)},F_2^{(3)})$, where $F_1^{(0)}={z^{(0)}}^2-{x_1^{(0)}}^3$, and
\[\begin{array}{ll}
F_2^{(1)} & \equiv 3{F_1^{(0)}}^2F_1^{(1)}\mbox{ mod }(x_2^{(0)},z^{(0)})\\

 & \equiv 0\mbox{ mod }(x_2^{(0)},z^{(0)},F_1^{(0)})\\
\\
F_2^{(2)} & \equiv 3{F_1^{(0)}}^2F_1^{(2)}+6F_1^{(0)}{F_1^{(1)}}^2\mbox{ mod }(x_2^{(0)},z^{(0)})\\

 & \equiv 0\mbox{ mod }(x_2^{(0)},z^{(0)},F_1^{(0)})\\
 \\
 F_2^{(3)} & \equiv {F_1^{(1)}}^3\mbox{ mod }(x_2^{(0)},z^{(0)},F_1^{(0)}),\\
 \end{array}\]
and then $(\pi_3^{-1}(Z_2))_{red}=V(x_2^{(0)},F_1^{(0)},F_1^{(1)})$. Notice that it is not a component of $(\pi_3^{-1}(X_{Sing}))_{red},$ since it is not irreducible. Indeed, it decomposes as
\[(\pi_3^{-1}(Z_2))_{red}=V(x_1^{(0)},x_2^{(0)},z^{(0)})\cup\overline{V(x_2^{(0)},F_1^{(0)},F_1^{(1)})\cap D(x_1^{(0)})}.\]
\label{Ex2}
\end{Exam}

We have seen in this example how, to give a minimal set of generators of $J_m^\nu$, we need to study the polynomials $F^{(l)}\mbox{ mod }J_{l-1}^\nu$.
Therefore we introduce the following definition.

\begin{Defi}
For $\nu\in \s_{Sing}\cap N_0$ and $0\leq l\leq m$, we denote by
\[F_{\nu}^{(l)}:=F^{(l)}\mbox{ mod }(I^\nu,J_{l-1}^\nu),\]
and for the approximated roots the notation is $F_{j,\nu}^{(l)}$.
\label{Fnu}
\end{Defi}

And now we obviously have that
\[J_m^\nu=Rad(F_\nu^{(0)},\ldots,F_\nu^{(m)}).\]

Regarding the claim in Remark \ref{RemTonto}, once we consider $F_\nu^{(l)}$, it is not true anymore that the polynomials are non-zero. But, whenever $F_\nu^{(l)}$ is non-zero, then it is quasi-homogeneous of degree $l$.

\vspace{3mm}

In general, the first approximated root which appears is not necessarily the first one, and the process does not finish with the last one.
To control, for a given $\nu$, all this behaviour, we defined the integers $i(\nu),c(\nu)$ and $m(\nu)$ in Definition \ref{defs}.
 Indeed, given $m\in\Z_{>0}$ and $\nu\in\s_{Sing}\cap N_0$ such that $l_{c(\nu)}(\nu)\leq m$, the approximated roots which will influence the defining ideal of $D_m^\nu$, are
\[f_{c(\nu)},\ldots,f_{j(m,\nu)},\]
where remember the convention $f_0=z$. Moreover, the moment when $f_i$ begins to influence the defining equations of $C_m^\nu$ (or in other words, the generators of $J_m^\nu$) for the first time is exactly at $m=l_i(\nu)$. This is the content of Corollary \ref{Corolario3}.
The meaning of the integer $i(\nu)$ is that, at $m=l_{i(\nu)}(\nu)$, $\nu$ does no longer give rise to an irreducible component (see Lemma \ref{LemD}).
The integer $j(m,\nu)$ will be useful to describe the component $C_m^\nu$ (see Proposition \ref{Cgeom}), while $j'(m,\nu)$ will be crucial when studying the inclusion $C_m^{\nu'}\subseteq C_m^\nu$ (see Proposition \ref{PropC1}).

\begin{Exam}
We revisit Example \ref{Ex1}. If we lift the component $C_3^{(0,2)}$ to level $4$, we have
\[\pi_{4,3}^{-1}(C_3^{(0,2)})=V(x_2^{(0)},x_2^{(1)},z^{(0)},z^{(1)},F^{(4)})\]
where $F^{(4)}\equiv {x_1^{(0)}}^3{x_2^{(2)}}^2\mbox{ mod }(x_2^{(0)},x_2^{(1)},z^{(0)},z^{(1)})$. Therefore
\[\pi_{4,3}^{-1}(C_3^{(0,2)})=C_4^{(1,2)}\cup C_4^{(0,3)}.\]
Then at level $m=4$, the vector $(0,2)$ does not give rise to an irreducible component any longer. The reason is that $(0,2)\notin N_1$ and $4=l_1(0,2)$.
\label{Ex1re}
\end{Exam}

This is the case in general as we claim in the next Lemma, whose proof is left to Section \ref{Proofs}.

\begin{Lem}
For $m\in\Z_{>0}$ and $\nu\in[0,m]^2\cap\s_{Sing}\cap N_0$, we have
\[D_m^\nu=\emptyset\mbox{ if and only if }m\geq l_{i(\nu)}(\nu).\]
\label{LemD}
\end{Lem}

As a consequence of this Lemma, we are going to prove an improvement of Lemma \ref{LemGamma}, namely, for $m\in\Z_{>0}$, to cover $\pi_m^{-1}(X_{Sing})$ it is enough to consider $\nu\in[0,m]^2\cap\s_{Sing}\cap N_0$ with $m<l_{i(\nu)}(\nu)$.

Notice that, if $l_{m(\nu)}(\nu)\leq m$, then $l_1(\nu)\leq l_{m(\nu)}(\nu)\leq m$, and by the previous Lemma, we have to ask $\nu\in N_1$ whenever $l_{m(\nu)}\leq m$.
\begin{Defi}
 Given $m\in\Z_{> 0}$ we define the sets:
\[\begin{array}{l}
H_m=\{\nu\in[0,m]^2\cap\s_{Sing}\cap N_0\ |\ l_{m(\nu)}(\nu)\geq m+1\},\\
\\
L_m=\{\nu\in[0,m]^2\cap\s_{Sing}\cap N_1\ |\ l_{m(\nu)}(\nu)\leq m<l_{i(\nu)}(\nu)\}.\\
\end{array}\]
It will be necessary later to subdivide the set $L_m$ as
\[\begin{array}{l}
L_m^==\{\nu\in L_m\ |\ l_{m(\nu)}(\nu)\leq m<\mbox{min}\{l_{m(\nu)}(\nu)+e_{m(\nu)},l_{i(\nu)}(\nu)\}\},\\
\\
L_m^<=\{\nu\in L_m\ |\ l_{m(\nu)}(\nu)+e_{m(\nu)}\leq m<l_{i(\nu)}(\nu)\}.\\
\end{array}\]
\label{DefsSets}
\end{Defi}

If we come back to Example \ref{Ex1}, we have that
\[\begin{array}{ll}
(\pi_1^{-1}(X_{Sing}))_{red}=C_1^{(1,0)}\cup C_1^{(0,1)} & \mbox{ with }(1,0),(0,1)\in H_1,\\
\\
(\pi_2^{-1}(X_{Sing}))_{red}=C_2^{(1,0)}\cup C_2^{(0,2)} & \mbox{ with }(1,0),(0,2)\in H_2,\\
\\
(\pi_3^{-1}(X_{Sing}))_{red}=C_3^{(1,0)}\cup C_3^{(0,2)} & \mbox{ with }(1,0)\in L_3^=\mbox{ and }(0,2)\in H_3.\\
\end{array}\]

\vspace{3mm}

\begin{Lem}
For $m\in\Z_{>0}$, we have that $H_m\cup L_m\neq \emptyset$, and
\[\pi_m^{-1}(X_{Sing})=\bigcup_{\nu\in H_m\cup L_m}C_m^\nu.\]
\label{LemTh}
\end{Lem}

{\em Proof.} The first claim follows because $(m,m)\in H_m$ for any $m\in\Z_{>0}$. Indeed, since $l_{m(\nu)}(\nu)\geq l_1(\nu)=n_1e_1\langle\nu,\gamma_1\rangle=e_1m(a_1+b_1)>m$, where the last inequality holds because $a_1+b_1>1$, since the branch is normalized.
 By Lemma \ref{LemGamma}
\[\bigcup_{\nu\in H_m\cup L_m}C_m^\nu\subseteq\bigcup_{\nu\in[0,m]^2\cap\s_{Sing}\cap N_0}C_m^\nu=\pi_m^{-1}(X_{Sing}).\]
We prove the other inclusion. Notice that $\nu\notin H_m\cup L_m$ is equivalent to $l_{i(\nu)}(\nu)\leq m$, because $m(\nu)\leq i(\nu)$. For any $\gamma\in\pi_m^{-1}(X_{Sing})$,

$\bullet$ if $x_i\circ\gamma\neq 0$ for $i=1,2$, then $\nu:=(\mbox{ord}_t(x_1\circ\gamma),\mbox{ord}_t(x_2\circ\gamma))\in[0,m]^2$ and $\gamma\in D_m^\nu$. Hence, by Lemma \ref{LemD} we have that $m<l_{i(\nu)}(\nu)$ and therefore $\nu\in H_m\cup L_m$.

Otherwise,

$\bullet$ if $x_i\circ\gamma=0$ for $i=1,2$ we saw in the proof of Lemma \ref{LemGamma} that $\gamma\in C_m^\nu$ for any $\nu\in[0,m]^2\cap\s_{Sing}\cap N_0$. Therefore, since $H_m\cap L_m\neq\emptyset$, there exists $\nu$ with $\gamma\in C_m^\nu$.

$\bullet$ if $x_1\circ\gamma=0$ and $x_2\circ\gamma\neq 0$, then by the proof of Lemma \ref{LemGamma}, we have that $\gamma\in C_m^{(m,\alpha)}$, where $\alpha=\mbox{ord}_t(x_2\circ\gamma)$. We have to prove that $\nu:=(m,\alpha)\in H_m$, and this follows since $l_{m(\nu)}(\nu)\geq l_1(\nu)=e_1(a_1m+b_1\alpha)>m$, again using that the branch is normalized.

$\bullet$ if $x_1\circ\gamma\neq 0$ and $x_2\circ\gamma=0$, then by the proof of Lemma \ref{LemGamma}, $\gamma\in C_m^\nu$ for any $\nu$ with $\nu_1\leq\alpha$ and $\nu_2\leq m$, where $\alpha=\mbox{ord}_t(x_1\circ\gamma)$. Hence we only have to prove that $([0,\alpha]\times[0,m])\cap(H_m\cup L_m)\neq\emptyset$. If $b_1\geq 1$ then $\nu:=(\alpha,m)\in H_m$. Indeed, if $b_1>1$ clearly $l_{m(\nu)}(\nu)\geq l_1(\nu)=e_1(a_1\alpha+b_1m)>m$. The same works if $b_1=1$ and $g>1$, because then $e_1>1$. If $b_1=1$ and $g=1$, then by Theorem \ref{CorSing} we deduce that $\r_2\nsubseteq\s_{Sing}$. Then $\alpha>0$ and it follows that $\nu\in H_m$.
 The case left is $b_1=0$. In this case, if we set $\nu:=(\alpha,m)$, we have $m(\nu)=g_1$ and $i(\nu)\geq g_1+1$. Then $l_{i(\nu)}(\nu)\geq l_{g_1+1}(\nu)=e_{g_1+1}(\alpha n_{g_1+1}\gamma_{g_1+1}^{(1)}+mn_{g_1+1}\gamma_{g_1+1}^{(2)})>m$, where we are using that if $\alpha=0$ then $\r_2\subseteq\s_{Sing}$ and $n_{g_1+1}\gamma_{g_1+1}^{(2)}>1$. If $\alpha\geq n_1$, the same argument shows that $(n_1,m)\in L_m$. Otherwise $(\alpha,m)\in L_m$, since $\gamma\in X_m$ and $\gamma_1=(\frac{a_1}{n_1},0)$, therefore $a_1\mbox{ord}_t(x_1\circ\gamma)=n_1\mbox{ord}_t(z\circ\gamma)$, which implies that $\alpha\frac{a_1}{n_1}\in\Z$ or in other words, $(\alpha,m)\in N_1$.
\hfill $\Box$

\vspace{3mm}

Given $\nu\in\s_{Sing}\cap N_0$, it gives rise to a candidate of irreducible component at level $m$, $C_m^\nu$, for
\[0=l_{c(\nu)}(\nu)\leq m<l_{i(\nu)}(\nu).\]

\begin{Rem}
\begin{enumerate}
\item[(i)] If $\nu\in H_m$, then $c(\nu)=0$, since otherwise $m(\nu)=c(\nu)$ and $l_{m(\nu)}(\nu)=0<m+1$, or in other words, if $g_1>0$ then $H_m\cap\r_2=\emptyset$.
\item[(ii)]  It is clear that
$j(m,\nu)=0  \mbox{ if and only if }\nu\in H_m$ and $j(m,\nu)\geq 1 \mbox{ if and only if }\nu\in L_m$.
\end{enumerate}
\label{RemDefHL}
\end{Rem}

For $\nu\in H_m$ the sets $C_m^\nu$ are very easy to describe, as we see in the next Proposition.

\begin{Pro}
For $m\in\Z_{>0}$ and $\nu\in H_m$ we have that
\[J_m^\nu=(z^{(0)},\ldots,z^{([m/n])}),\]
 and hence the set $C_m^\nu$ is defined by hyperplane coordinates in $\A_m^3$ as
\[C_m^\nu=V(x_i^{(0)},\ldots,x_i^{(\nu_i-1)},i=1,2;z^{(0)},\ldots,z^{([\frac{m}{n}])}).\]
\label{CHm}
\end{Pro}
{\em Proof.}
The proof is by induction on $m$. For $m=1$ we have

\[H_1=\left\{\begin{array}{ll}
\{(1,1)\} & \mbox{ if }\gamma_1=(\frac{1}{n},\frac{1}{n})\mbox{ and }g=1\\
\\
\{(1,0),(1,1)\} & \mbox{ if }g_1>0\mbox{ (recall that the branch is normalized)}\\
\\
\{(1,0),(1,1),(0,1)\} &\mbox{ otherwise}\\
\end{array}\right.\]
and the claim follows, since
\[C_1^{(1,0)}=V(x_1^{(0)},z^{(0)}),\ C_1^{(1,1)}=V(x_1^{(0)},x_2^{(0)},z^{(0)})\mbox{ and }\ C_1^{(0,1)}=V(x_2^{(0)},z^{(0)})\mbox{ if }g_1=0.\]
Suppose the claim is true for $m$ and we will prove it for $m+1$. Given $\nu\in H_{m+1}$, since $\nu\in H_m$, by induction hypothesis we have that
\[J_{m+1}^\nu=(z^{(0)},\ldots,z^{([m/n])},F_\nu^{(m+1)}).\]
We claim that
\[F_\nu^{(m+1)}=\left\{\begin{array}{cl}
0 & \mbox{ if }m+1\not\equiv 0\mbox{ mod }n\\
\\
{z^{(\frac{m+1}{n})}}^n & \mbox{ otherwise}\\
\end{array}\right.\]
which proves the result, since
\[\left[\frac{m+1}{n}\right]=\left\{\begin{array}{cl}
[\frac{m}{n}] & \mbox{ if }m+1\not\equiv 0\mbox{ mod }n\\
\\
\frac{m+1}{n} & \mbox{ otherwise}\\
\end{array}\right.\]
If $F_\nu^{(m+1)}\neq 0$ then it is a quasi-homogeneous polynomial of degree $m+1$. By the expansion (\ref{grgr1}) given in Lemma \ref{LemExpSR}
\[f=z^n+\sum c_{\underline{\alpha},\underline{r}}x_1^{\alpha_1}x_2^{\alpha_2}z^{r_1}\]
where $(\alpha_1,\alpha_2)+r_1\gamma_1>n_1e_1\gamma_1=n\gamma_1$. Then any monomial in $f-z^n$ verifies for any $\gamma\in D_m^\nu$
\[\begin{array}{ll}
\mbox{ord}_t(c_{\underline{\alpha},\underline{r}}x_1^{\alpha_1}x_2^{\alpha_2}\circ\gamma) & =\langle\nu,(\alpha_1,\alpha_2)\rangle+r_1\mbox{ord}_t(z\circ\gamma)\\
 & \geq\langle\nu,(\alpha_1,\alpha_2)\rangle+r_1\frac{m+1}{n}\\
 & \geq n\langle\nu,\gamma_1\rangle-r_1\langle\nu,\gamma_1\rangle+r_1\frac{m+1}{n}\\
 & =(n-r_1)\frac{l_1(\nu)}{n}+r_1\frac{m+1}{n}>(n-r_1)\frac{m+1}{n}+r_1\frac{m+1}{n}=m+1\\
 \end{array}\]
because $\nu\in H_{m+1}$ and by induction hypothesis ord$_t(z\circ\gamma)>[\frac{m}{n}]$ (and therefore $\geq\frac{m+1}{n}$). Hence these monomials do not contribute to $F_\nu^{(m+1)}$ and the result follows by the quasi-homogeneity of $F_\nu^{(m+1)}$.\hfill$\Box$

\vspace{3mm}

For $\nu\in L_m$ the geometry of $C_m^\nu$ is much more complicated, the ideal $J_m^\nu$ is described in Corollary \ref{Corolario3}.
In the next Proposition we compare jet schemes of a quasi-ordinary singularity with jet schemes of its approximated roots.
For $1\leq i\leq g$, $m\in\Z_{>0}$ and $\nu\in H_m\cup L_m$, we denote by $D_{i,m}^\nu$ the set
\[D_{i,m}^\nu=\{\gamma\in X_m^{(i)}\ |\ \mbox{ord}_t(x_k\circ\gamma)=\nu_k, k=1,2\}_{red},\]
and we denote $D_{g,m}^\nu$ simply by $D_m^\nu$ (see
Definition \ref{SemiX} for the definition of $X^{(i)}$).

\begin{Pro}
For $m\in\Z_{>0}$ and $\nu\in H_m\cup L_m$, we have that
\[D_m^\nu=(\pi_{m,[\frac{m}{e_j}]}^a)^{-1}(D_{j,[\frac{m}{e_j}]}^\nu)\]
where $j=j(m,\nu)$, and for $q>p$, $\pi_{q,p}^a:\A^3_q\longrightarrow\A^3_p$ is the projection on the jet schemes of the affine ambient space.
\label{Cgeom}
\end{Pro}

Hence, for $m\in\Z_{>0}$ and $\nu\in L_m$ with $j(m,\nu)=j$, the geometry of $C_m^\nu$ is determined
by the geometry of the $j$-th approximated root.

\vspace{2mm}

Before proving the Proposition we need the following technical result, whose proof is moved to Section \ref{Proofs}.

\begin{Lem}
For $m\in\Z_{>0}$ and $\nu\in H_m\cup L_m$, we have that for all $\gamma\in D_m^\nu$,
\[\begin{array}{l}
\mbox{ord}_t(f_k\circ\gamma)=\langle\nu,\gamma_{k+1}\rangle,\ \mbox{ for }0\leq k\leq j(m,\nu)-1\\
\\
\mbox{ord}_t(f_k\circ\gamma)>\frac{m}{e_k},\ \mbox{ for }j(m,\nu)\leq k\leq g\\
\end{array}\]
\label{Lemfk}
\end{Lem}

{\bf\em Proof of Proposition \ref{Cgeom}.} For $\nu\in H_m$ we have $j(m,\nu)=0$, and the claim follows by Proposition \ref{CHm}. For $\nu\in L_m$ it is enough to prove that, if $j(m,\nu)=j$ we have
\begin{equation}
D_m^\nu=\left\{\gamma\in\A_m^3\ |\ \mbox{ord}_t(x_i\circ\gamma)=\nu_i, i=1,2\mbox{ and }\mbox{ord}_t(f_j\circ\gamma)>\frac{m}{e_j}\right\}.
\label{Claim2}
\end{equation}
By Lemma \ref{Lemfk} it follows that
\[D_m^\nu\subseteq\left\{\gamma\in\A_m^3\ |\ \mbox{ord}_t(x_i\circ\gamma)=\nu_i,i=1,2\mbox{ and ord}_t(f_j\circ\gamma)>\frac{m}{e_j}\right\}\]
We prove the other inclusion. Let $\gamma$ be a jet with ord$_t(x_i\circ\gamma)=\nu_i$ for $i=1,2$ and ord$_t(f_j\circ\gamma)>\frac{m}{e_j}$. We want to prove that it is indeed an $m$-jet in $X$, or in other words, that ord$_t(f\circ\gamma)\geq m+1$. Notice that if $j=g$ there is nothing to prove.
Then $j<g$, and first we will prove that
\[\mbox{ord}_t(f_{j+1}\circ\gamma)>\frac{m}{e_{j+1}}.\]

Indeed, consider $f_j$ quasi-ordinary surface with $j$ characteristic exponents. If $\bar m:=[\frac{m}{e_j}]$ and $\bar\gamma=\pi_{m,\bar m}(\gamma)$, then we have that $\bar\gamma\in D_{j,\bar m}^\nu$. Moreover $n_j\langle\nu,\gamma_j\rangle\leq \bar m$, and then, by Lemma \ref{Lemfk} applied to $f_j$, we have
\[\mbox{ord}_t(f_k\circ\bar\gamma)=\langle\nu,\gamma_{k+1}\rangle,\ \mbox{ for }0\leq k\leq j-1.\]
Since $\langle\nu,\gamma_{k+1}\rangle\leq\langle\nu,\gamma_j\rangle<\bar m<m$, we deduce that ord$_t(f_k\circ\gamma)=\mbox{ord}_t(f_k\circ\bar\gamma)$.

Now we consider $f_{j+1}$. By Lemma \ref{Lema35} we have
\[f_{j+1}=f_j^{n_{j+1}}-c_{j+1}x_1^{\alpha_1^{(j+1)}}x_2^{\alpha_2^{(j+1)}}z^{r_1^{(j+1)}}\cdots f_{j-1}^{r_j^{(j+1)}}+\sum c_{\underline{\alpha},\underline{r}}x_1^{\alpha_1}x_2^{\alpha_2}z^{r_1}\cdots f_j^{r_{j+1}}\]
and using that ord$_t(f_k\circ\gamma)=\langle\nu,\gamma_{k+1}\rangle$ for $0\leq k\leq j-1$ we have
\[\begin{array}{rl}
\mbox{ord}_t(f_j^{n_{j+1}}\circ\gamma) & =n_{j+1}\mbox{ord}_t(f_j\circ\gamma)>\frac{m}{e_{j+1}}\\
\\
\mbox{ord}_t((x_1^{\alpha_1^{(j+1)}}x_2^{\alpha_2^{(j+1)}}z^{r_1^{(j+1)}}\cdots f_{j-1}^{r_j^{(j+1)}})\circ\gamma) & =\langle\nu,(\alpha_1^{(j+1)},\alpha_2^{(j+1)})+r_1^{(j+1)}\gamma_1+\cdots+r_j^{(j+1)}\gamma_j\rangle\\
 & =n_{j+1}\langle\nu,\gamma_{j+1}\rangle >\frac{m}{e_{j+1}}\\
\\
\mbox{ord}_t((c_{\underline{\alpha},\underline{r}}x_1^{\alpha_1}x_2^{\alpha_2}z^{r_1}\cdots f_j^{r_{j+1}})\circ\gamma) & =\langle\nu,(\alpha_1,\alpha_2)+r_1\gamma_1+\cdots+r_j\gamma_j\rangle+r_{j+1}\mbox{ord}_t(f_j\circ\gamma)\\
\end{array}\]
If ord$_t(f_{j+1}\circ\gamma)\leq\frac{m}{e_{j+1}}$, then there must exist $c_{\underline{\alpha},\underline{r}}\neq 0$ such that
\[\mbox{ord}_t(f_{j+1}\circ\gamma)=\langle\nu,(\alpha_1,\alpha_2)+r_1\gamma_1+\cdots+r_j\gamma_j\rangle+r_{j+1}\mbox{ord}_t(f_j\circ\gamma)\leq\frac{m}{e_{j+1}}\]
and we get the following inequalities
\[\begin{array}{rl}
(n_{j+1}-r_{j+1})\langle\nu,\gamma_{j+1}\rangle+r_{j+1}\frac{m}{e_j} & \leq\langle\nu,(\alpha_1,\alpha_2)+r_1\gamma_1+\cdots+r_j\gamma_j\rangle+r_{j+1}\frac{m}{e_j}\\
\\
 & <\langle\nu,(\alpha_1,\alpha_2)+r_1\gamma_1+\cdots+r_j\gamma_j\rangle+r_{j+1}\mbox{ord}_t(f_j\circ\gamma)\\
\\
 & \leq\frac{m}{e_{j+1}}\\
 \end{array}\]
Hence $(n_{j+1}-r_{j+1})\langle\nu,\gamma_{j+1}\rangle+r_{j+1}\frac{m}{e_j}<\frac{m}{e_{j+1}}=n_{j+1}\frac{m}{e_j}$, and then  $(n_{j+1}-r_{j+1})\langle\nu,\gamma_{j+1}\rangle<(n_{j+1}-r_{j+1})\frac{m}{e_j}$. Since $r_{j+1}<n_{j+1}$, we have
\[\langle\nu,\gamma_{j+1}\rangle<\frac{m}{e_j}\]
which is in contradiction with $j(m,\nu)=j$. Therefore we have just proved that ord$_t(f_{j+1}\circ\gamma)>\frac{m}{e_{j+1}}$.

 To finish, by Lemma \ref{LemExpSR} we have
\[f=f_{j+1}^{e_{j+1}}+\sum_{(i_1,i_2)+k\gamma_1>n_{j+1}e_{j+1}\gamma_{j+1}} c_{i_1i_2k}x_1^{i_1}x_2^{i_2}z^k,\]
and
\[\begin{array}{rl}
\mbox{ord}_t(f_{j+1}^{e_{j+1}}\circ\gamma) & =e_{j+1}\mbox{ord}_t(f_{j+1}\circ\gamma)>m\\
\\
\mbox{ord}_t((c_{i_1i_2k}x_1^{i_1}x_2^{i_2}z^k)\circ\gamma) & \geq n_{j+1}e_{j+1}\langle\nu,\gamma_{j+1}\rangle=l_{j+1}(\nu)>m\\
\end{array}\]
Hence ord$_t(f\circ\gamma)>m$ as we wanted to prove. \hfill$\Box$

\vspace{3mm}

As a consequence of Proposition \ref{Cgeom}, we have the following algebraic counterpart, where we explain how the equations of the approximated roots appear as generators of $J_m^\nu$, and therefore a minimal presentation of the ideal $J_m^\nu$ is given.

\begin{Cor} Given $m\in\Z_{>0}$ and $\nu\in H_m\cup L_m$, for $l_i(\nu)\leq l<l_{i+1}(\nu)$ (resp. $l_i(\nu)\leq l\leq m$) if $c(\nu)\leq i< j(m,\nu)$ (resp. $i=j(m,\nu)$), we have that
\[F_\nu^{(l)}=\left\{\begin{array}{cl}
{F_{i,\nu}^{(\frac{l}{e_i})}}^{e_i}  &  \mbox{ if }l\equiv 0\mbox{ mod }e_i\\
\\
0 & \mbox{ otherwise}\\
\end{array}\right.\]
Hence the ideal $J_m^\nu$ is generated by the polynomials
\[J_m^\nu=\left(F_{i,\nu}^{(\frac{l_i(\nu)}{e_i}+k_i)}\right)_i\]
for $c(\nu)\leq i\leq j(m,\nu)$ such that $l_i(\nu)<l_{i+1}(\nu)$, and $0\leq k_i<\frac{l_{i+1}(\nu)-l_i(\nu)}{e_i}\mbox{ if }i<j(m,\nu)$ and $0\leq k_{j(m,\nu)}\leq
[\frac{m-l_{j(m,\nu)}(\nu)}{e_{j(m,\nu)}}]$. Moreover,  we have
\begin{equation}
F_{i,\nu}^{(\frac{l_i(\nu)}{e_i})}={F_{i-1,\nu}^{(\frac{l_i(\nu)}{e_{i-1}})}}^{n_i}-c_i{x_1^{(\nu_1)}}^{\alpha_1^{(i)}}{x_2^{(\nu_2)}}^{\alpha_2^{(i)}}
{z^{(\langle\nu,\gamma_1\rangle)}}^{r_1^{(i)}}\cdots {F_{i-2,\nu}^{(\frac{l_{i-1}(\nu)}{e_{i-2}})}}^{r_{i-1}^{(i)}}+G_{i,\nu},
\label{EqP}
\end{equation}
where $G_{i,\nu}$ is the polynomial
\[G_{i,\nu}=\sum c_{\underline{\alpha},\underline{r}}{x_1^{(\nu_1)}}^{\alpha_1}{x_2^{(\nu_2)}}^{\alpha_2}
{z^{(\langle\nu,\gamma_1\rangle)}}^{r_1}\cdots {F_{i-2,\nu}^{(\frac{l_{i-1}(\nu)}{e_{i-2}})}}^{r_{i-1}}{F_{i-1,\nu}^{(\frac{l_i(\nu)}{e_{i-1}})}}^{r_i}\]
with $c_{\underline{\alpha},\underline{r}}$ are the coefficients appearing in the expansion given in Lemma \ref{Lema35}, and such that
\vspace{2mm}

\begin{enumerate}
\item[(i)] $c_{\underline{\alpha},\underline{r}}\neq 0$

\item[(ii)] $\langle\nu,(\alpha_1,\alpha_2)+r_1\gamma_1+\cdots+r_i\gamma_i\rangle=n_i\langle\nu,\gamma_i\rangle$
\end{enumerate}
\vspace{2mm}
Notice that condition (ii) does not hold when $\nu\notin\r_1\cup\r_2$, and hence $G_{i,\nu}=0$ in these cases.
\label{Corolario3}
\end{Cor}

{\em Proof.} It is a consequence of Proposition \ref{Cgeom}, applied to any $l_i(\nu)\leq l<l_{i+1}(\nu)$, and using the trivial observation that $J_{m'}^\nu\subseteq J_m^\nu$ for any $m'<m$. Indeed, for any $l_i(\nu)\leq l< l_{i+1}(\nu)$ we study the polynomials $F_\nu^{(l_i(\nu))},\ldots,F_\nu^{(l)}$ (note that we need $l_i(\nu)<l_{i+1}(\nu)$). We have that $j(l,\nu)=i$, and, by Proposition \ref{Cgeom}, $D_l^\nu=(\pi_{l,[\frac{l}{e_i}]}^a)^{-1}(D_{i,[\frac{l}{e_i}]}^\nu)$, or in other words
\[D_l^\nu=\left\{\gamma\in\A_l^3\ |\ \mbox{ord}_t(x_k\circ\gamma)=\nu_k,\ k=1,2,\mbox{ ord}_t(f_i\circ\gamma)>\left[\frac{l}{e_i}\right]\right\}\]
Then the ideal $J_l^\nu$ only depends on $f_i$ (and hence on its approximated roots). Moreover, by Lemma \ref{Lemfk}, we deduce that for $0\leq k<i$
\[F_{k,\nu}^{(r_k)}\in J_l^\nu,\ \mbox{ for }0\leq r_k<\langle\nu,\gamma_{k+1}\rangle.\]
By the expansion (\ref{grgr1}) in Lemma \ref{LemExpSR}
\[f=f_i^{e_i}+\sum c_{\underline{\alpha},\underline{r}}x_1^{\alpha_1}x_2^{\alpha_2}z^{r_1}\cdots f_i^{r_{i+1}}\]
where $(\alpha_1,\alpha_2)+r_1\gamma_1+\cdots+r_{i+1}\gamma_{i+1}>n_{i+1}e_{i+1}\gamma_{i+1}$.

The part $f_i^{e_i}$ contributes to $F^{(l)}$ with ${F_{i,\nu}^{(\frac{l}{e_i})}}^{e_i}$ and only when $l$ is divisible by $e_i$. While for the {\em monomials}
$x_1^{\alpha_1}x_2^{\alpha_2}z^{r_1}\cdots f_i^{r_{i+1}}$ the contribution is given by
\[{x_1^{(\nu_1)}}^{\alpha_1}{x_2^{(\nu_2)}}^{\alpha_2}{z^{(\langle\nu,\gamma_1\rangle)}}^{r_1}\cdots {F_{i-1,\nu}^{(\frac{l_i(\nu)}{e_{i-1}})}}^{r_i}{F_{i,\nu}^{(\frac{a}{r_{i+1}})}}^{r_{i+1}}\]
with $\langle\nu,(\alpha_1,\alpha_2)+r_1\gamma_1+\cdots+r_i\gamma_i\rangle+a=l$, and only when $a$ is divisible by $r_{i+1}$. Set $a=br_{i+1}$, we claim that $F_{i,\nu}^{(b)}$ belongs to $J_l^\nu$ and hence the {\em monomial} does not contribute to $F^{(l)}$. Indeed, at level $l'=be_i$ it appears as ${F_{i,\nu}^{(b)}}^{e_i}$, and we only need to prove that $l'<l$. Since
\[\begin{array}{ll}
l & =\langle\nu,(\alpha_1,\alpha_2)+r_1\gamma_1+\cdots+r_i\gamma_i\rangle+r_{i+1}b\\
 & =\langle\nu,(\alpha_1,\alpha_2)+r_1\gamma_1+\cdots+r_i\gamma_i\rangle+r_{i+1}\frac{l'}{e_i}\\
 \end{array}\]
we have
 \[\begin{array}{ll}
 l' & =\frac{e_i}{r_{i+1}}(l-\langle\nu,(\alpha_1,\alpha_2)+r_1\gamma_1+\cdots+r_i\gamma_i\rangle)\\
  & \leq\frac{e_i}{r_{i+1}}(l-l_{i+1}(\nu)+r_{i+1}\langle\nu,\gamma_{i+1}\rangle)\\
 & =\frac{e_i}{r_{i+1}}l-(\frac{e_i}{r_{i+1}}-1)l_{i+1}(\nu)\\
 \end{array}\]
 Therefore $r_{i+1}l'\leq e_il-(e_i-r_{i+1})l_{i+1}(\nu)$. Suppose that $l'\geq l$, then
 \[(e_i-r_{i+1})l_{i+1}(\nu)\leq e_il-r_{i+1}l'\leq (e_i-r_{i+1})l\]
 and since $e_i>r_{i+1}$ it contradicts the fact that $j(l,\nu)=i$.

Now equation (\ref{EqP}) follows by Lemma \ref{Lema35}.
\hfill$\Box$

\vspace{3mm}

From Corollary \ref{Corolario3} we deduce the following:
\begin{equation}
D_m^\nu\subset D(F_{0,\nu}^{(\frac{l_1(\nu)}{e_0})}\cdots F_{j-1,\nu}^{(\frac{l_j(\nu)}{e_{j-1}})}),\ \mbox{ for any }\nu\in L_m\mbox{ with }j(m,\nu)=j.
\label{eqO}
\end{equation}

\vspace{3mm}

To illustrate the description of $J_m^\nu$ given in Corollary \ref{Corolario3}, we consider some particular cases.

\vspace{2mm}

$\bullet$ First, the simplest case, when $\nu\notin\r_1\cup\r_2$ we have $G_{j,\nu}=0$ for any $j$, and $m(\nu)=1$. Hence (note that we also have $c(\nu)=0$):

\begin{equation}
    \begin{array}{lcl}
F_{0,\nu}^{(r)} & = &    z^{(r)}, \mbox{ for }0\leq r<\frac{l_{1}(\nu)}{e_0}\\
\\
F_{1,\nu}^{(\frac{l_1(\nu)}{e_1})} & = & {z^{(\langle\nu,\gamma_1\rangle)}}^{n_1}-{x_1^{(\nu_1)}}^{a_1}\\
\\
F_1^{(\frac{l_1(\nu)}{e_1}+r)} & &\mbox{ for }1\leq r<\frac{l_2(\nu)-l_1(\nu)}{e_1}\\

\ \ \ \ \ \vdots  & \vdots & \\

F_{g_1-1,\nu}^{(\frac{l_{g_1-1}(\nu)}{e_{g-1}})} & = & {F_{g_1-2,\nu}}^{(\frac{l_{g_1-1}(\nu)}{e_{g_1-2}})}-{x_1^{(\nu_1)}}^{n_{g_1-1}}{z^{(\langle\nu,\gamma_1\rangle)}}^{r_1^{(g_1-1)}}\cdots
                                    {F_{g_1-3,\nu}^{(\frac{l_{g_1-2}(\nu)}{e_{g_1-3}})}}^{r_{g_1-2}^{(g_1-1)}}\\
\\
F_{g_1-1,\nu}^{(\frac{l_{g_1-1}(\nu)}{e_{g-1}}+r)} & & \mbox{ for }1\leq r<\frac{l_{g_1}(\nu)-l_{g_1-1}(\nu)}{e_{g_1-1}}\\
\\
F_{g_1,\nu}^{(\frac{l_{g_1}(\nu)}{e_{g_1}})} & = & {F_{g_1-1,\nu}^{(\frac{l_{g_1}(\nu)}{e_{g_1-1}})}}^{n_{g_1}}-{x_1^{(\nu_1)}}^{\alpha_1^{(g_1)}}{z^{(\langle\nu,
\gamma_1\rangle)}}^{r_1^{(g_1)}}{F_{1,\nu}^{(\frac{l_2(\nu)}{e_1})}}^{r_2^{(g_1)}} \cdots {F_{g_1-2,\nu}^{(\frac{l_{g_1-1}(\nu)}{e_{g_1-2}})}}^{r_{g_1-1}^{(g_1)}}\\
\\
F_{g_1,\nu}^{(\frac{l_{g_1}(\nu)}{e_{g_1}}+r)} & &\mbox{ for }1\leq r<\frac{l_{g_1+1}(\nu)-l_{g_1}(\nu)}{e_{g_1}}\\
\\
F_{g_1+1,\nu}^{(\frac{l_{g_1+1}(\nu)}{e_{g_1+1}})} & = & {F_{g_1,\nu}^{(\frac{l_{g_1+1}(\nu)}{e_{g_1}})}}^{n_{g_1+1}}-{x_1^{(\nu_1)}}^{\alpha_1^{(g_1+1)}}
{x_2^{(\nu_2)}}^{\alpha_2^{(g_1+1)}}
{z^{(\langle\nu,\gamma_1\rangle)}}^{r_1^{(g_1+1)}}\cdots {F_{g_1-1,\nu}^{(\frac{l_{g_1}(\nu)}{e_{g_1-1}})}}^{r_{g_1}^{(g_1+1)}}\\

\ \ \ \ \ \vdots & \vdots &\\
 F_{j,\nu}^{(\frac{l_j(\nu)}{e_j})} & = & {F_{j-1,\nu}^{(\frac{l_j(\nu)}{e_{j-1}})}}^{n_j}-{x_1^{(\nu_1)}}^{\alpha_1^{(j)}}{x_2^{(\nu_2)}}^{\alpha_2^{(j)}}
 {z^{(\langle\nu,\gamma_1\rangle)}}^{r_1^{(j)}}\cdots {F_{j-2,\nu}^{(\frac{l_{j-1}(\nu)}{e_{j-2}})}}^{r_{j-1}^{(j)}}\\
\\
 F_{j,\nu}^{(\frac{l_j(\nu)}{e_j}+r)} & & \mbox{ for }1\leq r\leq[\frac{m-l_j(\nu)}{e_j}]\\
    \end{array}
    \label{ecsFs}
\end{equation}
Notice that the variable $x_2^{(\nu_2)}$ appears for the first time in the equation $F_{g_1+1}^{(\frac{l_{g_1+1}(\nu)}{e_{g_1+1}})}$, and raised to the power one or bigger depending on whether $g_2=g_1+1$ or $g_2=g_1$ respectively.

\vspace{2mm}

$\bullet$ Another example is when $j(m,\nu)=c(\nu)$, then
\[J_m^\nu=\left(F_{c(\nu),\nu}^{(0)},\ldots,F_{c(\nu),\nu}^{([\frac{m}{e_{c(\nu)}}])}\right)\]
and if moreover $c(\nu)=0$ this is just the content of Proposition \ref{CHm}. Notice that the description of $J_m^\nu$ for $\nu\in H_m$ given in Proposition \ref{CHm} is contained in Corollary \ref{Corolario3}, but we wanted to stress the fact that for $\nu\in H_m$ the description is particularly simple.

\vspace{3mm}

Now we can prove the irreducibility of the sets $C_m^\nu$.

\begin{Pro}
 For any $m\in\Z_{> 0}$ and $\nu\in H_m\cup L_m$, the set $C_m^\nu$ is irreducible and has codimension
\[\nu_1+\nu_2+\sum_{k=0}^{j-1}\frac{l_{k+1}(\nu)-l_k(\nu)}{e_k}+\left[\frac{m-l_{j(m,\nu)}(\nu)}{e_{j(m,\nu)}}\right]+1\]
\label{Prop1}
\end{Pro}

{\em Proof.}
 To simplify notation we will denote along this proof $j(m,\nu)$ just by $j$ and by $k_i(\nu)$, or simply by $k_i$ when $\nu$ is clear from the context, we denote
the quotient $\frac{l_{i+1}(\nu)-l_i(\nu)}{e_i}$.
First notice that, by definition of $c(\nu)$, $\sum_{k=0}^{j-1}\frac{l_{k+1}(\nu)-l_k(\nu)}{e_k}=\sum_{k=c(\nu)}^{j-1}\frac{l_{k+1}(\nu)-l_k(\nu)}{e_k}$.

If $\nu\in H_m$ we have $j(m,\nu)=0$ and the claim about the codimension follows by Proposition \ref{CHm}. The closed set $C_m^\nu$ is
irreducible since it is defined by hyperplane coordinates.

 If $\nu\in L_m^=$, then by Lemma \ref{LemCm}, we have that $C_m^\nu=V(I^\nu,J_m^\nu)$, where the ideal $J_m^\nu$ is described in equation (\ref{eqJ}). Then
 Codim$(C_m^\nu)=\nu_1+\nu_2+\frac{l_1(\nu)}{n}+1$. Notice that if $c(\nu)>0$ then $l_1(\nu)=0$. The
  irreducibility of $C_m^\nu$ follows from the irreducibility of $F_{m(\nu),\nu}^{(\frac{l_{m(\nu)}(\nu)}{e_{m(\nu)}})}$.

Now let $\nu$ be an element in $L_m^<$. We have to study carefully the generators of $J_m^\nu$ given in Corollary \ref{Corolario3}. Any $F_{i,\nu}^{(l)}$ is quasi-homogeneous of degree $l$, but the second property described in Remark \ref{RemTonto} is not true anymore once we consider the equations modulo $I^\nu$. We need to know when a certain variable $x_k^{(l)}$ or $z^{(l)}$ appear for the first time in the generators of $J_m^\nu$. Notice that for any $\gamma\in D_m^\nu$ we have ord$_t(x_i\circ\gamma)=\nu_i$ for $i=1,2$, and ord$_t(z\circ\gamma)=\langle\nu,\gamma_1\rangle$. It is clear that the variables $x_k^{(\nu_k)}$ for $k=1,2$ and $z^{(\langle\nu,\gamma_1\rangle)}$ appear for the first time in the first non-monomial equation $F_{m(\nu),\nu}^{(\frac{l_{m(\nu)}(\nu)}{e_{m(\nu)}})}$. In the next equation appear $x_k^{(\nu_k+1)}$, $z^{(\langle\nu,\gamma_1\rangle+1)}$, and so on. Then looking at the generators of $J_m^\nu$ described in Corollary \ref{Corolario3}, we deduce that, for $m(\nu)\leq i<j$ and $0\leq r<k_i$, or $0\leq r<[\frac{m-l_j(\nu)}{e_j}]$ when $i=j$,
\[(\ast)\ \ \ \ \ \ \ \mbox{ the variables }x_k^{(\nu_k+k_{m(\nu)}+\cdots+k_{i-1}+r)},z^{(\langle\nu,\gamma_1\rangle+k_{m(\nu)}+\cdots+k_{i-1}+r)}\mbox{ appear for the first time in }F_{i,\nu}^{(\frac{l_i(\nu)}{e_i}+r)}\]
Notice that for $1\leq l<m(\nu)$ we have $k_l=0$.
We divide the set of generators of $J_m^\nu$, given in Corollary \ref{Corolario3}, in two sets:
\[\begin{array}{l}
\mathcal C_1=\left\{F_{i,\nu}^{(\frac{l_i(\nu)}{e_i})}\right\}_{c(\nu)\leq i\leq j,\ l_i(\nu)<l_{i+1}(\nu)}\\
\\
\mathcal C_2=\left\{F_{i,\nu}^{(\frac{l_i(\nu)}{e_i}+r)}\right\}_{(i,r)\in A_2}\\
\end{array}\]
where $A_2=\{(i,r)\ |\ c(\nu)\leq i<j,\ l_i(\nu)<l_{i+1}(\nu),\ 0<r<k_i\}\cup\{(j,r)\ |\ 0<r<[\frac{m-l_j(\nu)}{e_j}]\}$.

\vspace{2mm}

We claim:
\begin{enumerate}
\item[(i)] $V(\mathcal C_1)\simeq Z^{\Gamma_m^\nu}$, the toric variety defined by the semigroup $\Gamma_m^\nu$ generated by
\[\{\gamma_i\}_{c(\nu)\leq i\leq j(m,\nu),\ l_i(\nu)<l_{i+1}(\nu)}\]
If $\nu\notin\r_1\cup\r_2$ then $\Gamma_m^\nu=\Gamma_{j(m,\nu)}$ and $V(\mathcal C_1)$ is isomorphic to the monomial variety associated to $X^{(j(m,\nu))}$ (see Definition \ref{defMV}).

\

\item[(ii)] any $F_{i,\nu}^{(\frac{l_i(\nu)}{e_i}+r)}\in\mathcal C_2$ is linear over $D(x_1^{(\nu_1)})\cap D(x_2^{(\nu_2)})$ with respect to one of the variables described in ($\ast$), which appears for the first time on this equation.
\end{enumerate}
Since any of these equations in $\mathcal C_2$ is linear in a different variable, and, by ($\ast$) we have that it appears for the first time in $\mathcal C_2$, we deduce
\[V(I^\nu,\mathcal C_2)\cap D(x_1^{(\nu_1)})\cap D(x_2^{(\nu_2)})\simeq\A^{\alpha(m,\nu)}\]
where $\alpha(m,\nu)=3(m+1)-\nu_1-\nu_2-\sum_{i=c(\nu),...,j-1,\ l_i(\nu)<l_{i+1}(\nu)}(k_i-1)-[\frac{m-l_j(\nu)}{e_j}]$, because $V(I^\nu,\mathcal C_2)\subseteq\A_m^3\simeq\A^{3(m+1)}$. Hence
\[D_m^\nu\simeq\left(Z^{\Gamma_m^\nu}\cap D(x_1^{(\nu_1)})\cap D(x_2^{(\nu_2)})\right)\times\A^{\alpha(m,\nu)}\]
 The toric variety $Z^{\Gamma_m^\nu}$ is irreducible and hence the irreducibility of $C_m^\nu$ follows by the previous isomorphism. Moreover $Z^{\Gamma_m^\nu}$ is complete intersection, hence the codimension equals the number of defining equations, which is the cardinal of $\mathcal C_1$.
Therefore
\[\begin{array}{ll}
\mbox{Codim}(C_m^\nu) & =\sharp\mathcal C_1+\nu_1+\nu_2+\sum_{c(\nu)\leq i<j,\ l_i(\nu)<l_{i+1}(\nu)}(k_i-1)+[\frac{m-l_j(\nu)}{e_j}]\\
\\
 & =\nu_1+\nu_2+\sum_{c(\nu)\leq i<j,\ l_i(\nu)<l_{i+1}(\nu)}k_i+[\frac{m-l_j(\nu)}{e_j}]+\sharp\mathcal C_1-\sharp\{c(\nu)\leq i<j\ |\ l_i(\nu)<l_{i+1}(\nu)\}\\
 \end{array}\]
Finally the statement about the codimension follows now by these two remarks:

$\bullet$ $\sum_{c(\nu)\leq i<j,\ l_i(\nu)<l_{i+1}(\nu)}k_i=\sum_{i=c(\nu)}^{j-1}k_i$, since $k_i=0$ whenever $l_i(\nu)=l_{i+1}(\nu)$.

\vspace{2mm}

$\bullet$ $l_j(\nu)<l_{j+1}(\nu)$ by definition of $j(m,\nu)$, and therefore
\[\sharp\{c(\nu)\leq i<j\ |\ l_i(\nu)<l_{i+1}(\nu)\}=\sharp\mathcal C_1-1.\]

\vspace{2mm}

 Now we prove the claim. To prove (i), notice that we can write equation (\ref{EqP}) as
\[F_{i,\nu}^{(\frac{l_i(\nu)}{e_i})}={F_{i-1,\nu}^{(\frac{l_i(\nu)}{e_{i-1}})}}^{n_i}-c_i{x_1^{(\nu_1)}}^{\alpha_1^{(i)}}{x_2^{(\nu_2)}}^{\alpha_2^{(i)}}
{z^{(\langle\nu,\gamma_1\rangle)}}^{r_1^{(i)}}\cdots{F_{i-2,\nu}^{(\frac{l_{i-1}(\nu)}{e_{i-2}})}}^{r_{i-1}^{(i)}}\cdot U\]
where, arguing as in the proof of Lemma \ref{TechLem} (i), we have $U\neq 0$. Then we have the isomorphism
\[V(F_{i,\nu}^{(\frac{l_i(\nu)}{e_i})})_{c(\nu)\leq i\leq j(m,\nu),\ l_i(\nu)<l_{i+1}(\nu)}\simeq V(h_i)_{c(\nu)\leq i\leq j(m,\nu),\ l_i(\nu)<l_{i+1}(\nu)}\]
where $h_i=w_i^{n_i}-x_1^{\alpha_1^{(i)}}
x_2^{\alpha_2^{(i)}}z^{r_1^{(i)}}w_1^{r_2^{(i)}}\cdots w_{i-2}^{r_{i-1}^{(i)}}$, with the relation  $n_i\gamma_i=(\alpha_1^{(i)},\alpha_2^{(i)})+r_1^{(i)}\gamma_1+\cdots+r_{i-1}^{(i)}\gamma_{i-1}$. And $V(h_i)_{c(\nu)\leq i\leq j(m,\nu),\ l_i(\nu)<l_{i+1}(\nu)}$ is  isomorphic to the toric variety $Z^{\Gamma_m^\nu}$.

\vspace{2mm}

To prove the claim (ii), we distinguish three cases, depending on $m$.

\vspace{2mm}

(a) For $m<l_{g_1+1}(\nu)+e_{g_1+1}$. In this case we have that $g_1>0$ and therefore $m(\nu)$ is either $g_1$ or $1$.
Suppose first that $m(\nu)=1$. Then
\[F_{1,\nu}^{(\frac{l_1(\nu)}{e_1})}= {z^{(\langle\nu,\gamma_1\rangle)}}^{n_1}-{x_1^{(\nu_1)}}^{a_1}+G_{1,\nu},\]
where $G_{1,\nu}$ is the polynomial $\sum c_{i_1i_2k}{x_1^{(\nu_1)}}^{i_1}{x_2^{(\nu_2)}}^{i_2}{z^{(\langle\nu,\gamma_1\rangle)}}^k$ with $\langle\nu,(i_1,i_2)+k\gamma_1\rangle=n_1\langle\nu,\gamma_1\rangle=a_1\nu_1$. Moreover, for $0<r<k_1$, $F_{1,\nu}^{(\frac{l_1(\nu)}{e_1}+r)}$ is a quasi-homogeneous polynomial of degree $\frac{l_1(\nu)}{e_1}+r=a_1\nu_1+r$, and using ($\ast$),
\[F_{1,\nu}^{(\frac{l_1(\nu)}{e_1}+r)}={x_1^{(\nu_1)}}^{a_1-1}x_1^{(\nu_1+r)}\cdot U_1^{(r)}+H_{1,\nu}^{(\frac{l_1(\nu)}{e_1}+r)}\]
where $H_{1,\nu}^{(\frac{l_1(\nu)}{e_1}+r)}$ is a polynomial in which $x_1^{(\nu_1+r)}$ do not appear, and $U_1^{(r)}\neq 0$. Analogously, for $1<i<j(m,\nu)$ and $0<r<k_i(\nu)$ we have
\[F_{i,\nu}^{(\frac{l_i(\nu)}{e_i}+r)}= n_i{F_{i-1,\nu}^{(\frac{l_i(\nu)}{e_{i-1}})}}^{n_i-1}F_{i-1,\nu}^{(\frac{l_i(\nu)}{e_{i-1}}+r)}+H_{i,\nu}^{(\frac{l_i(\nu)}{e_i}+r)},\]
where $H_{i,\nu}^{(\frac{l_i(\nu)}{e_i}+r)}$ is a polynomial in which the variable $x_1^{(\nu_1+k_1+\cdots+k_{i-1}+r)}$ does not appear. Moreover
\[F_{i-1,\nu}^{(\frac{l_i(\nu)}{e_{i-1}}+r)}= a_1{x_1^{(\nu_1)}}^{a_1-1}x_1^{(\nu_1+k_1+\cdots+k_{i-1}+r)}+H_{i-1,\nu}^{(r)},\]
where the variable $x_1^{(\nu_1+k_1+\cdots+k_{i-1}+r)}$ does not appear in the polynomial $H_{i-1,\nu}^{(r)}$. Then,
 by Lemma \ref{TechLem} (ii) it follows that for $1=m(\nu)\leq i<j(m,\nu)$ and $0<r<k_i(\nu)$ the equation
$F_i^{(\frac{l_i(\nu)}{e_i}+r)}$ is linear on $x_1^{(\nu_1+k_1+\cdots+k_{i-1}+r)}$ over $D(x_1^{(\nu_1)})$.
We still have to deal with the equations $F_j^{(\frac{l_j(\nu)}{e_j}+r)}$ with $1\leq r\leq [\frac{m-l_j(\nu)}{e_j}]$. Notice that $j(m,\nu)\leq g_1+1$. If $j(m,
\nu)<g_1+1$ then the argument is exactly as before, and if $j(m,\nu)=g_1+1$ then $\left[\frac{m-l_j(\nu)}{e_j}\right]=0$.

\vspace{2mm}

Suppose now that $m(\nu)=g_1>1$. Then $\nu=(0,\nu_2)$ and the generators of $J_m^\nu$ are
\[\begin{array}{ll}
F_{g_1}^{(r)} \mbox{ for }0\le r\leq[\frac{m-l_{g_1}(\nu)}{e_{g_1}}] & \mbox{ if }j(m,\nu)=g_1\\
\\
F_{g_1}^{(r)}\mbox{ for }0\leq r<k_{g_1}(\nu)\mbox{ and }F_{g_1+1}^{(\frac{l_{g_1+1}(\nu)}{e_{g_1+1}})} & \mbox{ if }j(m,\nu)=g_1+1\\
\end{array}\]
For $r>0$ we have that $F_{g_1}^{(r)}\equiv \frac{\partial F_{g_1}}{\partial x_1}(x_1^{(0)},z^{(0)})x_1^{(r)}+\frac{\partial F_{g_1}}{\partial z}(x_1^{(0)},
z^{(0)})z^{(0)}+H_r\mbox{ mod }I^\nu$, with $H_r$ a polynomial where the variables $x_1^{(r)}$ and $z^{(r)}$ do not appear. Since we are looking at $\{x_1^{(0)}
\neq 0\}$, we are outside the singular locus, and we deduce that these equations are linear either in $x_1^{(r)}$ or in $z^{(r)}$. The rest of the proof follows as in the previous part of this case.

\vspace{2mm}

(b) For $l_{g_1+1}(\nu)+e_{g_1+1}\leq m<l_{g_2+1}(\nu)+e_{g_2+1}$. This is only possible when $g_2=g_1+1$.
Notice that if $j(m,\nu)=g_2+1$ then $\left[\frac{m-l_j(\nu)}{e_j}\right]=0$. We just have to study the generators $F_{i,\nu}^{(\frac{l_i(\nu)}{e_i}+r)}$ for $i>g_1$
(which are all the generators of $J_m^\nu$ when $g_1=0$) since the others were studied in the previous case. That is, we have the equations $F_{g_2,\nu}^{(\frac{l_{g_2}
(\nu)}{e_{g_2}}+r)}$ for $0\leq r\leq k_{g_2}(\nu)-1$ if $j(m,\nu)=g_2+1$ and $0\leq r\leq\left[\frac{m-l_{g_2}(\nu)}{e_{g_2}}\right]$ otherwise. If $j(m,\nu)=g_2+1$
we also have the generator $F_{g_2+1,\nu}^{(\frac{l_{g_2+1}(\nu)}{e_{g_2+1}})}$.
\[\begin{array}{lcl}
F_{g_2,\nu}^{(\frac{l_{g_2}(\nu)}{e_{g_2}})} & = & {F_{g_1}^{(\frac{l_{g_2}(\nu)}{e_{g_1}})}}^{n_{g_2}}-{x_1^{(\nu_1)}}^{\alpha_1^{(g_2)}}x_2^{(\nu_2)}{z^{(\langle
\nu,\gamma_1\rangle)}}^{r_1^{(g_2)}}\cdots {F_{g_1-1}^{(\frac{l_{g_1}(\nu)}{e_{g_1-1}})}}^{r_{g_1}^{(g_2)}}+G_{g_2,\nu}\\
\\
F_{g_2,\nu}^{(\frac{l_{g_2}(\nu)}{e_{g_2}}+r)} & = & {x_1^{(\nu_1)}}^{\alpha_1^{(g_2)}}x_2^{(\nu_2+r)}{z^{(\langle\nu,\gamma_1\rangle)}}^{r_1^{(g_2)}}\cdots
{F_{g_1-1}^{(\frac{l_{g_1}(\nu)}{e_{g_1-1}})}}^{r_{g_1}^{(g_2)}}\cdot U_{g_2}^{(r)}+H_r\\
\end{array}\]
where $U_{g_2}^{(r)}$ is $1$ if $\nu\notin\r_1\cup\r_2$ and $\neq 0$ otherwise, and $H_r$ is a polynomial in which the variable $x_2^{(\nu_2+r)}$ does not appear. Then by Lemma \ref{TechLem} (ii) we deduce that every $F_{g_2,\nu}^
{(\frac{l_{g_2}(\nu)}{e_{g_2}}+r)}$ is linear on $x_2^{(\nu_2+r)}$ over $D(x_1^{(\nu_1)})$. And the rest of the argument goes as in the previous case.

(c) For $l_{g_2+1}(\nu)+e_{g_2+1}\leq m$, with the same arguments it is easy to see that for $i>g_2$ and $1\leq r<k_i(\nu)$ each  $F_{i,\nu}^{(\frac{l_i(\nu)}{e_i}+
r)}$ is linear on $x_1^{(\nu_1+k_1+\cdots+k_i+r)}$ over $D(x_1^{(\nu_1)})\cap D(x_2^{(\nu_2)})$.
\hfill$\Box$

\vspace{3mm}

In particular we have the following variation of the codimension of $C_m^\nu$ as $m$ grows.
\begin{Cor}
For $\nu\in H_m\cup L_m$ such that $\nu\in H_{m-1}\cup L_{m-1}$ we have that
\[\mbox{Codim}(C_m^\nu)=\left\{\begin{array}{ll}
\mbox{Codim}(C_{m-1}^\nu)+1 & \mbox{ if }m\equiv 0\mbox{ mod }e_{j(m-1,\nu)}\\
\\
\mbox{Codim}(C_{m-1}^\nu) & \mbox{ otherwise}\\
\end{array}\right.\]
\label{CorCodim}
\end{Cor}

\subsection{Inclusions among the $C_m^\nu$}
We have a collection of irreducible sets $\{C_m^\nu\ |\ \nu\in H_m\cup L_m\}$ covering $\pi_m^{-1}(X_{Sing})$, but in general it is not its decomposition in irreducible components. We have to study the inclusions
\begin{equation}
 C_m^{\nu'}\subseteq C_m^\nu\mbox{ for different }\nu,\nu'\in H_m\cup L_m.
\label{ContGen}
\end{equation}
We need to  see $C_m^\nu$ as the closure of a set, which is slightly different from $D_m^\nu$, though described by the ideals $I^\nu$ and $J_m^\nu$. For instance, when $j(m,\nu)=c(\nu)=g_1>0$, by Corollary \ref{Corolario3} $J_m^\nu=(F_{g_1,\nu}^{(0)},\ldots,F_{g_1,\nu}^{([\frac{m}{e_{g_1}}])})$, and
\[C_m^\nu=\overline{V(I^\nu,F_{g_1,\nu}^{(0)},\ldots,F_{g_1,\nu}^{([\frac{m}{e_{g_1}}])})\cap D(x_1^{(\nu_1)})},\]
because the polynomials $F_{g_1,\nu}^{(0)},\ldots,F_{g_1,\nu}^{([\frac{m}{e_{g_1}}])}$ do not depend on $x_2^{(\nu_2)}$, and hence, when taking the Zariski closure, we can drop the condition $D(x_2^{(\nu_2)})$ in the description of $D_m^\nu$ given in (\ref{eqast}). This is the description we are looking for, and it is the content of the next Lemma.

\begin{Lem}
For $m\in\Z_{>0}$ and $\nu\in H_m\cup L_m$ we have that $C_m^\nu=\overline{O_m^\nu}$, where
\[O_m^\nu:=\left\{\begin{array}{cl}
V(I^\nu,J_m^\nu) & \mbox{ if }j'(m,\nu)<m(\nu)\\
\\
V(I^\nu,J_m^\nu)\cap D(x_1^{(\nu_1)}) & \mbox{ if }m(\nu)\leq j'(m,\nu)\leq g_2\\
\\
V(I^\nu,J_m^\nu)\cap D(x_1^{(\nu_1)})\cap D(x_2^{(\nu_2)}) & \mbox{ if }j'(m,\nu)\geq g_2+1\\
\end{array}\right.\]
Notice that when $j'(m,\nu)<m(\nu)$ then $O_m^\nu=C_m^\nu$, and if $j'(m,\nu)\geq g_2+1$ then $O_m^\nu=D_m^\nu$.
\label{LemCm}
\end{Lem}

{\em Proof.} If $j'(m,\nu)<m(\nu)$ we have two possibilities regarding $j(m,\nu)$, either $j(m,\nu)<m(\nu)$ or $j(m,\nu)=m(\nu)$. Suppose that $j(m,\nu)<m(\nu)$, then $j(m,\nu)=0$, since $l_{j(m,\nu)}(\nu)<l_{j(m,\nu)+1}(\nu)$ and $l_1(\nu)=\cdots=l_{m(\nu)}(\nu)<l_{m(\nu)+1}(\nu)$. Then $\nu\in H_m$, and by Proposition \ref{CHm} we have that $C_m^\nu=V(I^\nu,J_m^\nu)$ where $J_m^\nu=(z^{(0)},\ldots,z^{([\frac{m}{n}])})$. If $j(m,\nu)=m(\nu)$ we have that $J_m^\nu$
is the ideal
\begin{equation}
J_m^\nu=\left\{\begin{array}{cl}
                 \left(F_{g_1,\nu}^{(0)}\right) & \mbox{ if }c(\nu)=g_1>0\\
                 \\
                 \left(z^{(0)},\ldots,z^{(\langle\nu,\gamma_1\rangle-1)},F_{m(\nu),\nu}^{(\frac{l_{m(\nu)}(\nu)}{e_{m(\nu)}})}\right) & \mbox{ otherwise}\\
                 \end{array}\right.
\label{eqJ}
\end{equation}
Therefore the conditions $D(x_1^{(\nu_1)})\cap D(x_2^{(\nu_2)})$ disappear when taking the Zariski closure.

Suppose now that $m(\nu)\leq j'(m,\nu)\leq g_2$. We prove that
\begin{equation}
\overline{V(I^\nu,J_m^\nu)\cap D(x_1^{(\nu_1)})\cap D(x_2^{(\nu_2)})}=\overline{V(I^\nu,J_m^\nu)\cap D(x_1^{(\nu_1)})},
\label{igD}
\end{equation}
or in other words, the open condition $x_2^{(\nu_2)}\neq 0$ is superfluous when taking the Zariski closure. This claim is obvious for $m(\nu)\leq j'(m,\nu)\leq g_1$, since $x_2^{(\nu_2)}$ appears in the generators of $J_m^\nu$ at most once (in $F_{g_1+1,\nu}^{(\frac{l_{g_1+1}(\nu)}{e_{g_1+1}})}$ if
$l_{g_1+1}(\nu)\leq m<l_{g_1+1}(\nu)+e_{g_1+1}$). For $g_1<j'(m,\nu)\leq g_2$, we are necessarily in the case $g_2=g_1+1$ and we prove the
equality (\ref{igD}) by induction on $m$. For $m=l_{g_2}(\nu)+e_{g_2}$, if the equality (\ref{igD}) does not hold, then
\[C:=\overline{V(I^\nu,J_m^\nu)\cap D(x_1^{(\nu_1)})\cap\{x_2^{(\nu_2)}=0\}}\nsubseteq C_m^\nu.\]
Notice that by Corollary \ref{Corolario3}
\[F_{g_2,\nu}^{(\frac{l_{g_2}(\nu)}{e_{g_2}})}={F_{g_2-1,\nu}^{(\frac{l_{g_2}(\nu)}{e_{g_2-1}})}}^{n_{g_2}}-c_{g_2}{x_1^{(\nu_1)}}^{\alpha_1^{(g_2)}}x_2^{(\nu_2)}
\cdots {F_{g_2-2,\nu}^{(\frac{l_{g_2-1}(\nu)}{e_{g_2-2}})}}^{r_{g_2-1}^{(g_2)}}+G_{g_2,\nu}\]
and by Lemma \ref{TechLem} (i) we deduce that, if $x_2^{(\nu_2)}=0$ then $F_{g_2-1,\nu}^{(\frac{l_{g_2}(\nu)}{e_{g_2-1}})}$. Moreover,
 since the polynomials $F_\nu^{(l)}$ are quasi-homogeneous
\[F_{g_2,\nu}^{(\frac{l_{g_2}(\nu)}{e_{g_2}}+1)}\equiv\bar c_{g_2}{x_1^{(\nu_1)}}^{\alpha_1^{(g_2)}}x_2^{(\nu_2+1)}\cdots {F_{g_2-2,\nu}^{(\frac{l_{g_2-1}(\nu)}{e_{g_2-2}})}}^{r_{g_2-1}^{(g_2)}}\mbox{ mod }x_2^{(\nu_2)}\]
and, by (\ref{eqO}), we deduce that $x_2^{(\nu_2+1)}=0$. Hence      \[C=V(I^\nu,J_{l_{g_2}(\nu)-1}^\nu,x_2^{(\nu_2)},x_2^{(\nu_2+1)},F_{g_2-1,\nu}^{(\frac{l_{g_2}(\nu)}{e_{g_2-1}})}),\]
 and now consider the closed
set $C':=\pi_{m,m-1}(C_m^\nu)=\overline{V(I^\nu,J_{m-1}^\nu)\cap D(x_1^{(\nu_1)})}$. We have that $\pi_{m,m-1}^{-1}(C')=\overline{V(I^\nu,
J_m^\nu)\cap D(x_1^{(\nu_1)})}=C_m^\nu\cup C$ with Codim$(C_m^\nu)=\mbox{Codim}(C')+1$ and Codim$(C)=\mbox{Codim}(C')+2$, which is a contradiction. Suppose it true for $m$ and we prove it for $m+1$. Consider $C':=\pi_{m+1,m}(C_m^\nu)=\overline{V(I^\nu,J_m^\nu)\cap D(x_1^{(\nu_1)})\cap D(x_2^{(\nu_2)})}$. By
induction hypothesis $C'=\overline{V(I^\nu,J_m^\nu)\cap D(x_1^{(\nu_1)})}$, and then $\pi_{m+1,m}^{-1}(C')=\overline{V(I^\nu,J_m^\nu,
F_\nu^{(m+1)})\cap D(x_1^{(\nu_1)})}$. If $F_\nu^{(m+1)}= 0$ then we are done. Otherwise, by Corollary \ref{Corolario3}, $F_\nu^{(m+1)}={F_{g_2,\nu}^
{(\frac{l_{g_2}(\nu)}{e_{g_2}}+r)}}^{e_{g_2}}$ where $r=\frac{m+1-l_{g_2}(\nu)}{e_{g_2}}$, and, as in the first step of
induction, if it subdivides as
\[C_{m+1}^\nu\cup \overline{V(I^\nu,J_m^\nu,F_\nu^{(m+1)})\cap D(x_1^{(\nu_1)})\cap\{x_2^{(\nu_2)}=0\}},\]
 then Codim$(\overline{V(I^\nu,J_m^\nu,F_\nu^{(m+1)})\cap D(x_1^{(\nu_1)})\cap\{x_2^{(\nu_2)}=0\}})=\mbox{Codim}(C')+2$ which is a contradiction.

Finally, if $j'(m,\nu)\geq g_2+1$ there is nothing to prove.
\hfill $\Box$

\vspace{3mm}

We will describe a set $F_m\subset H_m\cup L_m$ such that $\{C_m^\nu\ |\ \nu\in F_m\}$ is the set of irreducible components. The process of
defining $F_m$  as a subset of $H_m\cup L_m$ is done in two steps. The first reduction is easy. We consider the product ordering $\leq_p$ in $\Z^2$ given by:
\begin{equation}
 \nu\leq_p\nu'\mbox{ if and only if }\nu_i\leq\nu'_i\mbox{ for }i=1,2.
\label{ordenGen}
\end{equation}

\begin{Pro}
 For $\nu,\nu'\in H_m\cup L_m$ we have that
\begin{enumerate}
 \item [(i)] If $C_m^{\nu'}\subseteq C_m^\nu$ then $\nu\leq_p\nu'$.
 \\
\item [(ii)] Moreover if $\nu,\nu'\in H_m\cup L_m^=$ then we have
\[C_m^{\nu'}\subseteq C_m^\nu\Longleftrightarrow \nu\leq_p\nu'.\]

\end{enumerate}
\label{Lema2Gen}
\end{Pro}

{\em Proof.}

\begin{enumerate}
\item [(i)] Suppose that $\nu$ and $\nu'$ are not comparable. Then we can assume that $\nu_1<\nu_1'$ and $\nu_2>\nu_2'$. Then, since $C_m^\nu\subseteq V(I^\nu)$, and $C_m^{\nu'}\subset V(I^{\nu'})$, it follows that
\[C_m^\nu\nsubseteq C_m^{\nu'}\mbox{ and }C_m^{\nu'}\nsubseteq C_m^\nu.\]
\item[(ii)] The claim follows by (i) and the definition of $C_m^\nu$ for $\nu\in
 H_m\cup L_m^=$.
\end{enumerate}
\hfill$\Box$

\begin{Defi}
According to Proposition \ref{Lema2Gen} we define the set:
\[P_m=min_{\leq_p}\{H_m\cup L_m^=\}.\]
\label{defPm}
\end{Defi}
The second reduction, which defines the set $F_m\subseteq P_m\cup L_m^<$, is much more involved, and the singular locus of the approximated roots play a role now when studying the inclusions $C_m^{\nu'}\subseteq C_m^\nu$ for different elements $\nu$ and $\nu'$ in $P_m\cup L_m^<$. By Proposition \ref{Lema2Gen} (i)
we have to consider $\nu\lneq_p\nu'$, where, by definition of $P_m$, $\nu\in L_m^<$ and $\nu'\in P_m\cup L_m^<$.

\begin{Pro}
Given $m\in\Z_{>0}$, $\nu\in L_m^<$ and $\nu'\in P_m\cup L_m^<$ with $\nu\lneq_p\nu'$, such that $\nu'-\nu\in\s_{Reg,j'(m,\nu)}$ then $C_m^{\nu'}\subseteq C_m^\nu$.
\label{PropC1}
\end{Pro}

{\em Proof.} We simplify notation by setting $k_i(\nu)=\frac{l_{i+1}(\nu)-l_i(\nu)}{e_i}$, for $1\leq i\leq g$. By the
description of $\s_{Reg,j}$ given in (\ref{SigRegj}) we have to prove the inclusion $C_m^{\nu'}\subseteq C_m^\nu$ when $\nu'-\nu\in\s_{Reg,j'(m,\nu)}$ and
$1\leq j'(m,\nu)\leq g_2$. Then, by Lemma \ref{LemCm}, we have that
\[C_m^\nu=\overline{V(I^\nu,J_m^\nu)\cap D(x_1^{(\nu_1)})}.\]
Suppose first that $\s_{Reg,j'(m,\nu)}=\r_2$, then $\nu'=\nu+(0,\beta)$ with $\beta>0$. We distinguish two cases:

\vspace{2mm}

(i) If $\nu'\in H_m$, then by Proposition \ref{CHm},
\[C_m^{\nu'}=V(I^\nu,x_2^{(\nu_2)},\ldots,x_2^{(\nu_2+\beta-1)},z^{(0)},\ldots,z^{([m/n])}).\]
Note that $g_1=0$, because otherwise $l_1(\nu')=l_1(\nu)$ and $\nu'\notin H_m$. Then, since $1\leq j'(m,\nu)\leq g_2$, we deduce that $g_2=1$.
There exists $1\leq r<k_1(\nu)$ such that
\[l_1(\nu)+re_1\leq m<l_1(\nu)+(r+1)e_1\]
since $\nu\in L_m^<$. Then
\[J_m^\nu=\left(z^{(0)},\ldots,z^{(\langle\nu,\gamma_1\rangle-1)},F_{1,\nu}^{(\frac{l_1(\nu)}{e_1})},\ldots,F_{1,\nu}^{(\frac{l_1(\nu)}{e_1}+r)}\right).\]
Notice that $[\frac{m}{n}]=[\frac{l_1(\nu)+re_1}{n}]=\langle\nu,\gamma_1\rangle+\alpha$, where $\alpha=[\frac{r}{n_1}]$. Now, since $\nu'=\nu+(0,\beta)
\in H_m$ and $g_2=1$, we have that $l_1(\nu')=l_1(\nu)+e_1\beta\geq m+1$ and it follows that $\beta>r$.
Hence we have to prove
\begin{equation}
F_{1,\nu}^{(\frac{l_1(\nu)}{e_1}+l)}\equiv 0\mbox{ mod }(x_2^{(\nu_2)},\ldots,x_2^{(\nu_2+\beta-1)},F_{1,\nu}^{(\frac{l_1(\nu)}{e_1})},\ldots,F_{1,\nu}^{(\frac{l_1(\nu)}{e_1}+l-1)})
\label{eqqq}
\end{equation}
for $0\leq l\leq r$. By Corollary \ref{Corolario3} we have $F_{1,\nu}^{(\frac{l_1(\nu)}{e_1})}={z^{(\langle\nu,\gamma_1\rangle)}}^{n_1}-{x_1^{(\nu_1)}}^{a_1}x_2^{(\nu_2)}+G_{1,\nu}$.
And by Lemma \ref{TechLem} (i), if $x_2^{(\nu_2)}=0$ then $z^{(\langle\nu,\gamma_1\rangle)}=0$, and hence $G_{1,\nu}=0$.

By quasi-homogeneity we can write
\[F_{1,\nu}^{(\frac{l_1(\nu)}{e_1}+1)}=c_1{z^{(\langle\nu,\gamma_1\rangle)}}^{n_1-1}z^{(\langle\nu,\gamma_1\rangle+1)}+c_2{x_1^{(\nu_1)}}^{a_1-1}x_1^{(\nu_1+1)}
x_2^{(\nu_2)}+c_3{x_1^{(\nu_1)}}^{a_1}x_2^{(\nu_2+1)}+G_{1,\nu}^{(1)}\]
where $c_1,c_2,c_3$ are certain coefficients and $G_{1,\nu}^{(1)}$ is a quasi-homogeneous polynomial of degree $\frac{l_1(\nu)}{e_1}+1$. We have that  $G_{1,\nu}^{(1)}=0$ when $\nu\notin\r_1\cup\r_2$, and otherwise, we can apply the same arguments as in the proof of Lemma \ref{TechLem} we deduce that
\[F_{1,\nu}^{(\frac{l_1(\nu)}{e_1}+1)}\equiv 0\mbox{ mod }(x_2^{(\nu_2)},x_2^{(\nu_2+1)},z^{(\langle\nu,\gamma_1\rangle)}).\]
Again by quasi-homogeneity, if $n_1<r$,
\[\begin{array}{ll}
F_{1,\nu}^{(\frac{l_1(\nu)}{e_1}+n_1)} & =z^{(\langle\nu,\gamma_1\rangle)}h_{1,0}(z^{(\langle\nu,\gamma_1\rangle)},\ldots,z^{(\langle\nu,\gamma_1\rangle+n_1)})+
{z^{(\langle\nu,\gamma_1\rangle+1)}}^{n_1}+x_2^{(\nu_2)}h_{2,0}(x_1^{(\nu_1)},\ldots,x_1^{(\nu_1+n_1)})+\\
 \\
 & \ x_2^{(\nu_2+1)}h_{2,1}(x_1^{(\nu_1)},\ldots,x_1^{(\nu_1+n_1)})+\cdots+x_2^{(\nu_2+n_1)}h_{2,n_1}(x_1^{(\nu_1)},\ldots,x_1^{(\nu_1+n_1)})+G_{1,\nu}^{(n_1)}.\\
\end{array}\]
And analogously we prove that
\[F_{1,\nu}^{(\frac{l_1(\nu)}{e_1}+n_1)}\equiv 0\mbox{ mod }(x_2^{(\nu_2)},\ldots,x_2^{(\nu_2+n_1)},z^{(\langle\nu,\gamma_1\rangle)},z^{(\langle\nu,\gamma_1\rangle+1)}).\]
And in general, for $1\leq k\leq r$
\[\begin{array}{ll}
F_{1,\nu}^{(\frac{l_1(\nu)}{e_1}+k)} & = z^{(\langle\nu,\gamma_1\rangle)}h_{1,0}(z^{(\langle\nu,\gamma_1\rangle)},\ldots,z^{(\langle\nu,\gamma_1\rangle+k)})+
z^{(\langle\nu,\gamma_1\rangle+1)}
h_{1,1}(z^{(\langle\nu,\gamma_1\rangle+1)},\ldots,z^{(\langle\nu,\gamma_1\rangle+k)})\\
\\
 & \ +\ \cdots\ +z^{(\langle\nu,\gamma_1\rangle+[\frac{k}{n_1}])}h_{1,[k/n_1]}(z^{(\langle\nu,\gamma_1\rangle+[\frac{k}{n_1}])},\ldots,z^{(\langle\nu,\gamma_1
 \rangle+k)})+\\
 \\
 & +h_{2,0}(x_1^{(\nu_1)},\ldots,x_1^{(\nu_1+k)})x_2^{(\nu_2)}+h_{2,1}(x_1^{(\nu_1)},\ldots,x_1^{(\nu_1+k)})x_2^{(\nu_2+1)}+\cdots+\\
 \\
  & +h_{2,k}(x_1^{(\nu_1)},\ldots,x_1^{(\nu_1+k)})x_2^{(\nu_2+k)}\\
\end{array}\]
where $h_{1,0},h_{1,1},\ldots,h_{1,[k/n_1]},h_{2,0},\ldots,h_{2,k}$ are polynomials. And since $k\leq r<\beta$ and $[\frac{m}{n}]=\langle\nu,\gamma_1\rangle+\alpha$
with $\alpha=[\frac{r}{n_1}]\geq[\frac{k}{n_1}]$, we have proved (\ref{eqqq}) as we wanted.

\vspace{2mm}

(ii) If $\nu'\in L_m$, then by Lemma \ref{LemCm}
\[C_m^{\nu'}=\overline{V(I^{\nu'},J_m^{\nu'})\cap D(x_1^{(\nu_1')})}.\]
Since $\nu\lneq_p\nu'$ we have that $I^\nu\subseteq I^{\nu'}$. We are going to prove that the generators of
$J_m^\nu$ modulo $(I^{\nu'},z^{(l)})_{0\leq l<\frac{l_1(\nu')}{n}}$ belong to $J_m^{\nu'}$. Since $l_i(\nu)=l_i(\nu')$ for $1\leq i\leq g_1$, we have that
\[\begin{array}{ll}
F_{i,\nu}^{(\frac{l_i(\nu)}{e_i}+r_i)}=F_{i,\nu'}^{(\frac{l_i(\nu')}{e_i}+r_i)} & \mbox{ for }1\leq i\leq g_1-1
\mbox{ and }0\leq r_i<k_i(\nu)=k_i(\nu'),\\
\\
F_{g_1,\nu}^{(\frac{l_{g_1}(\nu)}{e_{g_1}}+r)}=F_{g_1,\nu'}^{(\frac{l_{g_1}(\nu')}{e_{g_1}}+r)} & \mbox{ for }
0\leq r<k_{g_1}(\nu)<k_{g_1}(\nu').\\
\end{array}\]
When $J_m^\nu$ has more generators, that is, when $j(m,\nu)>g_1$, then $j'(m,\nu)\geq g_1$. We distinguish two cases.
\begin{enumerate}
\item[(ii.a)] If $j'(m,\nu)=g_1$, then, by Corollary \ref{Corolario3},
\[J_m^\nu=\left(z^{(0)},\ldots,z^{(\frac{l_1(\nu)}{n}-1)};F_{i,\nu}^{(\frac{l_i(\nu)}{e_i}+r_i)},1\leq i\leq g_1,0\leq r_i<k_i(\nu);F_{g_1+1,\nu}^{(\frac{l_{g_1+1}(\nu)}{e_{g_1+1}})}\right)_{l_i(\nu)<l_{i+1}(\nu)},\]
and $l_{g_1+1}(\nu)\leq m<l_{g_1+1}(\nu)+e_{g_1+1}$. By (\ref{ecsFs}) we have that
\[F_{g_1+1,\nu}^{(\frac{l_{g_1+1}(\nu)}{e_{g_1+1}})}= {F_{g_1,\nu}^{(\frac{l_{g_1+1}(\nu)}{e_{g_1}})}}^{n_{g_1+1}},\]
since $\nu_2'>\nu_2$. Notice that $l_{g_1+1}(\nu)=l_{g_1}(\nu)+l_{g_1+1}(\nu)-l_{g_1}(\nu)=l_{g_1}(\nu')+\alpha\leq m$, with $\alpha>0$. Then $l_{g_1}(\nu')<m$ and
therefore $j(m,\nu')\geq g_1$.
If $j(m,\nu')>g_1$ then
\[F_{g_1,\nu}^{(\frac{l_{g_1+1}(\nu)}{e_{g_1}})}\in J_m^{\nu'},\]
 since $\frac{l_{g_1+1}(\nu)}{e_{g_1}}=
\frac{l_{g_1}(\nu')}{e_{g_1}}+k_{g_1}(\nu)$ and $k_{g_1}(\nu)<k_{g_1}(\nu')$. If $j(m,\nu')=g_1$ then
\[F_{g_1,\nu}^{(\frac{l_{g_1+1}(\nu)}{e_{g_1}})}\in J_m^{\nu'},\]
 since $\frac{l_{g_1+1}(\nu)}{e_{g_1}}=\frac{l_{g_1}(\nu')}{e_{g_1}}+k_{g_1}(\nu)$, and $k_{g_1}(\nu)\leq[\frac{m-l_{g_1}
(\nu')}{e_{g_1}}]$, because $l_{g_1}(\nu)=l_{g_1}(\nu')$ and $m\geq l_{g_1+1}(\nu)$.

\vspace{2mm}

\item[(ii.b)] If $j'(m,\nu)>g_1$, then we are in the case $g_2=g_1+1$ and $j'(m,\nu)=g_1+1$. There exists an integer $1\leq r<k_{g_1+1}(\nu)$ such that
\begin{equation}
l_{g_1+1}(\nu)+re_{g_1+1}\leq m<l_{g_1+1}(\nu)+(r+1)e_{g_1+1}.
\label{toro}
\end{equation}
Then $J_m^\nu=\left(z^{(0)},\ldots,z^{(\frac{l_1(\nu)}{n}-1)},F_{1,\nu}^{(\frac{l_1(\nu)}{e_1})},\ldots,F_{g_1,\nu}^{(\frac{l_{g_1}}{e_{g_1}}+k_{g_1}(\nu)-1)},
F_{g_1+1,\nu}^{(\frac{l_{g_1+1}(\nu)}{e_{g_1+1}})},\ldots,F_{g_1+1,\nu}^{(\frac{l_{g_1+1}(\nu)}{e_{g_1+1}}+r)}\right)$, where
\[F_{g_1+1,\nu}^{(\frac{l_{g_1+1}(\nu)}{e_{g_1+1}})}= {F_{g_1,\nu}^{(\frac{l_{g_1+1}(\nu)}{e_{g_1}})}}^{n_{g_1+1}}-{x_1^{(\nu_1)}}^{\alpha_1^{(g_1+1)}}x_2^{(\nu_2)}
{z^{(\langle\nu,\gamma_1\rangle)}}^{r_1^{(g_1+1)}}\cdots {F_{g_1-1,\nu}^{(\frac{l_{g_1-1}(\nu)}{e_{g_1-1}})}}^{r_{g_1}^{(g_1+1)}}.\]
And, analogously to case (i), we can write the polynomials $F_{g_1+1,\nu}^{(\frac{l_{g_1+1}(\nu)}{e_{g_1+1}}+k)}$ for $1\leq k\leq r$ as
\[\begin{array}{ll}
F_{g_1+1,\nu}^{(\frac{l_{g_1+1}(\nu)}{e_{g_1+1}}+k)} & =F_{g_1,\nu}^{(\frac{l_{g_1+1}(\nu)}{e_{g_1}})}h_{1,0}+F_{g_1,\nu}^{(\frac{l_{g_1+1}(\nu)}{e_{g_1}}+1)}h_{1,1}+
\cdots+F_{g_1,\nu}^{(\frac{l_{g_1+1}(\nu)}{e_{g_1}}+
[\frac{k}{n_{g_1+1}}])}h_{1,[k/n_{g_1+1}]}+\\
\\
 & h_{2,0}x_2^{(\nu_2)}+\cdots+h_{2,k}x_2^{(\nu_2+k)},\\
 \end{array}\]
for certain polynomials $h_{1,0},\ldots,h_{1,[k/n_{g_1+1}]},h_{2,0},\ldots,h_{2,k}$.

 If $j(m,\nu')=g_1$ then $l_{g_1+1}(\nu')>m$ and by (\ref{toro}) it follows that $\beta>r$. Moreover
 \[J_m^{\nu'}=\left(z^{(0)},\ldots,z^{(\frac{l_1(\nu')}{n}-1)},F_{1,\nu}^{(\frac{l_1(\nu')}{e_1})},\ldots,F_{g_1,\nu}^{(\frac{l_{g_1}(\nu')}{e_{g_1}})},\ldots,F_{g_1,\nu}^
 {(\frac{l_{g_1}(\nu')}{e_{g_1}}+
 \alpha)}\right)\]
 with $\alpha=[\frac{m-l_{g_1}(\nu')}{e_{g_1}}]$.
We have that $[\frac{m-l_{g_1}(\nu')}{e_{g_1}}]=[\frac{m-l_{g_1}(\nu)}{e_{g_1}}]=[\frac{m}{e_{g_1}}]-\frac{l_{g_1}(\nu)}{e_{g_1}}$, since $\frac{l_{g_1}(\nu)}
{e_{g_1}}$ is an integer, and by (\ref{toro}) we have that $\langle\nu,\gamma_{g_1+1}\rangle+\frac{r}{n_{g_1+1}}\leq\frac{m}{e_{g_1}}<\langle\nu,\gamma_{g_1+1}
\rangle+\frac{r+1}{n_{g_1+1}}$. Then $[\frac{m}{e_{g_1}}]=\langle\nu,\gamma_{g_1+1}\rangle+[\frac{r}{n_{g_1+1}}]$, and we have that
\[\frac{l_{g_1+1}(\nu)}{e_{g_1}}+\left[\frac{k}{n_{g_1+1}}\right]\leq\frac{l_{g_1}(\nu)}{e_{g_1}}+\left[\frac{m-l_{g_1}(\nu)}{e_{g_1}}\right].\]
It follows that $F_{g_1+1,\nu}^{(\frac{l_{g_1+1}(\nu)}{e_{g_1+1}}+k)}$ belongs to $J_m^{\nu'}$.

If $j(m,\nu')=g_1+1$ then $\nu'\in N_{g_+1}$ and therefore $\beta=qn_{g_1+1}$ with $q\in\Z_{>0}$. Moreover there exists an integer $0<r'<k_{g_1+1}(\nu')$ such that
\[l_{g_1+1}(\nu')+r'e_{g_1+1}\leq m<l_{g_1+1}(\nu')+(r'+1)e_{g_1+1}.\]
Or equivalently $l_{g_1+1}(\nu)+(r'+qn_{g_1+1})e_{g_1+1}\leq m<l_{g_1+1}(\nu)+(r'+qn_{g_1+1}+1)e_{g_1+1}$. Then by (\ref{toro}) it follows that
\[r=r'+qn_{g_1+1},\]
and therefore $r>\beta$. Moreover $k_{g_1}(\nu')=k_{g_1}(\nu)+qn_{g_1+1}$. Then
\[J_m^{\nu'}=\left(z^{(0)},\ldots,z^{(\frac{l_1(\nu')}{n}-1)},F_{1,\nu}^{(\frac{l_1(\nu')}{e_1})},\ldots,F_{g_1,\nu}^{(\frac{l_{g_1}(\nu')}{e_{g_1}}+k_{g_1}(\nu')-1)},
F_{g_1+1,\nu}^{(\frac{l_{g_1+1}(\nu')}{e_{g_1+1}})},\ldots,F_{g_1+1,\nu}^{(\frac{l_{g_1+1}(\nu')}{e_{g_1+1}}+r')}\right),\]
where notice that $\frac{l_{g_1}(\nu')}{e_{g_1}}+k_{g_1}(\nu')-1=\frac{l_{g_1}(\nu)}{e_{g_1}}+k_{g_1}(\nu)+qn_{g_1+1}-1$, and
\[J_m^\nu=\left(z^{(0)},\ldots,z^{(\frac{l_1(\nu)}{n}-1)},F_{1,\nu}^{(\frac{l_1(\nu)}{e_1})},\ldots,F_{g_1,\nu}^{(\frac{l_{g_1}(\nu)}{e_{g_1}}+k_{g_1}(\nu)-1)},
F_{g_1+1,\nu}^{(\frac{l_{g_1+1}(\nu)}{e_{g_1+1}})},\ldots,F_{g_1+1,\nu}^{(\frac{l_{g_1+1}(\nu)}{e_{g_1+1}}+r)}\right).\]
Since $I^{\nu'}=(x_1^{(0)},\ldots,x_1^{(\nu_1-1)},x_2^{(0)},\ldots,x_2^{(\nu_2+qn_{g_1+1}-1)})$ it follows
that for $0\leq k\leq r'$ and $s=qn_{g_1+1}+k$
\[F_{g_1+1,\nu}^{(\frac{l_{g_1+1}(\nu)}{e_{g_1+1}}+s)}= F_{g_1+1,\nu}^{(\frac{l_{g_1+1}(\nu')}{e_{g_1+1}}+k)}.\]
Then finally we have to prove that $F_{g_1+1,\nu}^{(\frac{l_{g_1+1}(\nu)}{e_{g_1+1}}+s)}\in J_m^{\nu'}$ for
$0\leq s<qn_{g_1+1}$. This follows as in the previous cases, since
\[\begin{array}{ll}
F_{g_1+1,\nu}^{(\frac{l_{g_1+1}(\nu)}{e_{g_1+1}})} & = {F_{g_1,\nu}^{(\frac{l_{g_1+1}(\nu)}{e_{g_1}})}}^{n_{g_1+1}}-{x_1^{(\nu_1)}}^{\alpha_1^{(g_1+1)}}x_2^{(\nu_2)}
{z^{(\langle\nu,\gamma_1\rangle)}}^{r_1^{(g_1+1)}}\cdots {F_{g_1-1,\nu}^{(\frac{l_{g_1-1}(\nu)}{e_{g_1-1}})}}^{r_{g_1}^{(g_1+1)}}\\
\\
 & \mbox{ and for }0\leq s<qn_{g_1+1}\\
 \\
 F_{g_1+1,\nu}^{(\frac{l_{g_1+1}(\nu)}{e_{g_1+1}}+s)} & = F_{g_1,\nu}^{(\frac{l_{g_1+1}(\nu)}{e_{g_1}})}h_{1,0}+\cdots+F_{g_1,\nu}^{(\frac{l_{g_1+1}(\nu)}{e_{g_1}}+
 [\frac{s}{n_{g_1+1}}])}h_{1,[s/n_{g_1+1}]}\\
 \\
 &\ -h_{2,0}x_2^{(\nu_2)} -\cdots-h_{2,s}x_2^{(\nu_2+s)}\\
 \end{array}\]
 for polynomials $h_{1,0},\ldots,h_{1,[s/n_{g_1+1}]},h_{2,0},\ldots,h_{2,s}$. And it is clear that $F_{g_1+1,\nu}^{(\frac{l_{g_1+1}(\nu)}{e_{g_1+1}}+s)}\in J_m^{\nu'}$.
\end{enumerate}
The key point in all the cases is that $\alpha_2^{(g_1+1)}=1$.

If $\nu'-\nu\in\r_1$ then we are in the case $\gamma_1=(\frac{1}{n_1},\frac{1}{n_1})$. Now $\nu_1'>\nu_1$ and $\nu_2'=\nu_2$, and similar arguments apply to this case to get the inclusion we want to prove.
\hfill$\Box$

\vspace{3mm}

The previous Proposition motivates the following definition.

\begin{Defi}
We consider the relation in $N_0$, depending on $m$ and denoted by $<_{R,m}$, given by
\[\nu\leq_{R,m}\nu'\mbox{ if and only if }\nu\leq_p\nu'\mbox{ and }\nu'-\nu\in\s_{Reg,j'(m,\nu)}.\]
We define the set $F_m=\mbox{min}_{\leq_{R,m}}\{P_m\cup L_m^<\}$.
\label{orden2}
\end{Defi}
Notice that, by (\ref{SigRegj}), for $m$ and $\nu$ such that $j'(m,\nu)>g_2$, this order is just equality.

It is worth pointing out that the inclusions which are described by this last relation in Proposition \ref{PropC1}, can be explained by the fact that even though a curve may be in the singular locus of a quasi-ordinary surface, it may not be part of the singular locus of its first approximated quasi-ordinary surfaces. And as Proposition \ref{Cgeom} explains, the geometry of $C_m^\nu$ is only determined by the geometry of one of its approximated roots, for $m$ small enough. Hence, the jets which project to the singular locus of the surface but not to the singular locus of the approximated surfaces will not give rise to irreducible components of the jet schemes for $m$ small enough, and they will be included in other components.

Now we prove that all possible inclusions are controlled by the relation defined in Definition \ref{orden2} and the product ordering, that is, in the set $F_m$.

\begin{Pro}
Given $m\in\Z_{>0}$ and $\nu,\nu'\in F_m$ with $\nu\leq_p\nu'$ then $C_m^{\nu'}\not\subseteq C_m^\nu$.
\label{PropC2}
\end{Pro}
{\em Proof.} We will prove that $C_m^{\nu'}\not\subseteq C_m^\nu$ by showing that
\begin{equation}
\mbox{Codim}(C_m^{\nu'})\leq\mbox{Codim}(C_m^\nu).
\label{eqCodim}
\end{equation}

First notice that $\nu\in L_m^<$, since otherwise there would not exist $\nu'\neq \nu$ such that $\nu'\in F_m$ and $\nu\leq_p\nu'$. Recall that $\s=\R_{\geq 0}^2$,
we define the set
\[\begin{array}{cl}
E(\nu)_m & =\{\nu'\in(\nu+\s)\cap(P_m\cup L_m^<)\ |\ \nu'\neq\nu\mbox{ and it is minimal with respect to }\leq_{R,m}\mbox{ in }\nu+\s\}\\
\\
 & =\{\nu'\in(\nu+\s)\cap(P_m\cup L_m^<)\ |\ \nu'\neq\nu\mbox{ and }\nexists\ \widetilde\nu\in(\nu+\s)\cap(P_m\cup L_m^<)\mbox{ such that }\widetilde\nu\leq_{R,m}\nu'\}.\\
 \end{array}\]

We claim that for any $\nu'\in E(\nu)_m$ we have that
\begin{equation}
\mbox{Codim}(C_{m_0}^{\nu'})\leq\mbox{Codim}(C_{m_0}^\nu)\ \mbox{ for }l_1(\nu)+e_1\leq m_0<l_{i(\nu)}(\nu).
\label{eqCodim0}
\end{equation}
We prove this claim by induction on $m$. For $m=l_1(\nu)+e_1$ we have that

(i) if $a_1=1$ then we are in the case $\gamma_1=(\frac{1}{n_1},\frac{1}{n_1})$ and $E(\nu)_m=\emptyset$, because $\nu+(1,0),\nu+(0,1)\notin N_1$ and $\nu+
(2,0),\nu+(0,2)\in H_m$ but $\nu\leq_{R,m}\nu+(2,0)$ and $\nu\leq_{R,m}\nu+(0,2)$.

(ii) If $a_1>1$ then $\nu+(1,0)\in H_m$ and it follows that in fact $E(\nu)_m=\{\nu+(1,0)\}$, because the only other possible $\nu'$ is
\[\begin{array}{ll}
\nu'=\nu+(0,1)\in P_m\cup L_m & \mbox{ if }b_1\equiv 0\mbox{ mod }n_1\\

\nu'=\nu+(0,2)\in P_m\cup L_m & \mbox{ otherwise}\\
\end{array}\]
and in both cases we have that $\nu\leq_{R,m}\nu'$. Now, by Lemma \ref{Prop1} we have that for $\nu'=\nu+(1,0)$,
\[\mbox{Codim}(C_m^{\nu'})=\nu_1+\nu_2+\frac{l_1(\nu)}{n}+2=\mbox{Codim}(C_m^\nu).\]

Suppose that the claim is true for $m-1$ and we prove it for $m$. Let $\nu'$ be an element in $E(\nu)_m$.

(i) If $\nu'\in E(\nu)_{m-1}$, by induction hypothesis, we have that $\mbox{Codim}(C_{m-1}^{\nu'})\leq\mbox{Codim}(C_{m-1}^\nu)$. By Corollary \ref{CorCodim}
we know that, passing from $m-1$ to $m$, the codimension of $C_m^\nu$ grows if and only if $m$ is divisible by $e_{j(m-1,\nu)}$, and it grows by one. But since
$\nu\leq_p\nu'$ we have that $j(m-1,\nu')\leq j(m-1,\nu)$ and therefore if $e_{j(m-1,\nu')}$ divides $m$ then $e_{j(m-1,\nu)}$ divides $m$, and it follows that
$\mbox{Codim}(C_m^{\nu'})\leq\mbox{Codim}(C_m^\nu)$.

\vspace{2mm}

(ii) If $\nu'\notin E(\nu)_{m-1}$, there exists $\widetilde\nu\in E(\nu)_{m-1}$ such that $\widetilde\nu\leq_{R,m-1}\nu'$ and $\widetilde\nu\not\leq_{R,m}
\nu'$. By induction hypothesis we have that $\mbox{Codim}(C_{m-1}^{\widetilde\nu})\leq\mbox{Codim}(C_{m-1}^\nu)$, and again, since $\nu\lneq\widetilde\nu$ then
$j(m,\nu)\geq j(m,\widetilde\nu)$ and therefore $\mbox{Codim}(C_m^{\widetilde\nu})\leq\mbox{Codim}(C_m^\nu)$. Now we are going to prove that $\mbox{Codim}
(C_m^{\nu'})\leq\mbox{Codim}(C_m^{\widetilde\nu})$. We have two possibilities, either $\widetilde\nu\in L_m^<$ or $\widetilde\nu\notin L_m^<$.

 If $\widetilde\nu\in L_m^<$, then $m=l_{g_2+1}(\widetilde\nu)+e_{g_2+1}$ and
    \[\pi_{m,m-1}^{-1}(O_{m-1}^{\widetilde\nu})=V(I^{\widetilde\nu},z^{(0)},\ldots,z^{(\langle\widetilde\nu,\gamma_1\rangle-1)},F_{1,\widetilde\nu}^
    {(\frac{l_1(\widetilde\nu)}{e_1})},\ldots,F_{g_2+1,\widetilde\nu}^{(\frac{l_{g_2+1}(\widetilde\nu)}{e_{g_2+1}})},
    F_{g_2+1,\widetilde\nu}^{(\frac{l_{g_2+1}(\widetilde\nu)}{e_{g_2+1}}+1)})\cap D(x_1^{(\widetilde\nu_1)}),\]
where $F_{g_2+1,\widetilde\nu}^{(\frac{l_{g_2+1}(\widetilde\nu)}{e_{g_2+1}})}= {F_{g_2,\widetilde\nu}^{(\frac{l_{g_2+1}(\widetilde\nu)}{e_{g_2}})}}^{n_{g_2+1}}-{x_1^{(\widetilde\nu_1)}}^
{\alpha_1^{(g_2+1)}}{x_2^{(\widetilde\nu_2)}}^{\alpha_2^{(g_2+1)}}\cdots F_{g_2-1,\widetilde\nu}^{(\frac{l_{g_2-1}(\widetilde\nu)}{e_{g_2-1}})}+G_{g_2+1,\widetilde\nu}$, with $\alpha_2^{(g_2+1)}>1$, and $F_{g_2+1,\widetilde\nu}^{(\frac{l_{g_2+1}(\widetilde\nu)}{e_{g_2+1}}+1)}= n_{g_2+1}
{F_{g_2,\widetilde\nu}^{(\frac{l_{g_2+1}(\widetilde\nu)}{e_{g_2}})}}^{n_{g_2+1}-1}F_{g_2,\widetilde\nu}^{(\frac{l_{g_2+1}(\widetilde\nu)}{e_{g_2}}+1)}-x_2^{(\widetilde\nu_2)}H$,
where $H$ is a polynomial in the variables
$H(x_1^{(\widetilde\nu_1)},x_1^{(\widetilde\nu_1+1)},x_2^{(\widetilde\nu_2)},x_2^{(\widetilde\nu_2+1)},\ldots,F_{g_2-1,\widetilde\nu}^{(\frac{l_{g_2-1}(\widetilde\nu)}{e_{g_2-1}})},
F_{g_2-1,\widetilde\nu}^{(\frac{l_{g_2-1}(\widetilde\nu)}{e_{g_2-1}}+1)})$. Then
\[(\pi_{m,m-1}^{-1}(C_{m-1}^{\widetilde\nu}))_{red}=\overline{V(I^{\widetilde\nu},J_m^{\widetilde\nu})\cap D(x_1^{(\widetilde\nu_1)})\cap
V(x_2^{(\widetilde\nu_2)})}\cup\overline{V(I^{\widetilde\nu},J_m^{\widetilde\nu})\cap D(x_1^{(\widetilde\nu_1)})\cap D(x_2^{(\widetilde\nu_2)})},\]
     and it is not difficult to see that $(\pi_{m,m-1}^{-1}(C_{m-1}^{\widetilde\nu}))_{red}=C_m^{\nu'}\cup C_m^{\widetilde\nu}$, where $\nu'=\widetilde\nu+(0,\alpha)$, with
    \[\alpha=\left\{\begin{array}{cl}
    1 & \mbox{ if }g_2=g_1\\
    \mbox{min}\{n_{g_1+1},k_{g_1+1}(\widetilde\nu)\} & \mbox{ otherwise}\\
    \end{array}\right.\]
where remember that $k_i(\widetilde\nu)$ denotes $\frac{l_{i+1}(\widetilde\nu)-l_i(\widetilde\nu)}{e_i}$. In both cases we have, by Proposition \ref{Prop1}, that $\mbox{Codim}(C_m^{\nu'})=\mbox{Codim}(C_{m-1}^{\widetilde\nu})+1=\mbox{Codim}(C_m^{\widetilde\nu})$.

If $\widetilde\nu\notin L_m^<$ then $m=l_{i(\widetilde\nu)}(\widetilde\nu)$ with $i(\widetilde\nu)\leq g_2+1$, since $j'(m-1,\widetilde\nu)\leq g_2$. We have that
$(\pi_{m,m-1}^{-1}(O_{m-1}^{\widetilde\nu}))_{red}=V(I^{\widetilde\nu},J_{m-1}^{\widetilde\nu},F_{i(\widetilde\nu),\widetilde\nu}^{(\frac
{l_{i(\widetilde\nu)}(\widetilde\nu)}{e_{i(\widetilde\nu)}})})\cap D(x_1^{(\widetilde\nu_1)})$, where
\[F_{i(\widetilde\nu),\widetilde\nu}^{(\frac
{l_{i(\widetilde\nu)}(\widetilde\nu)}{e_{i(\widetilde\nu)}})}={x_1^{(\widetilde\nu_1)}}^{\alpha_1^{(i(\widetilde\nu))}}{x_2^{(\widetilde\nu_2)}}^{\alpha_2^
{(i(\widetilde\nu))}} {z^{(\langle\widetilde\nu,\gamma_1\rangle)}}^{r_1^{(i(\widetilde\nu))}}\cdots {F_{i(\widetilde\nu)-2,\widetilde\nu}^{(\frac{l_{i(\widetilde\nu)-2}
(\widetilde\nu)}{e_{i(\widetilde\nu)-2}})}}^{r_{i(\widetilde\nu)-1}^{(i(\widetilde\nu))}}+G_{i(\widetilde\nu),\widetilde\nu}.\]
 Therefore, by Lemma \ref{TechLem}, $F_{i(\widetilde\nu),\widetilde\nu}^{(\frac{l_{i(\widetilde\nu)}(\widetilde\nu)}{e_{i(\widetilde\nu)}})}=0$ implies that $ x_2^{(\widetilde
 \nu_2)}=0$ because $i(\widetilde\nu)-2<g_2$. And, as before, if $g_2=g_1+1$ and $i(\widetilde\nu)=g_2+1$ then we have that $\nu'=\widetilde\nu+(0,\alpha)$ with
 $\alpha=\mbox{min}\{n_{g_1+1},k_{g_1+1}(\widetilde\nu)\}$. Otherwise $\nu'=\widetilde\nu+(0,1)$, and in both cases we have
 \[(\pi_{m,m-1}^{-1}(C_{m-1}^{\widetilde\nu}))_{red}=C_m^{\nu'}\]
 with $\mbox{Codim}(C_m^{\nu'})=\mbox{Codim}(C_{m-1}^{\widetilde\nu})+1$. Since $\widetilde\nu\in E(\nu)_{m-1}$, it follows that $j(m-1,\nu)>j(m-1,\widetilde\nu)=
 i(\widetilde\nu)-1$ and by Corollary \ref{CorCodim} we have that $\mbox{Codim}(C_m^\nu)=\mbox{Codim}(C_{m-1}^\nu)+1$, which finishes the proof.
\hfill $\Box$

\vspace{3mm}

Now we can prove the main theorem of this section.

\begin{The}
 For $m\in\Z_{>0}$ the decomposition of $\pi_m^{-1}(X_{Sing})$ in irreducible components is given by
\[(\pi_m^{-1}(X_{Sing}))_{red}=\bigcup_{\nu\in F_m}C_m^\nu.\]
\label{TheCaso1}
\end{The}

{\em Proof.}  The irreducibility of the sets $C_m^\nu$ was proven in Proposition \ref{Prop1}. And by Proposition \ref{Lema2Gen}, Proposition \ref{PropC1} and Proposition \ref{PropC2} we have that
\[\bigcup_{\nu\in H_m\cup L_m}C_m^\nu=\bigcup_{\nu\in F_m}C_m^\nu.\]
Hence the result follows by Lemma \ref{LemTh}.
\hfill$\Box$

\vspace{3mm}

\begin{Rem}
For $\nu\in N_g$, if $\nu\in F_{l_g(\nu)}$, then $\nu\in F_m$ for every $m\geq l_g(\nu)$, or in other words, $\nu$ gives rise to an irreducible component for any $m\geq l_g(\nu)$.
\label{RemarkFinal}
\end{Rem}

\begin{Rem}
When the equisingular dimension is $c=1$ (see Definition \ref{equi}), then $g_1=g_2=g$. Moreover we have the following properties for $1\leq i\leq g$
\[\begin{array}{l}
l_i(\nu)=l_i(\nu+(0,r)),\ \mbox{ for all }r\in\Z\\
\\
\mbox{if }\nu\in N_i\mbox{ then }\nu+(0,r)\in N_i,\ \mbox{ for all }r\in\Z\\
\end{array}\]
Hence we deduce that for any $m\in\Z_{>0}$ and $\nu\in H_m\cup L_m$ we have $\s_{Reg,j'(m,\nu)}=\r_2$, and therefore $F_m=(P_m\cup L_m^<)\cap\r_1$.

The behaviour of the jet schemes is exactly as the plane curve defined by the Puiseux pairs $\lambda_1^{(1)},\ldots,\lambda_g^{(1)}$.
In \cite{Hcur} the second author describes the irreducible components of jets through the origin in the case of plane curves .
\end{Rem}

The previous remark is the simplest evidence of the fact that the irreducible components are only affected by the topological type. This is proved in Theorem \ref{topTree}

\vspace{3mm}

To any quasi-ordinary surface singularity we can associate a weighted graph, containing information about the irreducible components of jet schemes and how they behave under truncation maps.

\begin{Defi}
The weighted graph of the jet schemes of $X$ is the leveled weighted graph $\Gamma$ defined as follows:

\begin{itemize}
\item for $m\geq 1$ we represent every irreducible components of $\pi_m^{-1}(X_{Sing})$ by a vertex $V_m$, the sub-index $m$ being the level of the vertex;

\item  we join the vertices $V_{m+1}$ and $V_m$ if the canonical morphism $\pi_{m+1,m}$ induces a morphism between the corresponding irreducible components;

\item we weight each vertex by the dimension of the corresponding irreducible component.
\end{itemize}

 We define  $E\Gamma$ to be the weighted graph that we obtain from $\Gamma$
by weighting any vertex of $\Gamma$  by the embedding dimensions of the corresponding irreducible components (note that by the definition of $\Gamma,$ these vertices are also weighted by their dimensions).
\label{defGrafo0}
\end{Defi}

Notice that the data of the codimension together with the embedded dimension permits to distinguish when the vertex corresponds to a hyperplane or a lattice component. Indeed, given a vertex of the graph, let $e$ be the embedded dimension and $c$ the codimension, then the vertex corresponds to a hyperplane component if and only if $e+c=3(m+1)$. Therefore we can extract from $E\Gamma$ a subgraph $\Gamma'$ as follows.

\begin{Defi}
We define a weighted subgraph $\Gamma'$ of $E\Gamma$ by adding the condition that we join  the vertices $V_m$ (corresponding to a certain component, say $C_m^{\nu'}$) and $V_{m-1}$ (corresponding to $C_m^\nu$) only if
\begin{itemize}
\item if $\nu\in L_{m-1}^<$ with $j(m-1,\nu)\leq g_2$ then $\nu'=\nu+(0,\alpha)$ with $\alpha$ minimal among the elements in $F_{m}$.

\item if $\nu\in L_{m-1}^<$ with $j(m-1,\nu)>g_2$ then $\nu'=\nu$.
\end{itemize}
\label{defGrafo}
\end{Defi}

The important thing about this new graph $\Gamma'$ is that, with the weights, we are able to detect when we pass from a hyperplane component at level $m$ to a lattice component at level $m+1$, as we also do in the graph $E\Gamma$, but now we can follow this component in a unique path in the graph as $m$ grows. This will be useful to prove the following result.

\begin{The}
The graph $\Gamma'$ determines and it is determined by the topological type of the singularity.
\label{topTree}
\end{The}
{\em Proof.} Obviously the graph is determined by the semigroup, and therefore, by \cite{Gau}, by the topological type.

To prove the converse we consider two different sets of generators of the semigroup $\{\gamma_1,\ldots,\gamma_g\}$ and $\{\gamma_1',\ldots,\gamma_{g'}'\}$ and we will prove that the corresponding weighted graphs are different too.

Given a weighted graph we can recover the number of characteristic exponents in the following way. Any vertex $V_m$ on the graph comes with the codimension $c(V_m)$ and the embedded dimension $e(V_m)$. Take an infinite branch, and consider the finite part that starts at
\[m_0=\mbox{max }\{m\ |\ V_{m-1}\mbox{ is a hyperplane component and }V_m\mbox{ is a lattice component}\},\]
and ends at
\[m_1=\mbox{min }\{m\ |\  c(V_m)=c(V_{m-1})+1\mbox{ for all }m>m_1\}.\]
We can read $e_0,\ldots,e_{g-1}$ making use of Corollary \ref{CorCodim}. Indeed, along the piece of branch, the vertex $V_m$ corresponds to a component $C_m^\nu$ with $\nu\in N_g$, $m_0=l_1(\nu)$ and $m_1=l_g(\nu)$. To read this data we consider only branches that projects into the component $Z_1$ of the singular locus, since otherwise we can only assure that $m_0=l_{m(\nu)}(\nu)$, and we do not have all the information whenever $m(\nu)>1$. Notice that $Z_1$ is always a component of the singular locus unless we are in the case $g=1$ and $\gamma=(\frac{1}{n},\frac{1}{n})$, which is very easy to recognize. Indeed, it is the only case when at level $m=1$ we have only one component, with codimension 3 and embedded dimension 0. Moreover the multiplicity $n$ equals the first time $m$ when we have a lattice component. Therefore this simple case is very easily understood in the graph. For the rest of the cases, since we know that $\nu\in N_g$, going backwards we look for the biggest $m'$ such that  $c(V_{m'})=c(V_{m_0})-1$. Then $n=m_0-m'$. Now, going from level $m_0$ to $m_1$, we know that the codimension grows by one exactly every $e_1$ steps at first, after every $e_2$ steps, and so on. Since $e_1>e_2>\cdots>e_g=1$ we can read these numbers on the graph. Notice that equivalently we get $n_1,\ldots,n_g$, and in particular we have $g$, the number of characteristic exponents.

Suppose now that the number of generators of the semigroups is the same, say $g$. We will prove by induction on $g$ that the graphs corresponding to different sets of generators, are different. We denote the vertices at level $m$ by $V_m(c(V_m),e(V_m))$. For $g=1$, the multiplicity is read from the graph as was explained before, and the situation for $m=1$ is:
\[\begin{array}{cl}
\bullet\ V_1(3,0) & \mbox{ if }\gamma=(\frac{1}{n},\frac{1}{n})\\
\\
\bullet\ V_1(2,1) & \mbox{ if }\gamma=(\frac{a}{n},\frac{1}{n}),\mbox{ with }a>1\\
\\
\bullet\ V_1(2,1)\ \ \bullet\ V_1(2,1) & \mbox{ if }\gamma=(\frac{a}{n},\frac{b}{n}),\mbox{ with }b>1\\
\end{array}\]
If we want to compare the graph associated to $\gamma$ and the graph associated to $\gamma'$, we just have to consider the cases $\gamma=(\frac{a}{n},\frac{1}{n})$, $\gamma'=(\frac{a'}{n},\frac{1}{n})$ with $a\neq a'$, and $\gamma=(\frac{a}{n},\frac{b}{n})$, $\gamma'=(\frac{a'}{n},\frac{b'}{n})$ with $b,b'\neq 1$ and $\gamma\neq\gamma'$. The first case is very easy to distinguish, since the first moment a component splits in two is at $m=a$ for one graph, and at $m=a'$ for the other. For the other case, first note that the graph of any quasi-ordinary with only one characteristic exponent $\gamma=(\frac{a}{n},\frac{b}{n})$ is the graph associated to $z^n-x_1^ax_2^b$. The key point is that, when $b>1$ we have in the graph a branch which generically corresponds to $x_2^{(0)}\neq 0$ (resp. $x_1^{(0)}\neq 0$), that is, it behaves like the graph of the curve $z^n-x_1^a$ (resp. $z^n-x_2^b$). Therefore comparing graphs associated to $\gamma$ and $\gamma'$ with $\gamma\neq\gamma'$, we deduce from Theorem 3.3 in \cite{Hcur}, that the graphs must be different.

Now, suppose it is true for $g-1$ characteristic exponents, and we will prove it for $g$. From Proposition \ref{Cgeom} we deduce that is sufficient to prove that the graphs associated to the sets  $\{\gamma_1,\ldots,\gamma_{g-1},\gamma_g\}$ and $\{\gamma_1,\ldots,\gamma_{g-1},\gamma_g'\}$ are different, since otherwise it holds by induction hypothesis. Moreover, since we read the integers $n_1,\ldots,n_g$ in the graph, we assume that $n_g'=n_g$. As in the case $g=1$, by looking at the singular locus (which is seen at $m=1$) we just have to consider the case
$\gamma_g^{(2)}=\gamma_g'^{(2)}=\frac{1}{n_g}$ and the case $\gamma_g^{(2)},\gamma_g'^{(2)}>\frac{1}{n_g}$. In the first case $\gamma_g^{(1)}\neq\gamma_g'^{(1)}$ and $\gamma_i^{(2)}=\gamma_i'^{(2)}=0$ for $1\leq i\leq g-1$. Therefore the graphs are the same till we get to level $m=\mbox{min }\{n_g\langle\nu,\gamma_g\rangle,n_g\langle\nu,\gamma_g'\rangle\}$, where $\nu=(\nu_1,0)\in\s_{Sing}\cap N_{g-1}$ with $\nu_1$ smallest with this property. Since $\langle\nu,\gamma_g\rangle\neq\langle\nu,\gamma_g'\rangle$ the graphs must differ at some moment.
Finally, when $\gamma\neq\gamma'$ with $\gamma_g^{(2)},\gamma_g'^{(2)}>\frac{1}{n_g}$, again by Proposition \ref{Cgeom}, the graphs must be the same for $\{\gamma_1,\ldots,\gamma_g\}$ and $\{\gamma_1,\ldots,\gamma_{g-1},\gamma_g'\}$, till the last approximated root, that is, $f$, starts playing a role in the definition of a component, say $C^\nu$. Since $\langle\nu,\gamma_g\rangle\neq\langle\nu,\gamma_g'\rangle$ we will see the difference on the graphs at level $m=\mbox{min }\{n_g\langle\nu,\gamma_g\rangle,n_g\langle\nu,\gamma_g'\rangle\}$.
\hfill$\Box$

\subsection{Log-canonical threshold}
In \cite{Mus}, Musta\c ta gave a formula of the log-canonical threshold in terms of the codimension of jet schemes, which in our setting can be stated as
\begin{equation}
lct(f)=\mbox{min}_{m\geq 0}\frac{\mbox{Codim}(X_m)}{m+1}.
\label{lct}
\end{equation}

Then, as an application to Theorem \ref{TheCaso1}, we can recover, for the case of surfaces, the result in \cite{BGG}.

\begin{Cor}
The log-canonical threshold of a quasi-ordinary surface singularity is given by:

\[lct_0(X,\A^3)=\left\{\begin{array}{cl}
\frac{1+\l_1^{(1)}}{e_0\l_1^{(1)}} & \ \mbox{ if }\l_1\neq(\frac{1}{n_1},\frac{1}{n_1})\\
\\
1 & \mbox{ if }\lambda_1=(\frac{1}{n_1},\frac{1}{n_1})\mbox{ and }g=1\\
\\
\frac{n_1(1+\l_2^{(1)})}{e_1(n_1(1+\l_2^{(1)})-1)} & \ \mbox{ if }\l_1=(\frac{1}{n_1},\frac{1}{n_1})\mbox{ and }g>1\\
\end{array}\right.\]
Moreover, the components that contribute to the log canonical threshold are
\[\begin{array}{ll}
C_{l_1(\nu)-1}^\nu & \mbox{ if }\gamma_1\neq(\frac{1}{n_1},\frac{1}{n_1})\mbox{ or }g=1\\
\\
C_{l_2(\nu)-1}^\nu & \mbox{ otherwise}\\
\end{array}\]
where $\nu=(l,0)\in N_1$ if $\gamma_1\neq (\frac{1}{n_1},\frac{1}{n_1})$ and $\nu=(l,0)\in N_2$ otherwise.
\label{Thlct}
\end{Cor}

{\em Proof.}
The case $\lambda_1=(\frac{1}{n_1},\frac{1}{n_1})$ and $g=1$ behaves as an $A_n$-singularity, and then $lct(f)=1$. For the rest of the cases, by Corollary \ref{CorCodim}, the codimension of a component grows faster as $m$ grows, for
bigger $j(m,\nu)$. Therefore, the smaller codimension  will be attached for $\nu\in P_m\cap F_m$, and more concretely for
$\nu\in H_m\cap F_m$ whenever $H_m\cap F_m\neq\emptyset$. If $g_1=0$, since $a_1\geq b_1$, we deduce that the minimal codimension among the elements in $P_m\cap F_m$
is attached for $\nu$ of the form $\nu=(l,0)$, while if $g_1>0$ then $P_m\cap F_m$ consists of just a point of the form $\nu=(l,0)$.

We want to minimize not just the codimension, but the quotient $\frac{\mbox{Codim}(X_m)}{m+1}$. That is, to find the biggest $m$ such that $\nu$ still
belongs to $P_m\cap F_m$. Then, when the first characteristic exponent is different from $(\frac{1}{n_1},\frac{1}{n_1})$, this is
attached for $m=l_1(\nu)-1$ such that $\nu\in L_{m+1}^=$. Then $m=l_1(l,0)-1$ and Codim$(C_m^\nu)=l+[\frac{m}{n}]+1$, and since $\nu\in L_{m+1}^=$, $(l,0)\in N_1$
and therefore Codim$(C_m^\nu)=l+l\frac{a_1}{n_1}$, which implies that $\frac{\mbox{Codim}(C_m^\nu)}{m+1}=\frac{a_1+n_1}{na_1}$.

If $\gamma_1=(\frac{1}{n_1},\frac{1}{n_1})$ and $g>1$, what happens is that when $m=l_1(\nu)$ there is no subdivision of the component and $\s_{Reg,1}=\r_1\cup\r_2$.
If we denote the second exponent by $\gamma_2=(\frac{\alpha_2}{n_1n_2},\frac{\beta_2}{n_1n_2})$, we look for $\nu$ of the form $(l,0)$ such that $m+1=l_2(\nu)$
with $\nu\in N_2$. Then $\mbox{Codim}(C_m^\nu)=l+\frac{l_1(\nu)}{n}+[\frac{m-l_1(\nu)}{e_1}]+1=l+\frac{l_1(\nu)}{n}+\frac{l_2(\nu)-l_1(\nu)}{e_1}$, and therefore
$\frac{\mbox{Codim}(C_m^\nu)}{m+1}=\frac{l+l\frac{1}{n_1}+\frac{1}{e_1}(e_2n_2l\frac{\alpha_2}{n_1n_2}-e_1n_1l\frac{1}{n_1})}{e_2n_2l\frac{\alpha_2}{n_1n_2}}=
\frac{1+\frac{1}{n_1}+\frac{\alpha_2}{n_1n_2}-1}{e_2\frac{\alpha_2}{n_1}}=\frac{1+\frac{\alpha_2}{n_2}}{e_1\frac{\alpha_2}{n_2}}$.
This coincides with the statement since $\l_2=(\frac{\alpha_2}{n_1n_2}-\frac{n_1-1}{n_1},\frac{\beta_2}{n_1n_2}-\frac{n_1-1}{n_1})$.
\hfill $\Box$

\vspace{3mm}

\begin{Rem}
Notice that $\frac{1+\lambda_1^{(1)}}{e_0\lambda_1^{(1)}}\leq 1$ except in the case $\lambda_1=(\frac{1}{n_1},\frac{1}{n_1})$ and $g=1$. Moreover in the case of
surfaces the condition $\l_1^{(1)}=\frac{1}{n_1}$ is equivalent to $\l_1=(\frac{1}{n_1},\frac{1}{n_1})$ since the branch is normalized. Then, see notations in
\cite{BGG}, $\ell_1=\ell_2$ and in Theorem \ref{Thlct} we recover, for the case of surfaces, the formula given in Theorem 3.1 in \cite{BGG}.
\end{Rem}

We now deduce a family of examples whose log canonical threshold can not be computed by a monomial valuation.

\begin{Cor}
Let $X$ be a quasi-ordinary surface singularity with $g>1$ characteristic exponents, and such that $\lambda_1=(\frac{1}{n_1},\frac{1}{n_1})$. Then $lct(X,\A^3)$ can not be contributed by monomial valuations in
any variables.
\label{CorMV}
\end{Cor}

{\em Proof.} It follows from Corollary \ref{Thlct} that $lct(X,\A^3)$ is contributed by $C_{l_1(\nu)}^\nu$, for $\nu$ as is made precise in the above statement.
 This is equivalent to say that the valuation
 \[\begin{array}{ll}
 \mathcal V_{C_{l_2(\nu)-1}^\nu} :&  \C[[x_1,x_2,z]]\longrightarrow\N\\
 \\
  & \ \ \ \ \ \ h\ \longmapsto\ \mbox{ord}_t(h\circ\eta)\\
  \end{array}\]
where $\eta$ is the generic point of $(\Psi_{l_2(\nu)-1}^{\A^3})^{-1}(C_{l_2(\nu)-1}^\nu)$. Note that $\nu$ can take all the values described in Corollary \ref{Thlct} but since ${z^{(\langle\nu,\gamma_1\rangle)}}^{n_1}-x_1^{(\nu_1)}x_2^{(\nu_2)}=0$ is one of the defining equations of $C_{l_2(\nu)-1}^\nu$, then $\mathcal V_{C_{l_2(\nu)-1}^\nu}(z^{n_1}-x_1x_2)>n_1\mathcal V_{C_{l_2(\nu)-1}^\nu}(z)$ and $\mathcal V_{C_{l_2(\nu)-1}^\nu}(z^{n_1}-x_1x_2)>\mathcal V_{C_{l_2(\nu)-1}^\nu}(x_1)+\mathcal V_{C_{l_2(\nu)-1}^\nu}(x_2)$. Therefore $\mathcal V_{C_{l_2(\nu)-1}^\nu}$ is not a monomial valuation.
\hfill$\Box$

\subsection{Example}
\label{exs}
Consider the quasi-ordinary surface $f=((z^2-x_1^3)^2-x_1^7x_2^3)^2-x_1^{11}x_2^5(z^2-x_1^3)$. The semigroup is generated by the vectors
\[\gamma_1=\left(\frac{3}{2},0\right),\ \gamma_2=\left(\frac{7}{2},\frac{3}{2}\right)\mbox{ and }\gamma_3=\left(\frac{29}{4},\frac{13}{4}\right).\]
We have that $g_1=g_2=1$.
The singular locus is reducible, of the form
\[X_{Sing}=\{z=x_1=0\}\cup\{x_2=z^2-x_1^3=0\}=Z_1\cup Z_2.\]
Then $\s_{Sing}=\R_{\geq 0}^2\setminus\{0\}$ and $\s_{Reg,1}=\r_2,\ \s_{Reg,2}=\s_{Reg,3}=\{(0,0)\}$.

The set $F_m$ describing the irreducible components is the following, for some $m$:
\[\begin{array}{l}
F_m=\{(1,0),(0,1)\}, \ \mbox{ for }1\leq m<6\\
\\
F_m=\{(1,0),(0,2)\}, \ \mbox{ for }6\leq m<12\\
\\
F_{12}=\{(2,0),(0,2)\}\\
\\
F_{13}=\{(2,0),(0,3)\} \\
\\
F_{18}=\{(2,0),(0,4)\} \\
\\
F_{26}=\{(2,0),(0,4),(0,5)\}  \\
\\
F_{28}=\{(3,0),(2,0),(0,4),(0,5)\} \\
\end{array}\]
and the result can be checked by lifting the components $Z_1$ and $Z_2$ of the singular locus to level $m$ as the following graph shows (we did not draw the weights of the vertices for clearness).

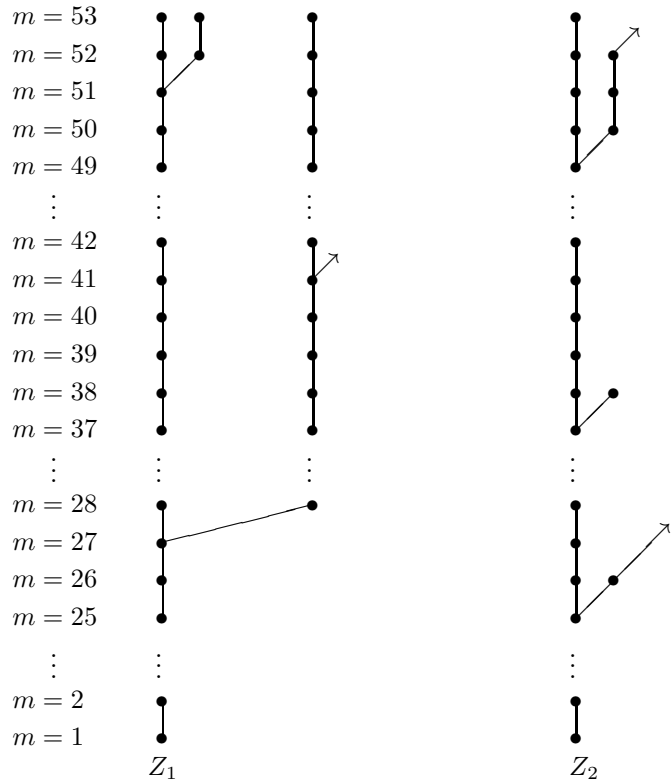
\begin{figure}[h]
\unitlength=1mm
\begin{center}
\begin{picture}(110,100)(0,3)
\linethickness{0.15mm}

\put(24,0){$Z_1$}\put(80,0){$Z_2$}

\put(6,4){$m=1$}\put(25,4){$\bullet$}\put(80,4){$\bullet$}
\put(26,5){\line(0,1){5}}\put(81,5){\line(0,1){5}}

\put(6,9){$m=2$}\put(25,9){$\bullet$}\put(80,9){$\bullet$}

\put(11,13){$\vdots$}\put(25,13){$\vdots$}\put(80,13){$\vdots$}

\put(6,20){$m=25$}\put(25,20){$\bullet$}\put(80,20){$\bullet$}
\put(26,21){\line(0,1){5}}\put(81,21){\line(0,1){5}}\put(81,21){\line(1,1){5}}

\put(6,25){$m=26$}\put(25,25){$\bullet$}\put(80,25){$\bullet$}\put(85,25){$\bullet$}
\put(26,26){\line(0,1){5}}\put(81,26){\line(0,1){5}}

\put(6,30){$m=27$}\put(25,30){$\bullet$}\put(80,30){$\bullet$}\put(86,26){\line(1,1){5}}\put(90,31){$\nearrow$}
\put(26,31){\line(0,1){5}}\put(26,31){\line(4,1){20}}\put(81,31){\line(0,1){5}}

\put(6,35){$m=28$}\put(25,35){$\bullet$}\put(45,35){$\bullet$}\put(80,35){$\bullet$}

\put(11,39){$\vdots$}\put(25,39){$\vdots$}\put(45,39){$\vdots$}\put(80,39){$\vdots$}

\put(6,45){$m=37$}\put(25,45){$\bullet$}\put(45,45){$\bullet$}\put(80,45){$\bullet$}
\put(26,46){\line(0,1){5}}\put(46,46){\line(0,1){5}}\put(81,46){\line(0,1){5}}\put(81,46){\line(1,1){5}}

\put(6,50){$m=38$}\put(25,50){$\bullet$}\put(45,50){$\bullet$}\put(80,50){$\bullet$}\put(85,50){$\bullet$}
\put(26,51){\line(0,1){5}}\put(46,51){\line(0,1){5}}\put(81,51){\line(0,1){5}}

\put(6,55){$m=39$}\put(25,55){$\bullet$}\put(45,55){$\bullet$}\put(80,55){$\bullet$}
\put(26,56){\line(0,1){5}}\put(46,56){\line(0,1){5}}\put(81,56){\line(0,1){5}}

\put(6,60){$m=40$}\put(25,60){$\bullet$}\put(45,60){$\bullet$}\put(80,60){$\bullet$}
\put(26,61){\line(0,1){5}}\put(46,61){\line(0,1){5}}\put(81,61){\line(0,1){5}}

\put(6,65){$m=41$}\put(25,65){$\bullet$}\put(45,65){$\bullet$}\put(80,65){$\bullet$}
\put(26,66){\line(0,1){5}}\put(46,66){\line(0,1){5}}\put(46,67){$\nearrow$}\put(81,66){\line(0,1){5}}

\put(6,70){$m=42$}\put(25,70){$\bullet$}\put(45,70){$\bullet$}\put(80,70){$\bullet$}

\put(11,74){$\vdots$}\put(25,74){$\vdots$}\put(45,74){$\vdots$}\put(80,74){$\vdots$}

\put(6,80){$m=49$}\put(25,80){$\bullet$}\put(45,80){$\bullet$}\put(80,80){$\bullet$}
\put(26,81){\line(0,1){5}}\put(46,81){\line(0,1){5}}\put(81,81){\line(0,1){5}}\put(81,81){\line(1,1){5}}

\put(6,85){$m=50$}\put(25,85){$\bullet$}\put(45,85){$\bullet$}\put(80,85){$\bullet$}\put(85,85){$\bullet$}
\put(26,86){\line(0,1){5}}\put(46,86){\line(0,1){5}}\put(81,86){\line(0,1){5}}\put(86,86){\line(0,1){5}}

\put(6,90){$m=51$}\put(25,90){$\bullet$}\put(45,90){$\bullet$}\put(80,90){$\bullet$}\put(85,90){$\bullet$}
\put(26,91){\line(0,1){5}}\put(26,91){\line(1,1){5}}\put(46,91){\line(0,1){5}}\put(81,91){\line(0,1){5}}\put(86,91){\line(0,1){5}}

\put(6,95){$m=52$}\put(25,95){$\bullet$}\put(30,95){$\bullet$}\put(45,95){$\bullet$}\put(80,95){$\bullet$}\put(85,95){$\bullet$}
\put(26,96){\line(0,1){5}}\put(31,96){\line(0,1){5}}\put(46,96){\line(0,1){5}}\put(81,96){\line(0,1){5}}\put(86,97){$\nearrow$}

\put(6,100){$m=53$}\put(25,100){$\bullet$}\put(30,100){$\bullet$}\put(45,100){$\bullet$}\put(80,100){$\bullet$}

\end{picture}
\end{center}
\caption{The graph of the surface defined by $f=((z^2-x_1^3)^2-x_1^7x_2^3)^2-x_1^{11}x_2^5(z^2-x_1^3)$.}
\label{figEx1}
\end{figure}

The arrows in the figure represent the behaviour explained in Remark \ref{RemarkFinal}.

\vspace{3mm}

Now we give some explanations to illustrate how Proposition \ref{Lema2Gen} and Proposition \ref{PropC1} work.

 For $m=1$, we have $H_1=\{(1,0),(0,1)\}$, $L_1^=\{(0,1)\}$ and $L_1^<=\emptyset$. The claim on $F_1$ in this case follows easily by Proposition \ref{Lema2Gen}.

At level $m=6$ we have $H_6=\{\nu\in[0,m]^2\cap N_0\ |\ \nu_1\neq 0\}$, $L_6^==\emptyset$ and $L_6^<=\{(0,\nu_2)\ |\ 2\leq\nu_2\leq 6\}$ with $j'(6,(0,2))=1$, then, by Proposition \ref{PropC1} $C_6^{\nu'}\subseteq C_6^{(0,2)}$ for $\nu'\in\{(0,3),(0,4),(0,5),(0,6)\}$, because $\s_{Reg,1}=\r_2=\{0\}\times\R_{\geq 0}$. By Proposition \ref{Lema2Gen}, only $\nu=(1,0)$ contributes to $F_6$ from the vectors in $H_6$, and the claim on $F_6$ follows. Note how at this level $\nu=(0,1)$ does no longer give rise to an irreducible component, since $l_2(0,1)=6$ and $(0,1)\notin N_2$. Then we have that $(0,2)\in F_6$ and the vertex associated with $C_5^{(0,1)}$ and the one associated with $C_6^{(0,2)}$ are joined in the graph $\Gamma'$.

\section{Technical results and proofs.}
\label{Proofs}
In this section we state and prove some results which are used along the paper but only in the proofs of other results, and can be skipped to read Section \ref{Sec4}. Moreover there are some proofs which we leave to this section.

\vspace{3mm}

Recall that we denote the first characteristic exponent by $\lambda_1=\gamma_1=(\frac{a_1}{n_1},\frac{b_1}{n_1})$ with $a_1\geq b_1,\ a_1>0\mbox{ and }b_1\geq 0$,
and if $b_1=0$ then $a_1>n_1$ (we always consider normalized branches). If $g_1>0$ then, by Lemma \ref{Lema35}, we have that $n_{g_1+1}\gamma_{g_1+1}=
(\alpha_1^{(g_1+1)},\alpha_2^{(g_1+1)})+(r_1^{(g_1+1)}\gamma_1^{(1)}+\cdots+r_{g_1}^{(g_1+1)}\gamma_{g_1}^{(1)},0)$, therefore $n_{g_1+1}\gamma_{g_1+1}^{(2)}\in\Z$, or in other words,  $\gamma_{g_1+1}=(\gamma_{g_1+1}^{(1)},\frac{b_{g_1+1}}{n_{g_1+1}})$
with $b_{g_1+1}\geq 1$.

\begin{Lem}
In the relation $n_{g_2+1}\gamma_{g_2+1}=(\alpha_1^{(g_2+1)},\alpha_2^{(g_2+1)})+r_1^{(g_2+1)}\gamma_1+\cdots+r_{g_2}^{(g_2+1)}\gamma_{g_2}$ given in Lemma
\ref{LemaPedro}, we have that $\alpha_2^{(g_2+1)}>1$.
\label{alphaM1}
\end{Lem}

{\em Proof.} If $g_2=g_1$ the claim is trivial since $\alpha_2^{(g_2+1)}=b_{g_2+1}>1$. Otherwise $n_{g_2+1}\gamma_{g_2+1}^{(2)}=\alpha_2^{(g_2+1)}+
r_{g_2}^{(g_2+1)}\frac{1}{n_{g_2}}$, and since, by Lemma \ref{LemaPedro}, $\gamma_{g_2+1}^{(2)}\geq n_{g_2}\gamma_{g_2}^{(2)}=1$ and $0\leq r_{g_2}^{(g_2+1)}<
n_{g_2}$, then $\alpha_2^{(g_2+1)}\geq n_{g_2+1}-\frac{r_{g_2}^{(g_2+1)}}{n_{g_2}}>1$, because $n_{g_2+1}\geq 2$. \hfill $\Box$

\vspace{3mm}

In Corollary \ref{Corolario3} we describe the generators of $J_m^\nu$ for $\nu\in H_m\cup L_m$. But we also need to describe the polynomial $F_{\nu}^{(l_{i(\nu)}(\nu))}$ (recall that by definition $\nu\notin H_{l_{i(\nu)}(\nu)}\cup L_{l_{i(\nu)}(\nu)}$).
 We do this in the next Lemma, but before we look at an example.

\begin{Exam}
Let $X$ be a quasi-ordinary surface defined by $f=((z^2-x_1^3x_2^2)^2-x_1^6x_2^4z)^3-x_1^{23}x_2^{14}z$. The generators of the semigroup are $\gamma_1=(\frac{3}{2},1),\ \gamma_2=(\frac{15}{4},\frac{5}{2}) \mbox{ and }\gamma_3=(\frac{49}{6},5)$. Notice that $\nu=(0,3)\notin N_2$, and $l_2(\nu)=l_3(\nu)$. At level $m=45$ we have the set
\[D_{45}^{(0,3)}=V(x_2^{(0)},x_2^{(1)},x_2^{(2)},z^{(0)},z^{(1)},z^{(2)},F_{1,\nu}^{(6)},F_{1,\nu}^{(7)},F_{3,\nu}^{(45)})\cap D(x_1^{(0)})\cap D(x_2^{(3)}),\]
where
\[\begin{array}{rl}
F_{3,\nu}^{(45)} & ={F_{2,\nu}^{(15)}}^3-{x_1^{(0)}}^{23}{x_2^{(3)}}^{14}z^{(3)}\\
 & =({x_1^{(0)}}^6{x_2^{(3)}}^4z^{(3)})^3-{x_1^{(0)}}^{23}{x_2^{(3)}}^{14}z^{(3)}\\
 = & ({x_1^{(0)}}^6{x_2^{(3)}}^4z^{(3)})^3(1-\frac{{x_1^{(0)}}^8{x_2^{(3)}}^3}{{z^{(3)}}^2}),\\
 \end{array}\]
 since $D_{45}^{(0,3)}\subset D(z^{(3)})$. Since $\gamma(0)=(x_1^{(0)},0,0)\in X$, and we are considering germs of quasi-ordinary singularities, we have that $|x_1^{(0)}|<<1$ and we deduce that $1-\frac{{x_1^{(0)}}^8{x_2^{(3)}}^3}{{z^{(3)}}^2}\neq 0$.
\label{Ex4}
\end{Exam}

\vspace{3mm}

This example illustrates the fact that we are looking at jet schemes of a germ of quasi-ordinary singularity, instead of jet schemes of the whole affine surface. If we looked at the whole surface there would be other irreducible components that we do not consider here. This is expectable because the components we consider are determined by the invariants of the topological type at the origin, so they describe only what happens in a small neighbourhood of zero. Actually the other components that may appear when looking at the whole affine surface, will project on closed points, different from the origin, of the singular locus.

\begin{Lem}
Given $m\in\Z_{>0}$ and $\nu\in H_m\cup L_m$ with $m+1=l_{i(\nu)}(\nu)$, then
\[F_\nu^{(l_{i(\nu)}(\nu))}= \left({x_1^{(\nu_1)}}^{\alpha_1^{(i(\nu))}}{x_2^{(\nu_2)}}^{\alpha_2^{(i(\nu))}}{F_{0,\nu}^{(\frac{l_1(\nu)}{e_0})}}^{r_1^{(i(\nu))}}{F_{1,\nu}^{(\frac{l_2(\nu)}{e_1})}}
^{r_2^{(i(\nu))}}\cdots{F_{i(\nu)-2,\nu}^{(\frac{l_{i(\nu)-1}(\nu)}{e_{i(\nu)-2}})}}^{r_{i(\nu)-1}^{(i(\nu))}}\right)^{e_{i(\nu)}}\cdot U,\]
where $U$ is a unit in $\C[[{x_1^{(\nu_1)}}^{\pm 1},{x_2^{(\nu_2)}}^{\pm 1},{F_0^{(\frac{l_1(\nu)}{e_0})}}^{\pm 1},\ldots,{F_{i(\nu)-2}^{(\frac{l_{i(\nu)-1}(\nu)}{e_{i(\nu)-2}})}}^{\pm 1}]]$.
When $\nu\notin\r_1\cup\r_2$, then $U=1$.
\label{LemF'}
\end{Lem}

{\em Proof.} We have that $j(m,\nu)=i(\nu)-1$, and, by Lemma \ref{Lemfk}, for any $\gamma\in D_m^\nu$
\[\begin{array}{ll}
\mbox{ord}_t(f_k\circ\gamma)=\langle\nu,\gamma_{k+1}\rangle & \mbox{ for }0\leq k\leq i(\nu)-2\\
\\
\mbox{ord}_t(f_{i(\nu)-1}\circ\gamma)>\frac{m}{e_{i(\nu)-1}}\\
\end{array}\]
Then ord$_t(f_{i(\nu)-1}\circ\gamma)\geq\frac{m+1}{e_{i(\nu)-1}}=\frac{l_{i(\nu)}(\nu)}{e_{i(\nu)-1}}=\langle\nu,\gamma_{i(\nu)}\rangle$, and since $\nu\notin N_{i(\nu)}$, $\langle\nu,\gamma_{i(\nu)}\rangle$ is not an integer. Hence
\[\mbox{ord}_t(f_{i(\nu)-1}\circ\gamma)>\langle\nu,\gamma_{i(\nu)}\rangle.\]
We have by Lemma \ref{Lema35}
\[f_{i(\nu)}=f_{i(\nu)-1}^{n_{i(\nu)}}-c_{i(\nu)}x_1^{\alpha_1^{(i(\nu))}}x_2^{\alpha_2^{(i(\nu))}}f_0^{r_1^{(i(\nu))}}\cdots f_{i(\nu)-2}^{r_{i(\nu)-1}^{(i(\nu))}}+\sum c_{\underline\alpha,\underline{r}}x_1^{\alpha_1}x_2^{\alpha_2}f_0^{r_1}\cdots f_{i(\nu)-1}^{r_{i(\nu)}},\]
and
\[\begin{array}{rl}
\mbox{ord}_t(f_{i(\nu)-1}^{n_{i(\nu)}}\circ\gamma) & >n_{i(\nu)}\langle\nu,\gamma_{i(\nu)}\rangle\\
\\
\mbox{ord}_t((c_{i(\nu)}x_1^{\alpha_1^{(i(\nu))}}x_2^{\alpha_2^{(i(\nu))}}\cdots f_{i(\nu)-2}^{r_{i(\nu)-1}^{(i(\nu))}})\circ\gamma) & =n_{i(\nu)}\langle\nu,\gamma_{i(\nu)}\rangle\\
\\
\mbox{ord}_t((c_{\underline\alpha,\underline{r}}x_1^{\alpha_1}x_2^{\alpha_2}\cdots f_{i(\nu)-1}^{r_{i(\nu)}})\circ\gamma) &
=\langle\nu,(\alpha_1,\alpha_2)+r_1\gamma_1+\cdots+r_{i(\nu)-1}\gamma_{i(\nu)-1}\rangle+r_{i(\nu)}\mbox{ord}_t(f_{i(\nu)-1}\circ\gamma)\\
 & >\langle\nu,(\alpha_1,\alpha_2)+r_1\gamma_1+\cdots+r_{i(\nu)}\gamma_{i(\nu)}\rangle\\
& \geq n_{i(\nu)}\langle\nu,\gamma_{i(\nu)}\rangle.\\
\end{array}\]
Then ord$_t(f_{i(\nu)}\circ\gamma)=n_{i(\nu)}\langle\nu,\gamma_{i(\nu)}\rangle=\frac{l_{i(\nu)}(\nu)}{e_{i(\nu)}}$, and
\[F_{i(\nu),\nu}^{(\frac{l_{i(\nu)}(\nu)}{e_{i(\nu)}})}= {x_1^{(\nu_1)}}^{\alpha_1^{(i(\nu))}}{x_2^{(\nu_2)}}^{\alpha_2^{(i(\nu))}}{z^{(\langle\nu,\gamma_1\rangle)}}^{r_1^{(i(\nu))}}{F_{1,\nu}^{(\frac{l_2(\nu)}{e_1})}}
^{r_2^{(i(\nu))}}\cdots{F_{i(\nu)-2,\nu}^{(\frac{l_{i(\nu)-1}(\nu)}{e_{i(\nu)-2}})}}^{r_{i(\nu)-1}^{(i(\nu))}}\]

 By Lemma \ref{LemExpSR}
\[f=f_{i(\nu)}^{e_{i(\nu)}}+\sum c_{ijk}^{(i(\nu))}x_1^ix_2^jz^k\]
where $(i,j)+k\gamma_1>n_{i(\nu)}e_{i(\nu)}\gamma_{i(\nu)}$. Then
\[\begin{array}{ll}
\mbox{ord}_t(f_{i(\nu)}^{e_{i(\nu)}}\circ\gamma)=l_{i(\nu)}(\nu)\\
\\
\mbox{ord}_t((c_{ijk}^{(i(\nu))}x_1^ix_2^jz^k)\circ\gamma)\geq l_{i(\nu)}(\nu)\\
\end{array}\]
and hence
\[F_\nu^{(l_{i(\nu)}(\nu))}={F_{i(\nu),\nu}^{(\frac{l_{i(\nu)}(\nu)}{e_{i(\nu)}})}}^{e_{i(\nu)}}+G_{i(\nu),\nu},\]
where
\[G_{i(\nu),\nu}=\sum c_{ijk}^{(i(\nu))}{x_1^{(\nu_1)}}^i{x_2^{(\nu_2)}}^j{z^{(\langle\nu,\gamma_1\rangle)}}^k\]
and the sum runs over $i,j,k$ such that
\begin{enumerate}
\item[(i)] $c_{ijk}^{(i(\nu))}\neq 0$
\item[(ii)] $\langle\nu,(i,j)+k\gamma_1\rangle=l_{i(\nu)}(\nu)$
\end{enumerate}
Notice that if $\nu\notin\r_1\cup\r_2$, then condition (ii) never holds and $G_{\nu,i(\nu)}=0$. In this case we are done. Otherwise,  from (\ref{eqO}) we deduce that $\gamma\in D_m^\nu\subset D(F_{\nu,i(\nu)}^{(\frac{l_{i(\nu)}(\nu)}{e_{i(\nu)}})})$, and hence
\[F_\nu^{(l_{i(\nu)}(\nu))}={F_{i(\nu),\nu}^{(\frac{l_{i(\nu)}(\nu)}{e_{i(\nu)}})}}^{e_{i(\nu)}}\left(1+G_{i(\nu),\nu}/{F_{i(\nu,\nu)}^{(\frac{l_{i(\nu)}(\nu)}{e_{i(\nu)}})}}
^{e_{i(\nu)}}\right).\]
\hfill$\Box$

\vspace{3mm}

{\bf\em Proof of Lemma \ref{LemD}.}
 If $m<l_{i(\nu)}(\nu)$, then we deduce from Proposition \ref{Prop1} that $D_m^\nu$ is non-empty.
Otherwise $m\geq l_{i(\nu)}(\nu)$, and by definition $F^{(l_{i(\nu)}(\nu))}\in J_m^\nu$, and by Lemma \ref{LemF'} (and its proof)
\[F_\nu^{(l_{i(\nu)}(\nu))}={F_{i_(\nu),\nu}^{(\frac{l_{i(\nu)}(\nu)}{e_{i(\nu)}})}}^{e_{i(\nu)}}\cdot U,\]
where
\[F_{i_(\nu),\nu}^{(\frac{l_{i(\nu)}(\nu)}{e_{i(\nu)}})}={x_1^{(\nu_1)}}^{\alpha_1^{(i(\nu))}}{x_2^{(\nu_2)}}^{\alpha_2^{(i(\nu))}}{F_{0,\nu}^{(\frac{l_1(\nu)}{e_0})}}^{r_1^{(i(\nu))}}{F_{1,\nu}^{(\frac{l_2(\nu)}{e_1})}}
^{r_2^{(i(\nu))}}\cdots{F_{i(\nu)-2,\nu}^{(\frac{l_{i(\nu)-1}(\nu)}{e_{i(\nu)-2}})}}^{r_{i(\nu)-1}^{(i(\nu))}}\]
with $(\alpha_1^{(i(\nu))},\alpha_2^{(i(\nu))})+r_1^{(i(\nu))}\gamma_1+\cdots+r_{i(\nu)-1}^{(i(\nu))}\gamma_{i(\nu)-1}=n_{i(\nu)}\gamma_{i(\nu)}$. And
\[U=1+G_{i(\nu),\nu}/{F_{i(\nu),\nu}^{(\frac{l_{i(\nu)}(\nu)}{e_{i(\nu)}})}}^{e_{i(\nu)}},\]
where
\[G_{i(\nu),\nu}=\sum c_{ijk}^{(i(\nu))}{x_1^{(\nu_1)}}^i{x_2^{(\nu_2)}}^j{z^{(\langle\nu,\gamma_1\rangle)}}^k\]
with $\langle\nu,(i,j)+k\gamma_1\rangle=n_{i(\nu)}e_{i(\nu)}\langle\nu,\gamma_{i(\nu)}\rangle$, though $(i,j)+k\gamma_1>n_{i(\nu)}e_{i(\nu)}\gamma_{i(\nu)}$.

Notice that at level $l_{i(\nu)}(\nu)-1$ we can apply Corollary \ref{Corolario3}, and deduce, as in (\ref{eqO}), that
\[D_{l_{i(\nu)-1}(\nu)-1}^\nu\subseteq D(F_{0,\nu}^{(\frac{l_1(\nu)}{e_0})}\cdots F_{i(\nu)-2,\nu}^{(\frac{l_{i(\nu)-1}(\nu)}{e_{i(\nu)-2}})})\]
Hence $D_{l_{i(\nu)}(\nu)}^\nu$ satisfies the same property, and since, by definition, $D_{l_{i(\nu)}(\nu)}^\nu\subseteq D(x_1^{(\nu_1)}\cdot x_2^{(\nu_2)})$, we just have to argue that $U\neq 0$.

In Lemma \ref{LemF'} we prove that $U=1$ when $\nu\notin\r_1\cup\r_2$. Suppose the contrary, then either $\nu_1=0$ or $\nu_2=0$. We work it all out and write $F_{i(\nu),\nu}^{(\frac{l_{i(\nu)}(\nu)}{e_{i(\nu)}})}$ in terms of $x_1^{(\nu_1)}$, $x_2^{(\nu_2)}$ and $z^{(\langle\nu,\gamma_1\rangle)}$. The key point is that when $G_{i(\nu),\nu}\neq 0$ is because $(i,j)+k\gamma_1>n_{i(\nu)}e_{i(\nu)}\gamma_{i(\nu)}$ but $\langle\nu,(i,j)+k\gamma_1\rangle=n_{i(\nu)}e_{i(\nu)}\langle\nu,\gamma_{i(\nu)}\rangle$, because $\nu\in\r_1\cup\r_2$. Therefore either $x_1^{(0)}$ or $x_2^{(0)}$ do appear in U (it do not cancel in the quotient $G_{i(\nu),\nu}/{F_{i(\nu),\nu}^{(\frac{l_{i(\nu)}(\nu)}{e_{i(\nu)}})}}^{e_{i(\nu)}}$), depending on whether $\nu_1=0$ or $\nu_2=0$. Then we have that $U$ depends on the origin of the jet $\gamma(0)\in X$, and since we are dealing with germs of quasi-ordinary singularities, we have that $|x_i^{(0)}|<<1$ for $i=1,2$, and we deduce that $U\neq 0$.
\hfill$\Box$

\vspace{3mm}

{\bf\em Proof of Lemma \ref{Lemfk}.} We distinguish the cases $\nu\in H_m$ and $\nu\in L_m$. For $\nu\in H_m$, $j(m,\nu)=0$ and we have to prove that $\mbox{ord}_t(f_k\circ\gamma)>\frac{m}{e_k}$ for $0\leq k\leq g$. By Proposition \ref{CHm} it is true for $k=0$. For $k=g$ the claim is obvious, and for $1\leq k\leq g-1$ we use the expansion (\ref{eqSRpqprima}) in Lemma \ref{LemExpSR},
\[f=f_k^{e_k}+\sum_{(i,j)+r\gamma_1>n_ke_k\gamma_k}c_{ijr}^{(k)}x_1^ix_2^jz^r.\]
Suppose that there exist $c_{ijr}^{(k)}\neq 0$ such that $\langle\nu,(i,j)\rangle+r\mbox{ ord}_t(z\circ\gamma)\leq m$. Then, using that ord$_t(z\circ\gamma)>\frac{m}{n}$ and $k\geq 1$, we have the following inequalities
\[l_1(\nu)-r\langle\nu,\gamma_1\rangle+r\frac{m}{n}<l_1(\nu)-r\langle\nu,\gamma_1\rangle+r\mbox{ ord}_t(z\circ\gamma)\leq l_k(\nu)-r\langle\nu,\gamma_1\rangle+r\mbox{ ord}_t(z\circ\gamma)\leq m.\]
Then $l_1(\nu)-r\langle\nu,\gamma_1\rangle<m(1-\frac{r}{n})$, and since $r<n$ this is equivalent to $\langle\nu,\gamma_1\rangle<\frac{m}{n}$, which contradicts the fact that $j(m,\nu)=0$. Hence ord$_t(c_{ijr}^{(k)}x_1^ix_2^jz^r\circ\gamma)>m$ and therefore ord$_t(f_k^{e_k}\circ\gamma)>m$ as we wanted to prove.

\vspace{2mm}

 For $\nu\in L_m$ the proof is by induction on $j(m,\nu)$. And we will make use repeatedly of the following observation. In general, for any function $f\in\C[x_1,x_2,z]$, and any $m$-jet $\gamma$, there is no relation among ord$_t(f\circ\gamma)$ and ord$_t(f\circ\pi_{m,m'}(\gamma))$ with $m'<m$. But if ord$_t(x_i\circ\gamma)\neq 0$ for $i=1,2$, ord$_t(z\circ\gamma)\neq 0$, and for $m'<m$, the $m'$-jet $\gamma':=\pi_{m,m'}(\gamma)$ verifies that ord$_t(x_i\circ\gamma')\neq 0$ and ord$_t(z\circ\gamma')\neq 0$, then
\[\mbox{ord}_t(f\circ\gamma)=\mbox{ord}_t(f\circ\gamma').\]

By Remark \ref{RemDefHL} (ii) the first case of induction is $j(m,\nu)=1$. Then in particular $l_1(\nu)<l_2(\nu)$. For $\gamma\in D_m^\nu$, set $\bar m:=l_1(\nu)-1<m$ and $\bar\gamma:=\pi_{m,\bar m}(\gamma)$. Then $j(\bar m,\nu)=0$ and $\bar\gamma\in D_{\bar m}^\nu$.  By Proposition \ref{CHm} we have
\[\mbox{ord}_t(z\circ\bar\gamma)>\frac{\bar m}{n}=\frac{l_1(\nu)-1}{n}=\langle\nu,\gamma_1\rangle-\frac{1}{n}\]
Hence ord$_t(z\circ\bar\gamma)\geq\langle\nu,\gamma_1\rangle$ and therefore ord$_t(z\circ\gamma)\geq\langle\nu,\gamma_1\rangle$. Suppose that the inequality is strict, ord$_t(z\circ\gamma)>\langle\nu,\gamma_1\rangle$. By Lemma \ref{LemExpSR}
\[f_1=z^{n_1}-x_1^{a_1}x_2^{b_1}+\sum_{(i_1,i_2)+k\gamma_1>n_1\gamma_1}x_1^{i_1}x_2^{i_2}z^k\]
Then
\[\begin{array}{rl}
\mbox{ord}_t(z^{n_1}\circ\gamma) & >n_1\langle\nu,\gamma_1\rangle\\
\\
\mbox{ord}_t(x_1^{a_1}x_2^{b_1}\circ\gamma) & =n_1\langle\nu,\gamma_1\rangle\\
\\
\mbox{ord}_t(x_1^{i_1}x_2^{i_2}z^k\circ\gamma) & >\langle\nu,(i_1,i_2)+k\gamma_1\rangle\geq n_1\langle\nu,\gamma_1\rangle\\
\end{array}\]
and hence ord$_t(f_1\circ\gamma)=n_1\langle\nu,\gamma_1\rangle$. Again by Lemma \ref{LemExpSR} we have
\[f=f_1^{e_1}+\sum_{(i_1,i_2)+k\gamma_1>n_1e_1\gamma_1}c_{i_1i_2k}x_1^{i_1}x_2^{i_2}z^k\]
and since
\[\begin{array}{rl}
\mbox{ord}_t(f_1^{e_1}\circ\gamma) & =e_1n_1\langle\nu,\gamma_1\rangle=l_1(\nu)\\
\\
\mbox{ord}_t(x_1^{i_1}x_2^{i_2}z^k\circ\gamma) & >\langle\nu,(i_1,i_2)+k\gamma_1\rangle\geq n_1e_1\langle\nu,\gamma_1\rangle=l_1(\nu)\\
\end{array}\]
we deduce ord$_t(f\circ\gamma)=l_1(\nu)\leq m$, which is a contradiction. Then ord$_t(z\circ\gamma)=\langle\nu,\gamma_1\rangle$.

Now we prove that ord$_t(f_1\circ\gamma)>\frac{m}{e_1}$. Suppose the contrary, ord$_t(f_1\circ\gamma)\leq\frac{m}{e_1}$. If $g=1$ there is nothing to prove.
If $g=2$, then we consider the expansion given in Lemma \ref{Lema35}
\[f=f_2=f_1^{n_2}-c_2x_1^{\alpha_1^{(2)}}x_2^{\alpha_2^{(2)}}z^{r_1^{(2)}}+\sum c_{\underline{\alpha},\underline{r}}x_1^{i_1}x_2^{i_2}z^{r_1}f_1^{r_2}\]
where $(i_1,i_2)+r_1\gamma_1+r_2\gamma_2>n_2\gamma_2$, and
\[\begin{array}{rl}
\mbox{ord}_t(f_1^{n_2}\circ\gamma) & \leq m,\\
\\
\mbox{ord}_t(x_1^{\alpha_1^{(2)}}x_2^{\alpha_2^{(2)}}z^{r_1^{(2)}}\circ\gamma) & =n_2\langle\nu,\gamma_2\rangle=l_2(\nu)>m\\
\end{array}\]
Then there must exist $c_{\underline{\alpha},\underline{r}}\neq 0$ such that
\[n_2\mbox{ord}_t(f_1\circ\gamma)=\langle\nu,(i_1,i_2)+r_1\gamma_1\rangle+r_2\mbox{ord}_t(f_1\circ\gamma)\]
or equivalently
\[(n_2-r_2)\mbox{ord}_t(f_1\circ\gamma)=\langle\nu,(i_1,i_2)+r_1\gamma_1\rangle\geq (n_2-r_2)\langle\nu,\gamma_2\rangle\]
And since $r_2<n_2$ we conclude ord$_t(f_1\circ\gamma)\geq\langle\nu,\gamma_2\rangle>\frac{m}{e_1}$, which is a contradiction.

If $g>2$, by Lemma \ref{LemExpSR}
\[f=f_1^{e_1}+\sum c_{\underline{\alpha},\underline{r}}x_1^{\alpha_1}x_2^{\alpha_2}z^{r_1}f_1^{r_2}\]
and since we are supposing that ord$_t(f_1^{e_1}\circ\gamma)\leq m$, there must exists $c_{\underline{\alpha},\underline{r}}\neq 0$ such that
$\mbox{ord}_t(f_1^{e_1}\circ\gamma)=\mbox{ord}_t(x_1^{\alpha_1}x_2^{\alpha_2}z^{r_1}f_1^{r_2}\circ\gamma)$,
hence
\[e_1\mbox{ord}_t(f_1\circ\gamma)=\langle\nu,(\alpha_1,\alpha_2)+r_1\gamma_1\rangle+r_2\mbox{ord}_t(f_1\circ\gamma)\]
or equivalently
\[(e_1-r_2)\mbox{ord}_t(f_1\circ\gamma)=\langle\nu,(\alpha_1,\alpha_2)+r_1\gamma_1\rangle\geq n_2e_2\langle\nu,\gamma_2\rangle-r_2\langle\nu,\gamma_2\rangle\]
and since $r_2<e_1$ we conclude that ord$_t(f_1\circ\gamma)\geq \langle\nu,\gamma_2\rangle$, which is a contradiction.

The rest is simple, by Lemma \ref{LemExpSR}, for $k>1=j(m,\nu)$, we have
\[f=f_k^{e_k}+\sum c_{i_1i_2k}x_1^{i_1}x_2^{i_2}z^k\]
and  since ord$_t(c_{i_1i_2k}x_1^{i_1}x_2^{i_2}z^k\circ\gamma)=\langle\nu,(i_1,i_2)+k\gamma_1\rangle\geq l_k(\nu)>m$, we deduce ord$_t(f_k^{e_k}\circ\gamma)>m$.

Suppose now that the claim is true for $j(m,\nu)=j$ and we will prove it for $j(m,\nu)=j+1$. Let $\gamma\in D_m^\nu$, with
$l_{j+1}(\nu)\leq m<l_{j+1}(\nu)$. We set $\bar m=l_{j+1}(\nu)-1$ and $\bar\gamma=\pi_{m,\bar m}(\gamma)$. Then $\gamma\in D_{\bar m}^\nu$ and $j(\bar m,\nu)=i\leq j$, where
\[l_i(\nu)\leq\bar m<l_{i+1}(\nu)=\cdots=l_{j+1}(\nu).\]
Then, by Lemma \ref{l-ord}, this is equivalent to
\[\begin{array}{c}
n_{i+1}\langle\nu,\gamma_{i+1}\rangle=\langle\nu,\gamma_{i+2}\rangle\\
\\
n_{i+2}\langle\nu,\gamma_{i+2}\rangle=\langle\nu,\gamma_{i+3}\rangle\\
\vdots\\
n_j\langle\nu,\gamma_j\rangle=\langle\nu,\gamma_{j+1}\rangle\\
\end{array}\]

By induction hypothesis we deduce ord$_t(f_k\circ\gamma)=\mbox{ord}_t(f_k\circ\bar\gamma)=\langle\nu,\gamma_{k+1}\rangle$ for $0\leq k<i$. We are going to prove that ord$_t(f_i\circ\gamma)=\langle\nu,\gamma_{i+1}\rangle$. By induction we have
\[\mbox{ord}_t(f_i\circ\gamma)\geq\mbox{ord}_t(f_i\circ\bar\gamma)>\frac{\bar m}{e_i}=\langle\nu,\gamma_{i+1}\rangle-\frac{1}{e_i}\]
Therefore
\[\mbox{ord}_t(f_i\circ\gamma)\geq\langle\nu,\gamma_{i+1}\rangle\]
Let us suppose that ord$_t(f_i\circ\gamma)>\langle\nu,\gamma_{i+1}\rangle$. By Lemma \ref{LemExpSR}
\[f_{i+1}=f_i^{n_{i+1}}-c_{i+1}x_1^{\alpha_1^{(i+1)}}x_2^{\alpha_2^{(i+1)}}z^{r_1^{(i+1)}}\cdots f_{i-1}^{r_i^{(i+1)}}+\sum c_{\underline{\alpha},\underline{r}}x_1^{\alpha_1}x_2^{\alpha_2}z^{r_1}\cdots f_i^{r_{i+1}}\]
and we have
\[\begin{array}{rl}
\mbox{ord}_t(f_i^{n_{i+1}}\circ\gamma) & > n_{i+1}\langle\nu,\gamma_{i+1}\rangle\\
\\
\mbox{ord}_t(x_1^{\alpha_1^{(i+1)}}x_2^{\alpha_2^{(i+1)}}z^{r_1^{(i+1)}}\circ\gamma) & =n_{i+1}\langle\nu,\gamma_{i+1}\rangle\\
\\
\mbox{ord}_t(c_{\underline{\alpha},\underline{r}}x_1^{\alpha_1}x_2^{\alpha_2}z^{r_1}\cdots f_i^{r_{i+1}}\circ\gamma) & >\langle\nu,(\alpha_1,\alpha_2)+r_1\gamma_1+\cdots+r_{i+1}\gamma_{i+1}\rangle\geq n_{i+1}\langle\nu,\gamma_{i+1}\rangle\\
\end{array}\]
Therefore ord$_t(f_{i+1}\circ\gamma)=n_{i+1}\langle\nu,\gamma_{i+1}\rangle$. By Lemma \ref{LemExpSR}
\[f=f_{i+1}^{e_{i+1}}+\sum c_{\underline{\alpha},\underline{r}}x_1^{\alpha_1}x_2^{\alpha_2}z^{r_1}\cdots f_i^{r_{i+1}}\]
where $(\alpha_1,\alpha_2)+r_1\gamma_1+\cdots+r_{i+1}\gamma_{i+1}>n_{i+1}e_{i+1}\gamma_{i+1}$. Then
\[\begin{array}{rl}
\mbox{ord}_t(f_{i+1}^{e_{i+1}}\circ\gamma) & =l_{i+1}(\nu)\\
\\
\mbox{ord}_t(c_{\underline{\alpha},\underline{r}}x_1^{\alpha_1}x_2^{\alpha_2}z^{r_1}\cdots f_i^{r_{i+1}}\circ\gamma) & >\langle\nu,(\alpha_1,\alpha_2)+r_1\gamma_1+\cdots+r_{i+1}\gamma_{i+1}\rangle\geq l_{i+1}(\nu)\\
\end{array}\]
Then ord$_t(f\circ\gamma)=l_{i+1}(\nu)=l_{j+1}(\nu)\leq m$, which is a contradiction. Therefore ord$_t(f_i\circ\gamma)=\langle\nu,\gamma_{i+1}\rangle$.

We can prove that ord$_t(f_k\circ\gamma)=\langle\nu,\gamma_{k+1}\rangle$ for $i<k\leq j$ one after the other exactly as the proof of ord$_t(f_i\circ\gamma)=\langle\nu,\gamma_{i+1}\rangle$.

We prove now that ord$_t(f_{j+1}\circ\gamma)>\frac{m}{e_{j+1}}$. Suppose the contrary, ord$_t(f_{j+1}\circ\gamma)\leq\frac{m}{e_{j+1}}$.

If $j+2=g$ then by Lemma \ref{LemExpSR}
\[f=f_{j+2}=f_{j+1}^{n_{j+2}}-c_{j+2}x_1^{\alpha_1^{(g)}}x_2^{(g)}z^{r_1^{(g)}}\cdots f_j^{r_{j+1}^{(g)}}+\sum c_{\underline{\alpha},\underline{r}} x_1^{\alpha_1}x_2^{\alpha_2}z^{r_1}\cdots f_{j+1}^{r_{j+2}}\]
while, if $j+2<g$ we have the expansion
\[f=f_{j+1}^{e_{j+1}}+\sum c_{\underline{\alpha},\underline{r}} x_1^{\alpha_1}x_2^{\alpha_2}z^{r_1}\cdots f_{j+1}^{r_{j+2}}\]
with $(\alpha_1,\alpha_2)+r_1\gamma_1+\cdots+r_{j+2}\gamma_{j+2}>n_{j+2}e_{j+2}\gamma_{j+2}$. In both cases we have
\[f=f_{j+1}^{e_{j+1}}+\sum c_{\underline{\alpha},\underline{r}}x_1^{\alpha_1}x_2^{\alpha_2}z^{r_1}\cdots f_{j+1}^{r_{j+2}}\]
with $(\alpha_1,\alpha_2)+r_1\gamma_1+\cdots+r_{j+2}\gamma_{j+2}\geq n_{j+2}e_{j+2}\gamma_{j+2}$. Looking at the expansion, since ord$_t(f_{j+1}^{e_{j+1}}\circ\gamma)\leq m$, there must exist $c_{\underline{\alpha},\underline{r}}\neq 0$ such that
ord$_t(f_{j+1}^{e_{j+1}}\circ\gamma)=\langle\nu,(\alpha_1,\alpha_2)+r_1\gamma_1+\cdots+r_{j+1}\gamma_{j+1}\rangle+r_{j+2}\mbox{ord}_t(f_{j+1}\circ\gamma)$, or equivalently
\[(e_{j+1}-r_{j+2})\mbox{ord}_t(f_{j+1}\circ\gamma)=\langle\nu,(\alpha_1,\alpha_2)+r_1\gamma_1+\cdots+r_{j+1}\gamma_{j+1}\rangle\geq n_{j+2}e_{j+2}\langle\nu,\gamma_{j+2}\rangle-r_{j+2}\langle\nu,\gamma_{j+2}\rangle\]
And since $r_{j+2}<e_{j+1}$ we deduce that ord$_t(f_{j+1}\circ\gamma)\geq\langle\nu,\gamma_{j+2}\rangle>\frac{m}{e_{j+1}}$, which is a contradiction.

Finally we prove that ord$_t(f_k^{e_k}\circ\gamma)>m$. By Lemma \ref{LemExpSR}, $f=f_k^{e_k}+\sum c_{i_1i_2r}x_1^{i_1}x_2^{i_2}z^r$
with $(i_1,i_2)+r\gamma_1>n_ke_k\gamma_k$. Then ord$_t(x_1^{i_1}x_2^{i_2}z^r\circ\gamma)=\langle\nu,(i_1,i_2)+r\gamma_1\rangle\geq l_k(\nu)>m$, and the result follows.
\hfill$\Box$

\vspace{3mm}

\begin{Lem} For $m\in\Z_{>0}$ and $\nu\in L_m$, we have the following.
\begin{enumerate}
\item[(i)] If $i\leq j(m,\nu)$, then
\[F_{i-1,\nu}^{(\frac{l_i(\nu)}{e_{i-1}})}=0\mbox{ if and only if }{x_1^{(\nu_1)}}^{\alpha_1^{(i)}}{x_2^{(\nu_2)}}^{\alpha_2^{(i)}}{z^{(\langle\nu,\gamma_1\rangle)}}^{r_1^{(i)}}\cdots {F_{i-2,\nu}^{(\frac{l_{i-1}(\nu)}{e_{i-2}})}}^{r_{i-1}^{(i)}}=0\]
Roughly speaking the part $G_{i,\nu}$ in equation (\ref{EqP}) is not meaningful.

\

\item[(ii)] For $m(\nu)\leq j\leq g_1$ we have
\[V(I^\nu,F_{i,\nu}^{(\frac{l_i(\nu)}{e_i})})_{m(\nu)\leq i\leq j}\cap D(x_1^{(\nu_1)})\subset D(z^{(\langle\nu,\gamma_1\rangle)})\cap
D(F_{1,\nu}^{(\frac{l_2(\nu)}{e_1})})\cap\cdots\cap D(F_{j-1,\nu}^{(\frac{l_j(\nu)}{e_{j-1}})}),\]
and for $g_1<j\leq g$ we have
\[V(I^\nu,F_{i,\nu}^{(\frac{l_i(\nu)}{e_i})})_{m(\nu)\leq i\leq j}\cap D(x_1^{(\nu_1)})\cap D(x_2^{(\nu_2)})\subset D(z^{(\langle\nu,\gamma_1
\rangle)})\cap D(F_{1,\nu}^{(\frac{l_2(\nu)}{e_1})})\cap\cdots\cap D(F_{j-1,\nu}^{(\frac{l_j(\nu)}{e_{j-1}})}).\]
\end{enumerate}
\label{TechLem}
\end{Lem}

{\em Proof.}

(i) First observe that if $j(m,\nu)\geq i$, then $F_{i,\nu}^{(\frac{l_i(\nu)}{e_i})}\in J_m^\nu$ and it has the form given in equation (\ref{EqP}) in Corollary \ref{Corolario3}. The proof is obvious if $G_{i,\nu}=0$. Suppose  the contrary. Then $\nu\in\r_1\cup\r_2$, and the claim is not obvious if $\bar G_{i,\nu}\neq 0$, where
\[\bar G_{i,\nu}=\sum c_{\underline{\alpha},\underline{r}}{x_1^{(\nu_1)}}^{\alpha_1}{x_2^{(\nu_2)}}^{\alpha_2}{z^{(\langle\nu,\gamma_1\rangle)}}^{r_1}\cdots
{F_{i-2,\nu}^{(\frac{l_{i-1}(\nu)}{e_{i-2}})}}^{r_{i-1}}\]
with the conditions $c_{\underline{\alpha},\underline{r}}\neq 0$, $\langle\nu,(\alpha_1,\alpha_2)+r_1\gamma_1+\cdots+r_i\gamma_i\rangle=n_i\langle\nu,\gamma_i\rangle$ and $r_i=0$. Then, $F_{i-1,\nu}^{(\frac{l_i(\nu)}{e_{i-1}})}=0$ if and only if
\[c_i{x_1^{(\nu_1)}}^{\alpha_1^{(i)}}{x_2^{(\nu_2)}}^{\alpha_2^{(i)}}{z^{(\langle\nu,\gamma_1\rangle)}}^{r_1^{(i)}}\cdots {F_{i-2,\nu}^{(\frac{l_{i-1}(\nu)}{e_{i-2}})}}^{r_{i-1}^{(i)}}+\bar G_{i,\nu}=0\]
where, remember that any term $c_{\underline{\alpha},\underline{r}}{x_1^{(\nu_1)}}^{\alpha_1}{x_2^{(\nu_2)}}^{\alpha_2}{z^{(\langle\nu,\gamma_1\rangle)}}^{r_1}\cdots
{F_{i-2,\nu}^{(\frac{l_{i-1}(\nu)}{e_{i-2}})}}^{r_{i-1}}$ in particular appears in $F_{i,\nu}^{(\frac{l_i(\nu)}{e_i})}$, and hence it satisfies that $(\alpha_1^{(i)},\alpha_2^{(i)})+r_1\gamma_1+\cdots+r_{i-1}^{(i)}\gamma_{i-1}<(\alpha_1,\alpha_2)+r_1\gamma_1+\cdots+r_{i-1}\gamma_{i-1}$. Then we can write last equation as
\[{x_1^{(\nu_1)}}^{\alpha_1^{(i)}}{x_2^{(\nu_2)}}^{\alpha_2^{(i)}}{z^{(\langle\nu,\gamma_1\rangle)}}^{r_1^{(i)}}\cdots {F_{i-2,\nu}^{(\frac{l_{i-1}(\nu)}{e_{i-2}})}}^{r_{i-1}^{(i)}}\left(c_i+P(x_1^{(\nu_1)},x_2^{(\nu_2)},z^{(\langle\nu,\gamma_1\rangle)},\ldots,F_{i-2,\nu}^
{(\frac{l_{i-1}(\nu)}{e_{i-2}})})\right)\]
where $P$ is a polynomial non-unit. If we work it all out, then we can write $P$ as a polynomial in $x_1^{(\nu_1)}$, $x_2^{(\nu_2)}$ and $z^{(\langle\nu,\gamma_1\rangle)}$. Now we use that $\nu\in\r_1\cup\r_2$, and hence either $\nu_1=0$ or $\nu_2=0$. Then $P(x_1^{(\nu_1)},x_2^{(\nu_2)},z^{(\langle\nu,\gamma_1\rangle)})$ depends on the origin $\gamma(0)$ of the jets, and, since we are dealing with germs of quasi-ordinary singularities we can always suppose that $|x_i^{(\nu_i)}|<<1$ for $i=1$ or $2$ depending on whether $\nu_1=0$ or $\nu_2=0$. As a consequence we can always suppose that $P<<c_i$.

\vspace{2mm}

(ii) The inclusions follow directly from Corollary \ref{Corolario3} when $G_{j,\nu}=0$. If $G_{j,\nu}\neq 0$, the proof is by induction on $j$. For $j=m(\nu)$ the claim says that
 \[V(I^\nu,F_{m(\nu),\nu}^{(\frac{l_{m(\nu)}(\nu)}{e_{m(\nu)}})})\cap D(M)\subset D(z^{(\langle\nu,\gamma_1\rangle)}),\]
 where
 \[M=\left\{\begin{array}{cl}
 x_1^{(\nu_1)} & \mbox{ if }g_1\geq 1\\
 \\
 x_1^{(\nu_1)}x_2^{(\nu_2)} & \mbox{ if }g_1=0\\
 \end{array}\right.\]
 and it follows by Lemma \ref{Lemfk}. Suppose it is true for $j-1$ and we prove it for $j$. We only have to prove that $F_{j-1,\nu}^{(\frac{l_j(\nu)}{e_{j-1}})}\neq 0$. But this follows by (i).

 \hfill$\Box$

     \bibliographystyle{amsplain}
\def\cprime{$'$}
\providecommand{\bysame}{\leavevmode\hbox to3em{\hrulefill}\thinspace}
\providecommand{\MR}{\relax\ifhmode\unskip\space\fi MR }
\providecommand{\MRhref}[2]{%
  \href{http://www.ams.org/mathscinet-getitem?mr=#1}{#2}
}
\providecommand{\href}[2]{#2}

\end{document}